\input lanlmac
\def\href#1#2{{#2}}

\input epsf.tex

\overfullrule=0mm

\newcount\figno
\figno=0
\def\fig#1#2#3{
\par\begingroup\parindent=0pt\leftskip=1cm\rightskip=1cm\parindent=0pt
\baselineskip=11pt
\global\advance\figno by 1
\midinsert
\epsfxsize=#3
\centerline{\epsfbox{#2}}
\vskip 12pt
{\bf Fig.\ \the\figno:} #1\par
\endinsert\endgroup\par
}
\def\figlabel#1{\xdef#1{\the\figno}}
\def\encadremath#1{\vbox{\hrule\hbox{\vrule\kern8pt\vbox{\kern8pt
\hbox{$\displaystyle #1$}\kern8pt}
\kern8pt\vrule}\hrule}}


\magnification=\magstep1
\baselineskip=12pt
\hsize=6.3truein
\vsize=8.7truein
 at 8truept
 at 8truept
 at 10truept

\vbox{\hfill IPhT-t13/039}

\bigskip\bigskip

\font\bigrm=cmr12 at 14pt \centerline{\bigrm On irreducible maps and slices} 
 
\bigskip\bigskip

\centerline{J. Bouttier$^{1,2}$ and E. Guitter$^1$}
  \smallskip
\centerline{$^1$ Institut de Physique Th\'eorique}
  \centerline{CEA, IPhT, F-91191 Gif-sur-Yvette, France}
  \centerline{CNRS, URA 2306}
\centerline{$^2$ D\'epartement de Math\'ematiques et Applications}
\centerline{\'Ecole normale sup\'erieure, 45 rue d'Ulm, F-75231 Paris Cedex 05}  
\centerline{\tt jeremie.bouttier@cea.fr}
\centerline{\tt emmanuel.guitter@cea.fr}

  \bigskip

     \bigskip\bigskip

     \centerline{\bf Abstract}
     \smallskip
     {\narrower\noindent

We consider the problem of enumerating $d$-irreducible maps, i.e.\ planar maps whose all
cycles have length at least $d$, and such that any cycle of length $d$ is the boundary
of a face of degree $d$. We develop two approaches in parallel: the natural approach 
via substitution where these maps are obtained from general maps by a replacement
of all $d$-cycles by elementary faces, and a bijective approach via slice decomposition
which consists in cutting the maps along shortest paths. Both lead to explicit expressions 
for the generating functions of $d$-irreducible maps with controlled face degrees,
summarized in some elegant ``pointing formula''. We provide an equivalent description 
of $d$-irreducible slices in terms of so-called  $d$-oriented trees. We finally show that irreducible maps 
give rise to a hierarchy of 
discrete integrable equations which include equations encountered previously in the 
context of naturally embedded trees.   

\par}

     \bigskip

\nref\TutteCPT{W.T. Tutte, {\it A census of planar triangulations}, Canad.\ J.\ Math.\ 
{\bf 14} (1962) 21-38.}
\nref\TutteCHP{W.T. Tutte, {\it A Census of Hamiltonian polygons}, Canad.\ J.\ Math. {\bf 14} (1962) 402-417.}
\nref\TutteCS{W.T. Tutte, {\it A Census of slicings},
Canad.\ J.\ Math. {\bf 14} (1962) 708-722.}
\nref\TutteCPM{W.T. Tutte, {\it A Census of planar maps}, Canad.\ J.\ Math.\ 
{\bf 15} (1963) 249-271.}
\nref\OHRMT{{\it The Oxford Handbook of Random Matrix Theory}, Oxford University Press (2011)
G. Akemann, J. Baik and P. Di Francesco Editors.}
\nref\Schae{G. Schaeffer, {\it Conjugaison d'arbres
et cartes combinatoires al\'eatoires}, PhD Thesis, Universit\'e 
Bordeaux I (1998).}
\nref\BFa{O. Bernardi and \'E. Fusy, {\it A bijection for triangulations, quadrangulations, pentagulations, etc.},
J.\ Combin.\ Theory Ser.\ A {\bf 119} (2012) 218-244, arXiv:1007.1292 [math.CO].}
\nref\BFb{O. Bernardi and \'E. Fusy, {\it Unified bijections for maps with prescribed degrees and girth},
J.\ Combin.\ Theory Ser.\ A {\bf 119} (2012) 1351-1387, arXiv:1102.3619 [math.CO].}
\nref\AlPo{M. Albenque and D. Poulalhon, {\it Generic method for bijections between blossoming trees and planar maps}, arXiv:1305.1312 [math.CO].}
\nref\PMCF{J. Bouttier and E. Guitter, {\it Planar maps and continued fractions}, Commun.\ Math.\ Phys.\ {\bf 309}:3 
(2012) 623-662, arXiv:1007.0419 [math.CO].}
\nref\LG{J.-F. Le Gall {\it Uniqueness and universality of the Brownian map}, Ann.\ Probab.\ to appear,
arXiv:1105.4842 [math.PR].}
\nref\MuSch{R. Mullin and P. Schellenberg, {\it The enumeration of c-nets via quadrangulations}, 
J.\ Combin.\ Theory {\bf 4} (1968) 259-276.}
\nref\FusyT{\'E. Fusy, {\it Transversal structures on triangulations: A combinatorial study and straight-line drawings},
Discrete Mathematics {\bf 309} (2009) 1870-1894, arXiv:math/ 0602163 [math.CO].}
\nref\FPS{\'E. Fusy, D. Poulalhon and G. Schaeffer, {\it Dissections, orientations, and trees, with applications 
to optimal mesh encoding and to random sampling}, Trans. Algorithms {\bf 4}(2) (2008) Art 19,
arXiv:0810.2608 [math.CO].}
\nref\FLAJFRAC{P. Flajolet, {\it Combinatorial aspects of continued
fractions}, Discrete Mathematics 32 (1980), 125-161.
Reprinted in the 35th Special Anniversary Issue of Discrete
Mathematics, Vol.\ {\bf 306}, Issue 10-11, Pages 992-1021 (2006).}
\nref\BaFlScSo{C. Banderier, P. Flajolet, G. Schaeffer and M. Soria, {\it 
Random Maps, Coalescing Saddles, Singularity Analysis, and Airy Phenomena},
Random Structures and Algorithms {\bf 19} (2001) 194-246.}
\nref\GEOD{J. Bouttier, P. Di Francesco and E. Guitter, {\it Geodesic
distance in planar graphs}, Nucl.\ Phys.\ {\bf B663}[FS] (2003) 535-567, 
arXiv:cond-mat/0303272.}
\nref\DFG{P. Di Francesco and E. Guitter, {\it Integrability of graph 
combinatorics via random walks and heaps of dimers}, 
J.\ Stat.\ Mech.\ (2005) P09001, arXiv:math/0506542 [math.CO].}
\nref\PDFRaman{P. Di Francesco, {\it Geodesic Distance in Planar Graphs: 
An Integrable Approach}, The Ramanujan Journal {\bf 10} (2005) 153-186,
arXiv:math/0506543 [math.CO].}
\nref\Kuba{M. Kuba, {\it A note on naturally embedded ternary trees}, Electr.\ J.\ Comb.\ {\bf 18}(1) (2011) P142,
arXiv:0902.2646 [math.CO].}
\nref\BM{M. Bousquet-M\'elou, {\it Limit laws for embedded trees. Applications to the integrated superBrownian 
excursion}, Random Structures and Algorithms {\bf 29}(4) (2006) 475-523, arXiv:math/0501266 [math.CO].}
\nref\AB{M. Albenque and J. Bouttier, {\it Constellations and multicontinued fractions: application to Eulerian triangulations},
DMTCS Proceedings AR 805-816 (2012) 24th International Conference on Formal Power Series and Algebraic Combinatorics (FPSAC 2012), arXiv:1112.6379 [math.CO].}
\nref\PSEUDOQUAD{J. Bouttier and E. Guitter, {\it Distance statistics
in quadrangulations with a boundary, or with a self-avoiding loop},
J. Phys.\ A: Math.\ Theor.\ {\bf 42} (2009) 465208, arXiv:0906.4892
[math-ph].}
\nref\Brown{W.G. Brown, {\it Enumeration of quadrangular dissections of the disk}, Canad.\ J.\ Math.\ 
{\bf 17} (1965) 302-317.}
\nref\FS{P. Flajolet and R. Sedgewick, {\it Analytic Combinatorics}, Cambridge 
University Press (2009).}
\nref\Comtet{L. Comtet, {\it Advanced Combinatorics:
The Art of Finite and Infinite Expansions}, Springer, 1974.}
\nref\GouldJack{I.P. Goulden and D.M.
Jackson, {\it Combinatorial Enumeration}, John Wiley \& Sons, New
York, 1983, republished by Dover, New York, 2004.}
\nref\DSQ{J. Bouttier and E. Guitter, {\it Distance statistics in quadrangulations with no multiple edges 
and the geometry of minbus}, J.\ Phys.\ A: Math.\ Theor.\ {\bf 43} (2010) 205207,
arXiv:1002.2552 [math-ph].}
\nref\GJKP{I.P. Goulden and D.M.
Jackson, {\it The KP hierarchy, branched covers, and triangulations}, 
Adv.\ Math\ {\bf 219}(3) (2008) 932-951, arXiv:0803.3980 [math.CO].}
\nref\AlFuPo{M. Albenque, \'E. Fusy and D. Poulalhon, {\it On symmetric quadrangulations and triangulations}, European Journal of Combinatorics {\bf 35} (2014) 13-31, arXiv:1103.3657 [math.CO].}
\nref\FusyNB{\'E. Fusy, {\it New bijective links on planar maps via orientations}, 
European Journal of Combinatorics {\bf 31}(1) (2010) 145-160, arXiv:0810.2607 [math.CO].}
\nref\Viennot{X.G. Viennot, {\it Une th\'eorie combinatoire des polyn\^omes 
orthogonaux}, Lecture Notes UQAM, 217p., Publication du LACIM, Universit\'e 
du Qu\'ebec \`a Montr\'eal (1984), re\'ed. 1991.}
\nref\CF{G. Collet and \'E. Fusy, {\it A simple formula for the series of bipartite and quasi-bipartite maps 
with boundaries}, FPSAC 2012, Nagoya, Japan, DMTCS proc. {\bf AR} (2012) 607-618, arXiv:1205.5215 [math.CO].}
\nref\WIP{J. Bouttier and E. Guitter, {\it A note on irreducible maps with several boundaries}, arXiv:1305.4816 [math.CO].}

\newsec{Introduction}

\subsec{General introduction}

The enumeration of planar maps has been of constantly renewed interest in combinatorics since Tutte's seminal 
papers [\xref\TutteCPT-\xref\TutteCPM], some of its developments arising from theoretical physics or probability theory. Among the various enumeration
techniques used so far, let us mention Tutte's original recursive decomposition, the matrix integral approach 
\OHRMT\ and 
the more recent use of bijections with trees \Schae. While the first two approaches are fairly systematic (i.e.\ 
allow to translate almost automatically the counting problem into equations), a weakness of the 
third one is that it relies on some preliminary divination of the family of trees with which a bijection is to be found.
In practice, one usually solves the counting problem via another approach and, guided by the knowledge
of the result, discovers the bijection afterwards. For this reason, the numerous bijections found
in the literature might appear as a myriad of rather {\it ad hoc} tools. Recently, some authors have undertaken 
the task of understanding the general underlying principles of the bijective approach and providing a unified framework,
the ``master bijection'', in which all the previously known bijections appear as special cases [\xref\BFa-\xref\AlPo].      

On some other side, it was realized that one of the desirable feature of trees, namely that they are easy
to enumerate thanks to their natural recursive structure, could be directly achieved at the level of the
maps themselves via the so-called {\it slice decomposition} \PMCF. This approach, close in spirit to Tutte's approach,
has the merit of relying on a simple systematic construction, which consists in cutting a map along some shortest paths.
The parts obtained in this decomposition are maps with geodesic boundaries, called slices for short, and may 
themselves be iteratively cut into smaller slices, reproducing a recursive tree-like structure. We observe that slices
were also used in \LG\ under the name of DMGB (discrete maps with geodesic boundaries).

In this paper, we extend the slice decomposition formalism to the case of maps with a girth constraint, i.e. 
with a control on the minimal length of their cycles. Such maps were already considered in the 
master bijection framework and we shall indeed recover some of the results of [\xref\BFa,\xref\BFb].  We actually consider the slightly
more general case of {\it irreducible maps}, i.e.\ maps with a girth constraint and without separating shortest cycles
(as we shall see below, it is indeed more general since the former may be recovered by setting some parameter to $0$). 
Irreducible triangulations and quadrangulations were first enumerated respectively by Tutte \TutteCPT, and by 
Mullin and Schellenberg  \MuSch\ using a substitution approach (these authors used the denomination ``simple'' instead of irreducible, which is slightly 
misleading since a simple map is nowadays understood as a map without loops or multiple edges). These results
were later promoted to bijections with trees [\xref\FusyT,\xref\FPS]. 

We actually start by extending the Tutte-Mullin-Schellenberg 
substitution approach to arbitrary irreducible maps. As a key ingredient, we use
an expression for the generating function of maps with a boundary which originates combinatorially from the slice 
decomposition. It turns out that slice decomposition somehow ``commutes'' with the substitution approach.
We are therefore led to studying irreducible slices and their recursive decomposition.

Before presenting in more details the outline of the paper, we would
like to discuss the influence of Philippe Flajolet on this
work. Making an exhaustive list is an impossible task so let us focus
on two particular examples. First, the slice decomposition \PMCF\ was
inspired by the combinatorial theory of continued fractions
\FLAJFRAC. This theory is far more general than the context of planar
maps, but it proved to be the key to understanding the phenomenon of
``discrete integrability'', on which we will return below. Second,
Philippe was no stranger to counting maps via substitution: in
\BaFlScSo, he and his coauthors have shown that a universal ``Airy
phenomenon'' occurs when, generally speaking, one decomposes a map
into (multi)connected ``cores''. Figuring out whether the same
phenomenon subsists in the present context of irreducible maps is an
intriguing open question.

\subsec{Definitions}

A {\it planar map} is an embedding of a connected graph in the sphere
without edge crossing, considered up to continuous deformation.
It is made of {\it vertices}, {\it edges} and {\it faces}. A {\it rooted} map is a map with a distinguished oriented edge,
the {\it root edge}. The face on the right of the root edge is called the {\it outer face} (whose degree is called the
{\it outer degree} of the map), the other ones being referred
to as {\it inner faces}. For $n$ and $d$ two positive integers, a {\it $d$-angular dissection of the $n$-gon} is a
rooted map whose outer degree is $n$ and where all inner faces have degree $d$. 

The {\it girth} of a map is the minimal length (number of edges) of its cycles, a cycle being a simple closed
path on the map. Note that, by this definition, trees have an infinite girth since they contain no cycle. In 
a map not reduced to a tree, the degree of every face is larger than or equal to the girth (which is finite).

Given a nonnegative integer $d$, we say that a rooted map is {\it $d$-irreducible} if its girth is at least
$d$ and any cycle of length $d$ is the boundary of an inner face of degree $d$. 
Note that every rooted map
is $0$-irreducible. Furthermore, by definition, a $d$-irreducible map with outer degree smaller 
than or equal to $d$ is either a tree (with one face of even degree $\leq d$) or is reduced to a cycle of length $d$ (delimiting
two faces of degree $d$). A $d$-irreducible $d$-angular dissection will be called an irreducible $d$-angular dissection for short.
Note finally that maps of girth at least $d$ are nothing but $(d-1)$-irreducible maps containing no $(d-1)$-valent face,
and in this sense, the class of irreducible maps is more general than that of maps with controlled girth.

Let us denote by $F_n^{(d)}(z;x_{d+1},x_{d+2},\ldots )$ the generating function of $d$-irreducible maps with 
outer degree $n$, counted with a weight $z$ per inner face of degree $d$ and, for all $i \geq d+1$, a weight
$x_i$ per inner face of degree $i$. Our motivation for choosing a different notation for the weight of faces of degree $d$
is that it plays a very different role in the forthcoming expressions. From the above remark, we have, for $n \leq d$
\eqn\Fsmalln{F_n^{(d)}(z;x_{d+1},x_{d+2},\ldots )=\left\{\matrix{ {\rm Cat}(n/2) &   \hbox{for $n < d$} \cr
 {\rm Cat}(d/2) + z &   \hbox{for $n = d$} }\right.}
 where ${\rm Cat}(k)$ is equal to ${2k \choose k}/(k+1)$ (the $k$-th Catalan number, counting rooted trees 
 with $k$ edges, hence an outer degree $2k$) for integer $k$
 and $0$ for noninteger $k$.

Let us now discuss a few interesting specializations of
$F_n^{(d)}$. First, by taking all $x_i$, $i\geq d+1$, to $0$, we
obtain the generating function $f_n^{(d)}(z)=F_n^{(d)}(z;0,0,\ldots )$
of irreducible $d$-angular dissections of the $n$-gon, depending on
the single variable $z$ coupled to the number of inner faces. Second,
by conversely taking $z=0$, we forbid all faces of degree $d$ hence
all cycles of length $d$, so that
$F_n^{(d)}(0;x_{d+1},x_{d+2},\ldots)$ coincides with the generating
function of rooted maps of girth at least $d+1$ and outer degree $n$,
as studied in \BFb. Note that, in the particular case $d=0$,
$F_n(x_1,x_2,\ldots)=F_n^{(0)}(0;x_{1},x_{2},\ldots)$ is nothing but
the generating function of arbitrary maps with outer degree $n$.
Finally, a third specialization concerns {\it bipartite planar maps},
i.e.\ maps whose all faces have even degrees: it is obtained by taking
$n$ and $d$ even, and setting all odd $x_i$ to $0$.
  
\subsec{Overview of the main results}

Our main result is a general expression for $F_n^{(d)}$ in terms of
two auxiliary quantities, which we denote by $R^{(d)}$ and $S^{(d)}$,
and which may be interpreted as $d$-irreducible slice generating
functions. As such, these quantities are themselves determined by an
explicit system of two equations, which is algebraic whenever we
impose a bound on the face degrees (i.e.\ $x_i=0$ for $i$ large
enough). A particularly elegant expression for $F_n^{(d)}$ is via a
{\it pointing formula} which amounts to counting {\it annular maps},
i.e.\ rooted maps having a distinguished inner face of degree $d$.
This pointing formula has a clear combinatorial interpretation from
the slice decomposition.  We now mention a number of other interesting
results appearing on the way.

We find that $R^{(d)}$ and $S^{(d)}$ are particular members
(corresponding essentially to the first two values $k=-1,0$) of a
larger family of generating functions $V_k^{(d)}$ ($k\geq -1$) of
so-called $d$-irreducible $k$-slices, where $k$ controls some excess
boundary length of the slice (a more precise definition will come in
due time). We provide a closed system of equations for $V_k^{(d)}$
which results from an elementary recursive decomposition of
$k$-slices.  We are then able to explicitly eliminate all $V_k^{(d)}$
with $k\geq 1$, yielding the wanted system of two equations
determining $R^{(d)}$ and $S^{(d)}$.

While all the enumeration is carried out in terms of slices, we also
discuss, in the case of irreducible $d$-angular dissections, an
equivalent formulation in terms of trees: more precisely, $V_k^{(d)}$
may in this case be interpreted as the generating function of
so-called $d$-oriented $k$-trees, reminiscent of the $d/(d-2)$-trees
considered in [\xref\BFa,\xref\BFb]. In the particular cases $d=3$ and $d=4$,
$d$-oriented trees reduce respectively to ternary and binary trees and
we recover the bijections of [\xref\FusyT,\xref\FPS].

Finally, we consider slices subject to a control of an extra
parameter, namely their maximal length.  In the case of
non necessarily irreducible maps, generating functions of these
objects are known to be solution of a hierarchy of ``discrete
integrable equations'' [\xref\GEOD-\xref\PDFRaman,\xref\PMCF]. We show that this phenomenon subsists in
the $d$-irreducible case and explicit the corresponding equations. A
particular attention is paid to the cases $d=3$ and $d=4$ for which we
recover integrable equations describing so-called naturally embedded
trees, respectively in their ternary \Kuba\ and binary \BM\ flavour.
In particular, answering a question raised by Bousquet-M\'elou, we
provide a combinatorial explanation of Proposition 25 in \BM, in the
same spirit as those previously found for well labeled trees \PMCF\ and very well labeled 
trees \AB.

\subsec{Outline of the paper}

We now come to the detailed plan of the paper. We begin by describing
the substitution approach to the enumeration of $d$-irreducible
maps. For pedagogical reasons, we first address the simplest cases
$d=4$ (Sect.~2.1) and $d=3$ (Sect.~2.2). The general case is
treated in Sect.~3.  We proceed by induction on $d$ and explain in
Sect.~3.1 how $F_n^{(d)}$ is related to $F_n^{(d-1)}$. We deduce in
Sect.~3.2 a relation between $F_n^{(d)}$ and the generating function
$F_n$ of arbitrary maps. We then exploit known expressions for $F_n$
to obtain $F_n^{(d)}$, first in the simpler bipartite case (Sect.~3.3), 
then in the general case (Sect.~3.4). This yields the first
derivation of our main result, together with the pointing formula.  The
quantities $R^{(d)}$, $S^{(d)}$ and $V_k^{(d)}$ appear in this
derivation as mere intermediate products. In Sect.~4, we take the
time to discuss their combinatorial significance as slice generating
functions. We recall in Sect.~4.1 the definition of slices and
extend it to what we call $k$-slices. We then show in Sect.~4.2 that
$V_k^{(d)}$ is nothing but the generating function of $d$-irreducible
$k$-slices. In Sect.~5, we obtain a recursive decomposition of these
slices, actually in two variants (Sects.~5.1 and 5.2).  The closed
system for the $V_k^{(d)}$ is deduced in Sect.~5.3 and we discuss in
Sects.~5.4 (bipartite case) and 5.5 (general case) how to eliminate
the $V_k^{(d)}$ with $k\geq 1$.  This yields a second route to the
system of equations satisfied by $R^{(d)}$ and $S^{(d)}$. Sect.~6
discusses the equivalent formulation in terms of trees (for the case
of irreducible $d$-angular dissections): we define $d$-oriented
$k$-trees in Sect.~6.1 and exhibit their one-to-one correspondence
with slices. An alternative description of the correspondence as a
closure algorithm is given in Sect.~6.2.  Some simplifications
occuring in the bipartite case are mentioned in Sect.~6.3. Sect.~7
is devoted to the bijective proof of the pointing formula: we explain
in Sect.~7.1 how to build an annular map out of slices, and present
in Sect.~7.2 the inverse mapping (involving the notion of
``lift''). Combined with the results of Sect.~5, this provides a
second derivation of our main result. Sect.~8 is devoted to discrete
integrable equations: we first discuss the particular cases of
irreducible quadrangular (Sect.~8.1) and triangular (Sect.~8.2)
slices, related to naturally embedded trees, and we then write down
the general equations in Sect.~8.3. Sect.~9 discusses other aspects
of irreducibility: in Sect.~9.1 we relax the definition of $d$-irreducibility 
for maps with outer degree $d$ and solve the corresponding enumeration
problem. Sect.~9.2 deals with the enumeration of $d$-irreducible maps
with two marked faces of degree strictly larger than $d$. Finally, we
consider in Sect.~9.3 a generalized notion of annular maps, leading to
beautiful identities extending the pointing formula.

\newsec{First simple cases}
\seclab\firstsimp
As a preamble to Sect.~3 where we shall explain in details the substitution approach
to $d$-irreducible maps, let us  discuss in the simplest case of quadrangular and triangular
dissections how this approach allows one to obtain expressions for generating functions of irreducible
maps at no cost. More precisely, substitution tells us that these generating functions may be 
obtained from those of arbitrary maps by a simple renormalization of the weights $x_i$. 
This turns out to be sufficient to determine them fully.

\fig{Left: an irreducible triangular dissection of the square. Right: an irreducible quadrangular
dissection of the hexagon.}{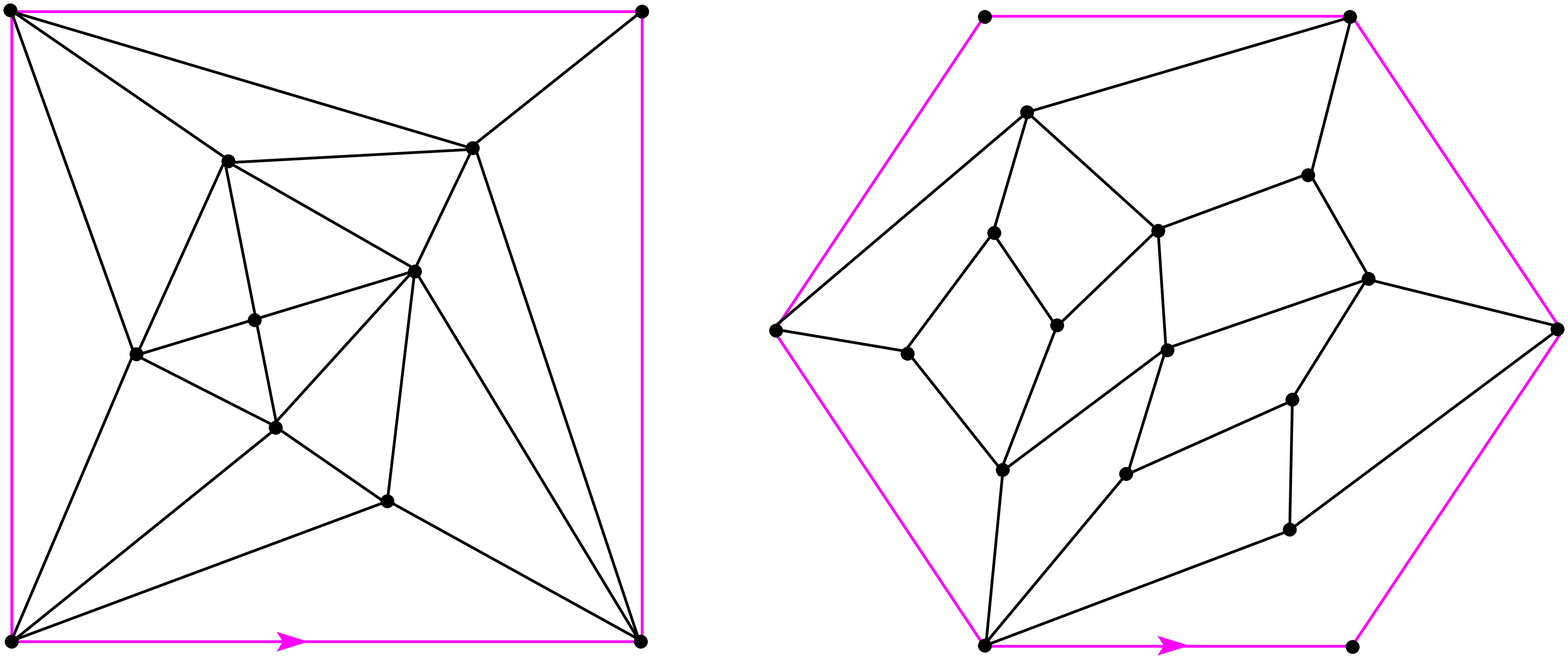}{10.cm}
\figlabel\triquaddiss

\subsec{Irreducible quadrangular dissections}
\subseclab\irredquad

We have at our disposal a number of expressions for the generating functions of general (non necessarily
irreducible) quadrangular dissections, counted with a weight $x_4$ per face. Recall that $F_n(x_1,x_2,\ldots)$
denotes the generating function of arbitrary planar maps with outer degree $n$, counted with weights $x_i$ per
inner face of degree $i$. In this section, we denote by $F_n=F_n(0,x_2,0,x_4,0,0,\ldots)$ its specialization to $x_i=0$ 
for $i$ different from $2$ or $4$, that is to say we consider maps with only bivalent or tetravalent inner faces. 
The expression for $F_n$ takes a particularly simple form if we introduce the series $R\equiv R(x_2,x_4)$, solution of 
 $R=1+x_2 R+3 x_4 R^2$  (which itself may be understood as a generating function of some kind). 
 We have in particular [\xref\DFG,\xref\PMCF]
\eqn\Fnquad{\eqalign{F_2&=R - \left(x_4 R^3\right),\cr
F_4&=2 R^2-3 R \left(x_4 R^3\right),\cr
F_6&= 5 R^3-9 R^2 \left(x_4 R^3\right).\cr
}}  

We now claim that there exists two formal power series $X_2 \equiv
X_2(z)$ and $X_4 \equiv X_4(z)$ such that the generating function
$f_n^{(4)}\equiv f_n^{(4)}(z)$ of {\it irreducible} quadrangular
dissections of the $n$-gon, counted with a weight $z$ per face, is
obtained from $F_n$ by the substitution $x_2 \to X_2(z)$ and
$x_4 \to X_4(z)$. Furthermore, the ``renormalized face weights'' $X_2$ and
$X_4$ are entirely determined by the condition \Fsmalln\ for
$d=4$. We will prove these statements in greater generality
in Section 3 below. At this stage, let us simply justify them
heuristically by noting that a general quadrangular dissection may be
obtained from an irreducible one by a renormalization procedure which,
so to say, consists in replacing each (four-valent) face of the latter
by a more general quadrangular dissection with outer degree
$4$. Actually, we also need to eliminate multiple edges which is
achieved by introducing bivalent faces and renormalizing them out.

Defining $r^{(4)}=R(X_2(z),X_4(z))$, we deduce from \Fnquad\ that
\eqn\Fnquadd{\eqalign{f_2^{(4)}&=1=r^{(4)} - \left(X_4(z) (r^{(4)})^3\right),\cr
f_4^{(4)}&=2+z=2 (r^{(4)})^2-3\, r^{(4)} \left(X_4(z) (r^{(4)})^3\right),\cr
f_6^{(4)}&= 5 (r^{(4)})^3-9 (r^{(4)})^2 \left(X_4(z) (r^{(4)})^3\right).\cr
}}
Now it is interesting to note that we do not need any precise expression for $X_2$ or $X_4$. Indeed, 
eliminating $X_4$ from the first two equations of \Fnquadd\ yields:
\eqn\algebRquad{z + (r^{(4)})^2 - 3 r^{(4)} + 2 = 0 ,}
which fully determines $r^{(4)}$ as a function of $z$. Once $r^{(4)}$ is known, we may get 
$f_6^{(4)}$ from the third equation which, after elimination of $X_4$, reads simply:
\eqn\fsixquad{f_6^{(4)}=9 (r^{(4)})^2-4 (r^{(4)})^3\ .}
Differentiating both equations \algebRquad\ and \fsixquad\ with respect to $z$, we deduce in particular
$(3-2 r^{(4)}){d r^{(4)}\over d z}=1$ and
\eqn\dfsixquaddz{{d f_6^{(4)}\over d z}=6 r^{(4)}(3-2 r^{(4)}){d r^{(4)}\over d z} = 6 r^{(4)}.}
This result is a particular example of a more general pointing formula which will be discussed in
details later.

Eq.~\algebRquad\ is more transparent upon setting 
\eqn\rtoT{r^{(4)}(z)=1+z \,T(z) ,}
as it then reads
\eqn\eqforT{T=1+z\, T^2.}
This allows to identify $T$ with the generating function of binary trees (with a weight $z$ per inner vertex). 
Using $[z^n]T={\rm Cat}(n)$, we arrive at
\eqn\fsixquadzn{[z^{n+2}] f_6^{(4)}={6\over n+2}\, {\rm Cat}(n) }
for the number of irreducible quadrangular dissections of the hexagon with $(n+2)$ squares, $n\geq 0$
(see Fig.~\triquaddiss\ for an example with $n=12$).
We recover here a result of [\xref\MuSch,\xref\FPS].

More general formulas are obtained along the same lines. From the expression
\eqn\Fnquadgen{F_{2m} = {\rm Cat}(m) R^m - {3 \over m-1} {2m \choose m-2} x_4 R^{m+2}}
for the generating function of arbitrary quadrangular dissections of
the $2m$-gon [\xref\DFG,\xref\PMCF], we readily deduce that the generating function of
irreducible ones reads
\eqn\fnquad{f_{2m}^{(4)}(z) = {2m \choose m-2} \left( {3 \over m-1} (r^{(4)})^{m-1} - {2 \over m} (r^{(4)})^{m}\right)}
and satisfies the pointing formula $d f_{2m}^{(4)}/dz = {2m \choose
m-2} (r^{(4)})^{m-2}$.  However the general coefficient of
$f_{2m}^{(4)}(z)$ does not seem to be ``nice''. In contrast, Mullin and
Schellenberg \MuSch\ obtained a nice general coefficient when
considering irreducible dissections whose outer boundary is simple. We
may recover their formula from the expression [\xref\PSEUDOQUAD, Eq.~(5.16)] 
\eqn\Fnquadsim{{\tilde F}_{2p} = {(3p-3)! \over p! (2p-1)!} \left( p x_4^{p-1} R^{3p-2} + (2-3p) x_4^p R^{3p} \right)}
for the generating function of arbitrary quadrangular dissections of
the $2p$-gon with a simple outer boundary. It is not difficult to check that this formula remains valid if we also
allow for bivalent faces provided $R$ is taken as the solution $R=1+x_2R+3x_4R^2$ as before. Substituting $x_4\to X_4(z)$ and
$R\to r^{(4)}$ in  \Fnquadsim\ and noting that $X_4(z) (r^{(4)})^3=z \, T(z)$ and $[z^k]\, T(z)^p= {p
\over 2k+p} {2k+p \choose k}$, we obtain after some algebra that the
number of irreducible quadrangular dissections of the $2m$-gon with $k$
inner faces and a simple outer boundary reads
\eqn\Fnquadsimcoef{[z^k] {\tilde f}^{(4)}_{2p} = {(3p-3)! \over (p-3)!(2p-1)!} {(2k-p-1)! \over k! (k-p+1)!}}
which, by the reparametrization $m=p-2$ and $n=k-p+1$, coincides for $m>0$ 
with Mullin and Schellenberg's formula. Let us observe that, in contrast to the present approach,
these authors started directly from Brown's formula \Brown\ for the number of quadrangular dissections 
of the $2p$-gon which both are simple (i.e.\ have no multiple edges)
and have a simple outer boundary.

\subsec{Irreducible triangular dissections}
\subseclab\irredtri

We may now play the same game with triangular dissections. As before,
we have simple expressions for the generating functions of general
triangular dissections, with weight $x_3$ per face. In this section,
$F_n=F_n(x_1,x_2,x_3,0,0,\ldots)$ denotes the generating function of
planar maps with outer degree $n$ and all inner faces of degree at
most $3$ (as obtained by specializing $x_i=0$ for $i > 3$).
Introducing the series $R,S$ in the variables $x_1,x_2,x_3$
specified by the equations $R=1+x_2 R + 2 x_3 R S$, $S=x_1+x_2 S+
x_3(S^2+2 R)$, we have \PMCF
\eqn\Fntri{\eqalign{F_1&=S - \left(x_3 R^2\right),\cr
F_2&= S^2+R-2 S \left(x_3 R^2\right),\cr
F_3&= S^3+3 R S-(3 S^2+2 R) \left(x_3 R^2\right), \cr
F_4&= S^4+6 R S^2+2 R^2-(4 S^3+8 R S) \left(x_3 R^2\right).\cr
}}  

We now claim that there exists three formal power series $X_1 \equiv
X_1(z)$, $X_2 \equiv X_2(z)$, $X_3 \equiv X_3(z)$ such that the
generating function $f_n^{(3)} \equiv f_n^{(4)}(z)$ of {\it
irreducible} triangular dissections of the $n$-gon, counted with a
weight $z$ per inner face, is obtained from $F_n$ by the substitution
$x_1 \to X_1(z)$, $x_2 \to X_2(z)$ and $x_3 \to X_3(z)$. Again, the
renormalized face weights $X_1$, $X_2$ and $X_3$ are entirely
determined by \Fsmalln\ for $d=3$. Setting
$r^{(3)}(z)=R(X_1(z),X_2(z),X_3(z))$ and
$s^{(3)}(z)=S(X_1(z),X_2(z),X_3(z))$, we can now write
\eqn\Fntrid{\eqalign{f^{(3)}_1&=0=s^{(3)} - \left(X_3(z) (r^{(3)})^2\right),\cr
f^{(3)}_2&= 1=(s^{(3)})^2+r^{(3)}-2 s^{(3)} \left(X_3(z)(r^{(3)})^2\right),\cr
f^{(3)}_3&= z=(s^{(3)})^3+3 r^{(3)} s^{(3)}-(3 (s^{(3)})^2+2 r^{(3)}) \left(X_3(z) (r^{(3)})^2\right), \cr
f^{(3)}_4&= (s^{(3)})^4+6 r^{(3)} (s^{(3)})^2+2 (r^{(3)})^2-(4 (s^{(3)})^3+8 r^{(3)} s^{(3)}) \left(X_3(z) (r^{(3)})^2\right).\cr
}}  
Eliminating $X_3$, the first three equations lead to the following algebraic system determining $r^{(3)}$ and $s^{(3)}$
in terms of $z$:
\eqn\algebRStri{r^{(3)}=1+(s^{(3)})^2\ , \qquad z+(s^{(3)})^3-s^{(3)}=0\ .}
As for $f_4^{(3)}$, we deduce from the fourth equation
\eqn\ffourtri{f_4^{(3)}=2+2 (s^{(3)})^2-3(s^{(3)})^4}
and, upon differentiating with respect to $z$
\eqn\dffourtridz{{d f_4^{(3)}\over d z}= 4 s^{(3)} (1-3 (s^{(3)})^2) {d s^{(3)}\over d z} = 4\, s^{(3)}}
since, from \algebRStri, $(1-3 (s^{(3)})^2) {d s^{(3)}\over d z} =1$.

Eq.~\algebRStri\ for $s^{(3)}$ is more transparent upon setting 
\eqn\stoT{s^{(3)}(z)=z \,T(z) ,}
as it then reads
\eqn\eqforTbis{T=1+z^2\, T^3.}
This allows to identify $T$ with the generating function of ternary trees (with a weight $z^2$ per inner vertex). 
Using $[z^{2n}]T={3n \choose n}/(2n+1)$, we arrive at
\eqn\ffourtrizn{[z^{2n+2}] f_4^{(3)}={2\over n+1} \times {{3n\choose n}\over 2n+1}}
for the number of irreducible triangular dissections of the square with $2n+2$ triangles, $n\geq 0$ (the number
of triangles in a dissection of the square must be even - see Fig.~\triquaddiss\ for an example
with $n=7$).
We recover here a result of \TutteCPT.

\newsec{Approach by substitution}

A natural approach to irreducible maps is via
substitution. Intuitively speaking, a $d$-irreducible map is obtained
by erasing the contents of all cycles of length $d$ in a general
map. This naive viewpoint can be made more precise by combining the
following two observations:
\item{(i)} rooted maps of girth at least $d$ and outer degree $n$ are
obtained from $(d-1)$-irreducible maps with outer degree $n$ by
forbidding all inner faces of degree $d-1$;
\item{(ii)} rooted maps of girth at least $d$ and outer degree $n$ are
alternatively obtained from $d$-irreducible maps with outer degree $n$
by {\it substituting} each inner face of degree $d$ with an arbitrary
rooted map of girth $d$ and outer degree $d$.
\par\noindent
Observation (i), which was already made in Sect.~1.2, implies that
the generating function of rooted maps of girth at least $d$ and outer
degree $n$ is equal to
$F_n^{(d-1)}(0;x_d,x_{d+1},\ldots)$. Observation (ii), which we will justify
in the forthcoming subsection, implies that the same generating
function is equal to
$F_n^{(d)}(G_d(x_d,x_{d+1},\ldots);x_{d+1},\ldots)$ where
\eqn\Gd{G_d(x_{d},x_{d+1},\ldots )=F_{d}^{(d-1)}(0;x_d,x_{d+1},\ldots
)-{\rm Cat}(d/2)}
is the generating function of rooted map of girth $d$ and outer degree
$d$ (indeed a map of girth at least $d$ and outer degree $d$ has girth
exactly $d$ unless it is reduced to a tree). 
Since we are expressing the
same quantity in two manners, we get the basic identity
\eqn\substF{F_n^{(d-1)}(0;x_d,x_{d+1},\ldots)=F_n^{(d)}(G_d(x_d,x_{d+1},\ldots);x_{d+1},\ldots).}
Let us now complete the proof of this identity before explaining how it allows
to compute practically $F_n^{(d)}$.

\subsec{The basic substitution relation}
\subseclab\basicsub

We now justify the observation (ii) made above.  More precisely, we
shall prove that we have a bijection between, on the one hand, the set
of rooted maps of girth at least $d$ and outer degree $n$ and, on the
other hand, the set of pairs of the form $({\cal M},(m_f)_{f \in {\cal
F}_d({\cal M})})$ where ${\cal M}$ is a $d$-irreducible map of outer
degree $n$, ${\cal F}_d({\cal M})$ is the set of its $d$-valent inner
faces and, for $f \in {\cal F}_d({\cal M})$, $m_f$ is a rooted map of
girth $d$ and outer degree $d$.

Starting from such a pair $({\cal M},(m_f)_{f \in {\cal F}_d({\cal
M})})$, we define a rooted map ${\cal M}'$ by ``gluing'' inside each
face $f \in {\cal F}_d({\cal M})$ the map $m_f$. More precisely, we
identify clockwise each edge of $m_f$ incident to the outer face with
an edge of ${\cal M}$ incident to $f$, starting from the root edge of
$m_f$ which is identified with an edge of $f$ selected in some
canonical manner (for instance by breadth-first search from the root
of ${\cal M}$). Note that the boundary of $f$ and
that of the outer face of $m_f$ are both simple, thus by
identification they yield a cycle of ${\cal M'}$.  ${\cal M}'$ is a
rooted map (with the same root as ${\cal M}$) of outer degree $n$ and
we claim that its girth is at least $d$.  The proof relies on two
lemma, the first of which will also be useful later on.

\fig{Illustration of the no-shortcut lemma: $v$ and $v'$ are two
vertices $m$ edges away in one direction (hence $d-m$ in the other direction) 
along the outer boundary (blue) of a map of girth $d$ and outer degree $d$. A
simple path from $v$ to $v'$ (red) containing at least one inner edge
necessarily has length $n \geq \max(m,d-m)$.}{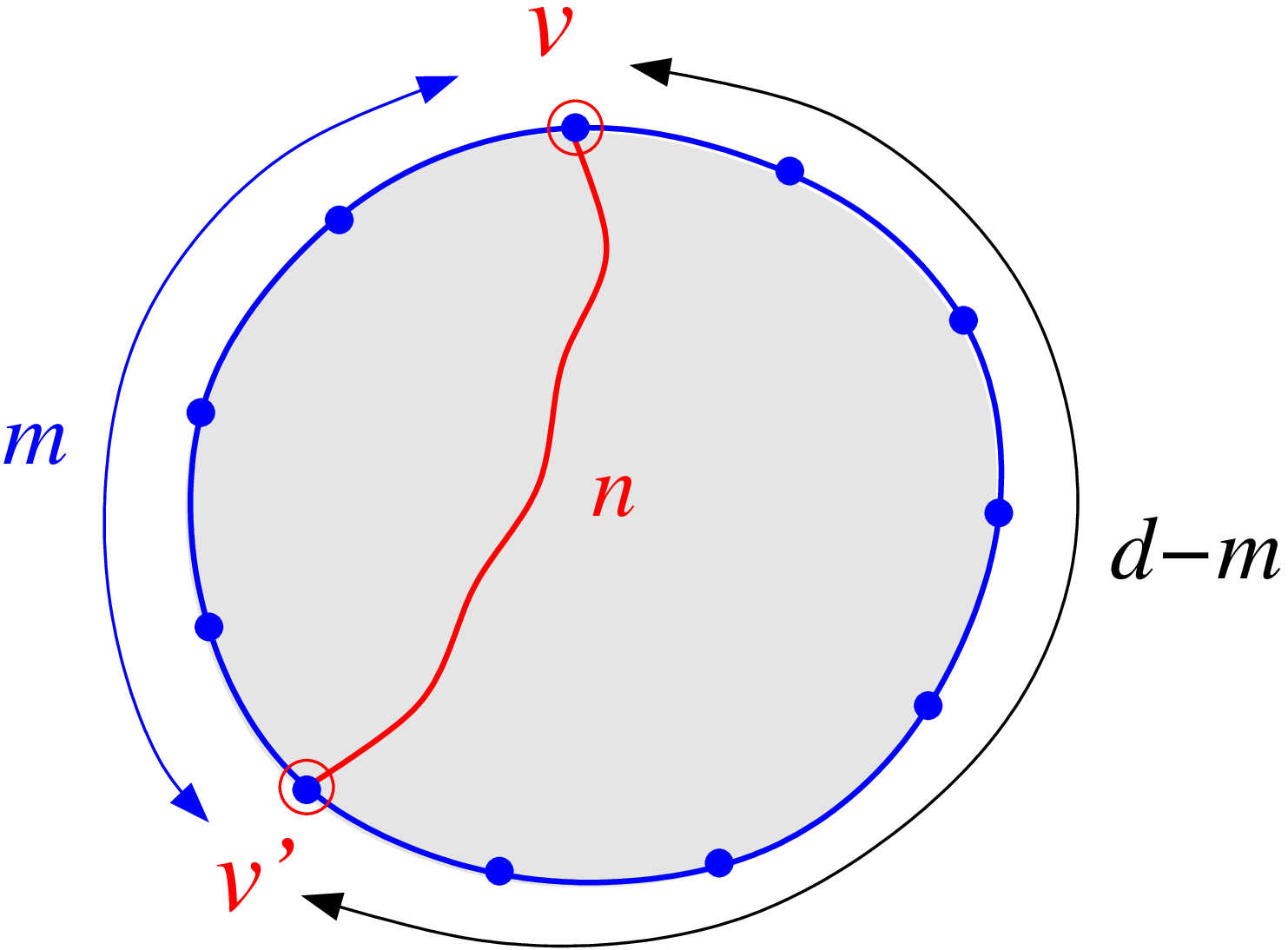}{6.cm} 
\figlabel\geodirr
\smallskip

\item{} {\bf No-shortcut lemma:} in a map of girth $d$ and outer
degree $d$, if $v$ and $v'$ are two outer vertices $m$ edges away in
one direction or the other along the boundary, then any simple path from
$v$ to $v'$ containing at least one inner edge has length larger than
or equal to $\max(m,d-m)$ (which itself is larger than or equal to $d/2$). \smallskip

\noindent {\it Proof:} consider a simple path from $v$ to $v'$
containing at least one inner edge (see Fig.~\geodirr), and denote
its length by $n$. Adding the shorter boundary between $v'$ and $v$,
we obtain a closed path of length $n+\min(m,d-m)$, which is not necessarily
simple but encircles at least one inner face. Thus, by possibly
removing some edges we obtain a cycle of length at most $n+\min(m,d-m)$.
Since the map has girth $d$,  we deduce $n\geq d-\min(m,d-m)=\max(m,d-m)$,
Q.E.D.
\smallskip

\item{} {\bf Encircling lemma:} given a cycle $C'$ of ${\cal M'}$,
  there exists a cycle $C$ of ${\cal M}$ such that $C'$ lies in the
  closed region bounded by $C$, and such that the length of $C$ is at
  most that of $C'$.  \smallskip

\noindent
{\it Proof:} let $C'$ be a cycle of ${\cal M'}$ of length $\ell$. If
$C'$ remains within a single $d$-valent face $f$ of ${\cal M}$, then
the statement is clearly true: $C'$ is a cycle of $m_f$, so $\ell \geq
d$ and we may take $C$ to be the boundary of $f$, which has length $d$. Otherwise, $C'$
necessarily visits at least two vertices of ${\cal M}$. We orient $C'$
counterclockwise, and denote by $v_1,v_2,\ldots,v_p$ ($2 \leq p \leq
\ell$) the successive vertices of ${\cal M}$ visited along $C'$. We
then let $\gamma_i'$ ($1 \leq i \leq p$) be the part of $C'$ comprised
between $v_i$ and $v_{i+1}$ (with $v_{p+1}=v_1$). $\gamma_i'$ is a
simple open oriented path which is either reduced to a single edge of
${\cal M}$, or made of edges not in ${\cal M}$ that are all inside a
same $d$-valent face $f_i$ of ${\cal M}$. In the former case we let
$\gamma_i=\gamma_i'$, while in the latter case we let $\gamma_i$ be
the part of the boundary of $f_i$ going from $v_i$ to $v_{i+1}$ in the
counterclockwise direction: the no-shortcut lemma ensures that the
length of $\gamma_i$ is at most that of $\gamma_i'$. Upon
concatenating $\gamma_1,\cdots,\gamma_p$ together we obtain a closed
path ${\tilde C}$ on ${\cal M}$ of length at most $\ell$. However, we
are not assured that it is a cycle encircling $C'$ so we cannot yet
conclude.  Instead, we let $R'$ be the closed region bounded by $C'$
and, for $1 \leq i \leq p$, $r_i$ be the closed region bounded by
$\gamma_i \cup \gamma'_i$ if $\gamma_i \neq \gamma'_i$, or the empty
set otherwise. Then, $R = R' \cup r_1 \cup \cdots \cup r_p$ is a
simply connected closed region whose boundary is the cycle $C$ we are
looking for. Indeed, note first that it obviously contains
$C'$. Secondly, if $r_i$ is nonempty, then $\gamma'_i$ has $r_i$ on
its right and $R'$ on its left, thus the interior of $R$ is simply
connected and, furthermore, any edge of $C'$ not in ${\cal M}$ cannot
be on the boundary of $R$. Therefore, the boundary of $R$ is a cycle
$C$ whose edges form a subset of those of ${\tilde C}$, thus its
length is at most $\ell$, Q.E.D.  \smallskip

The encircling lemma immediately implies that the girth of ${\cal M'}$
is at least $d$ (more precisely it is equal to the girth of ${\cal M}$
which is at least $d$, with equality iff ${\cal F}_d({\cal M})$ is
nonempty). Furthermore, since ${\cal M}$ is irreducible, this lemma also
implies that any cycle of
length $d$ in ${\cal M'}$ necessarily remains within a single
$d$-valent face of ${\cal M}$: the boundaries of $d$-valent faces of
${\cal M}$ are thus precisely the {\it outermost} cycles of length $d$
in ${\cal M}'$. This shows that, starting from ${\cal M'}$, we may
recover ${\cal M}$ by erasing all the edges and vertices that are
interior to the outermost $d$-cycles. Then, for a given $f \in {\cal
F}_d({\cal M})$, the edges and vertices of ${\cal M'}$ that lie within
$f$ form the map $m_f$. In conclusion, the mapping that maps $({\cal
M},(m_f)_{f \in {\cal F}_d({\cal M})})$ to ${\cal M'}$ is
injective. It remains to check that it is surjective.

\bigskip
\fig{Illustration of the proof that two outermost cycles of length $d$ in
a map of girth $d$ do not overlap. If two cycles $C_1$ and $C_2$ 
of length $d$ overlap (hence we may find two vertices $v$ and $v'$ 
such that the part of $C_1$ between $v$ and $v'$, of length $n_1$
lies in the interior of $C_2$ and vice versa), we can build a cycle $C_{12}$ (dashed line)
encircling them both, and 
having length at most $2d-(n_1+n_2)\leq d$, hence equal to $d$.}{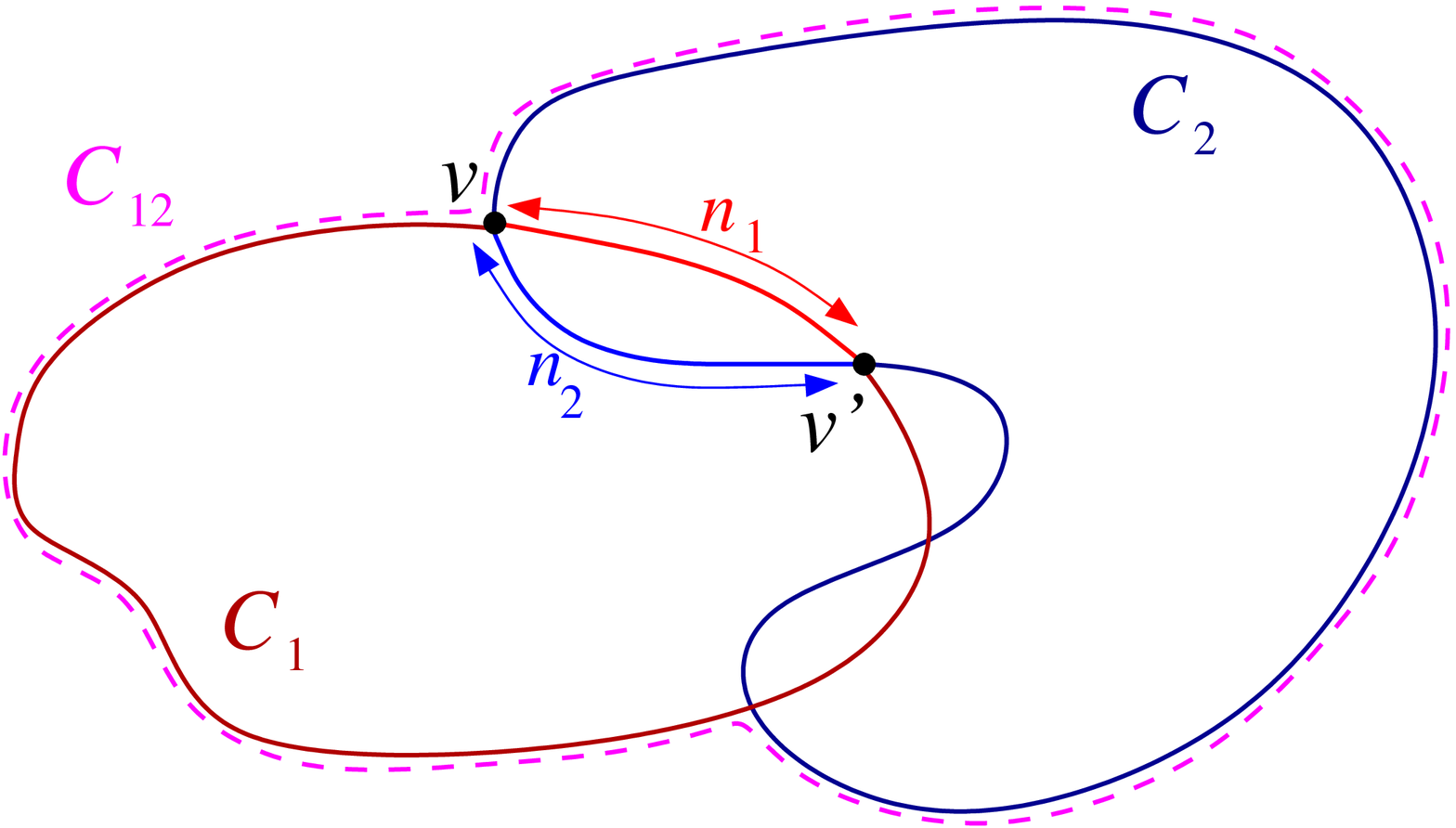}{8.cm}
\figlabel\twocycles

Let us now start conversely with an arbitrary rooted map ${\cal M}'$
of outer degree $n$ and girth at least $d$ and consider the set ${\cal
C}_d({\cal M}')$ of its {\it outermost} cycles of length $d$ , i.e.\
those cycles of length $d$ whose interior is not strictly included
into the interior of another cycle of length $d$. For $C \in {\cal
C}_d({\cal M}')$, the edges and vertices of ${\cal M'}$ that lie
within $C$ form a map $m_C$ which clearly has girth $d$ and outer
degree $d$ (we discuss its rooting below).
Now we have the crucial property that two distinct cycles in ${\cal
C}_d({\cal M}')$, say $C_1$ and $C_2$, cannot overlap, i.e.\ the
intersection of their interiors is necessarily empty.  Indeed since
the interior of one cycle cannot be included in the interior of the
other, if we assume that these interiors have a nonempty
intersection, there exist two vertices $v$ and $v'$ at the
intersection of $C_1$ and $C_2$ such that one of the parts of $C_1$
between $v$ and $v'$ lies in the interior of $C_2$ and one of the
parts of $C_2$ between $v$ and $v'$ lies in the interior of $C_1$ (see
Fig.~\twocycles).  Calling $n_1$ and $n_2$ the lengths of these
parts, we have $n_1\geq d/2$ and $n_2\geq d/2$ by the no-shortcut
lemma applied to $m_{C_2}$ and $m_{C_1}$ respectively. We can then
build a cycle $C_{12}$ by following the outer boundary of the union
of the interiors of $C_1$ and $C_2$ (note that this union is connected but
not necessarily simply connected -- see Fig.~\twocycles).  This cycle has length at most
$(d-n_1)+(d-n_2)\leq d$, hence it has length $d$ from the girth
condition on ${\cal M}'$. The interiors of $C_1$ and $C_2$ are then
strictly included in the interior of the cycle $C_{12}$ of length $d$,
a contradiction.  Since outermost cycles of length $d$ do not overlap,
we may unambiguously replace the content of each such cycle by a
simple face of degree $d$, resulting in a rooted map ${\cal M}$ of
outer degree $n$ and girth at least $d$, such that all cycles of
length $d$ are the boundary of an inner face of degree $d$, hence a
$d$-irreducible map (note that the outer face of ${\cal M}'$ is
unaffected by the substitution since it cannot belong to the interior
of a cycle). In particular, we may identify ${\cal F}_d({\cal M})$
with ${\cal C}_d({\cal M}')$. For each $f \in {\cal F}_d({\cal M})$,
associated with $C\in {\cal C}_d({\cal M}')$,
we select an edge of ${\cal M}$ incident to $f$ by the same canonical
procedure as before: this provides a canonical rooting of
$m_f\equiv m_C$. Obviously, applying the previous construction to $({\cal
M},(m_f)_{f \in {\cal F}_d({\cal M})})$ restores ${\cal M'}$, hence
the mapping from $({\cal M},(m_f)_{f \in {\cal F}_d({\cal M})})$ to ${\cal M'}$ is
surjective. It is therefore a bijection.

As a final remark, note that this bijection preserves the following
parameters:
\item{-} the number of $d$-valent inner faces of ${\cal M'}$ is equal
to the total number of $d$-valent inner faces in all $m_f$, $f \in
{\cal F}_d({\cal M})$,
\item{-} for each $k > d$, the number of $k$-valent inner faces of
${\cal M'}$ is equal to the total number of $k$-valent inner faces in
${\cal M}$ and all $m_f$, $f \in {\cal F}_d({\cal M})$.
\par \noindent
It follows that the generating function of rooted map of girth at
least $d$ and outer degree $n$ is indeed equal to
$F_n^{(d)}(G_d(x_d,x_{d+1},\ldots);x_{d+1},\ldots)$, which concludes
the proof of the basic identity \substF.

\subsec{General strategy}

We turn to discussing the practical use of the basic identity \substF\
in computing $F_n^{(d)}$.  We first claim that there exists a series
$X_d(z;x_{d+1},\ldots)$ such that
\eqn\Gcompinv{G_d(X_d(z;x_{d+1},\ldots),x_{d+1},\ldots)=z,}
in other words $G_d$ admits a {\it compositional inverse} with respect
to its first variable. Observe indeed that $G_d$ is a series in
$x_d,x_{d+1},\ldots$ whose constant term is zero (since any map
contributing to $G_d$ contains at least one inner face) and whose
coefficient of $x_d$ is $1$ (corresponding to the map reduced to a
single $d$-gon), i.e.
\eqn\Gdecomp{G_d(x_d,x_{d+1},\ldots) = A(x_{d+1},\ldots) + x_d (1 + B(x_{d+1},\ldots)) + x_d^2 C(x_d,x_{d+1},\ldots)}
where $A,B$ are formal power series in $x_{d+1},x_{d+2},\ldots$
without constant term. Since $1+B$ is invertible in the ring of formal
power series in $x_{d+1},x_{d+2},\ldots$, the series $x_d (1+B) +
x_d^2 C$, viewed as a series in $x_d$ whose coefficients
are series in $x_{d+1},x_{d+2},\ldots$, admits a compositional inverse
$D$ with respect to the variable $x_d$,
satisfying $D(0,x_{d+1},x_{d+2},\ldots)=0$.
The wanted series is then
$X_d(z;x_{d+1},\ldots)=D(z-A(x_{d+1},\ldots),x_{d+1},\ldots)$, the
substitution being well defined as $A$ has no constant term. Note
that, when specializing $x_i=0$ for all odd $i$, we have
$G_d(0,x_{d+1},\ldots)=A(x_{d+1},\ldots)=0$ for $d$ odd (since there
are no maps with odd outer degree and all inner faces of even degree),
hence $X_d(0;x_{d+1},\ldots)=D(0,x_{d+1},\ldots)=0$.

Replacing $x_d$ by $X_d(z;x_{d+1},\ldots)$ in \substF, we obtain the
reciprocal identity
\eqn\substFrecip{F_n^{(d)}(z;x_{d+1},\ldots) = F_n^{(d-1)}(0;X_d(z;x_{d+1},\ldots),x_{d+1},\ldots).}
By iterating this relation $d$ times, we relate $F_n^{(d)}$ to the
generating function $F_n(x_1,x_2,\ldots)=F_n^{(0)}(0;x_1,x_2,\ldots)$
of arbitrary maps with outer degree $n$. Namely, we have the
general substitution relation
\eqn\FndFnexpr{F_n^{(d)} = F_n(X_1^{(d)},X_2^{(d)},\ldots,X_d^{(d)},x_{d+1},\ldots)}
where the series $X_j^{(d)}(z;x_{d+1},\ldots)$, $1 \leq j \leq d$, are
defined inductively by
\eqn\Xkddef{X_j^{(d)}(z;x_{d+1},\ldots) = \left\{\matrix{
X_d(z;x_{d+1},\ldots) & \hbox{for $j = d$,} \cr
X_j(0;X_{j+1}^{(d)},\ldots,X_d^{(d)},x_{d+1},\ldots) & \hbox{for $j <
d$.} }\right.}
Note that, in the bipartite case ($n,d$ even, $x_i=0$ for $i$ odd), we
showed above that $X_j(0;x_{j+1},\ldots)=0$ for $j$ odd so that, by
\Xkddef, $X_j^{(d)}=0$. In other words, \FndFnexpr\ relates the
generating function of bipartite $d$-irreducible maps to the
generating function of bipartite maps without irreducibility
constraints (alternatively, this can be shown by writing the bipartite
analogue of \substF, relating directly $F_n^{(d-2)}$ and $F_n^{(d)}$).

While $F_n$ is a well studied quantity for which convenient
expressions are known (see below), we have {\it a priori} no such
expressions for $X_1^{(d)},\ldots,X_d^{(d)}$ which appear in
\FndFnexpr. However, these $d$ unknown quantities may in principle be {\it
determined} by the conditions \Fsmalln, and then {\it eliminated} from
the expression of $F_n^{(d)}$. We have seen in Sect.~2 two cases where
this elimination can been carried out smoothly, and this can be done in
general, as we shall see in the following subsections.
We first concentrate
on the bipartite case where the expressions are somewhat simpler.

\subsec{Elimination in the bipartite case}
\subseclab\elimbip

In the bipartite case, both the outer degree $n$ and the girth $d$ are
even integers, hence we write $n=2m$, $d=2b$. The generating function $F_{2m}$ of bipartite,
non necessarily irreducible, maps is given by [\xref\DFG,\xref\PMCF]
\eqn\Fnbip{F_{2m} = {\rm Cat}(m) R^m - \sum_{j \geq 1} \sum_{k = 1}^{\min(m,j-1)}
{2k+1 \over 2m+1} {2m+1 \choose m-k} {2j-1 \choose j+k} x_{2j} R^{m+j}}
where $R$ is the formal power series determined by the equation
\eqn\eqforR{R = 1 + \sum_{j \geq 1} {2j-1 \choose j} x_{2j} R^j.}
For concreteness, let us mention that the general coefficient of $R$ reads explicitly
\eqn\coeffR{\Big[ \prod_{j \geq 1} (x_{2j})^{n_j} \Big] R =
  { \left( \sum_{j \geq 1} j n_j \right)! \over \left( 1 + \sum_{j \geq 1} (j-1) n_j \right)!}
  \prod_{j \geq 1} {1 \over n_j!} {2j-1 \choose j}^{n_j}}
as seen by applying the Lagrange inversion formula [\xref\FS-\xref\GouldJack].
We may rewrite \Fnbip\ in a more compact form by introducing  the shorthand notations
\eqn\Amj{A_{m,k} = {2k+1 \over 2m+1} {2m+1 \choose m-k}}
(note that ${\rm Cat}(m)=A_{m,0}$) and
\eqn\Ujdef{U_k = \sum_{j \geq k+1} {2j-1 \choose j+k} x_{2j} R^{j+k}}
so that
\eqn\Fnbipcomp{F_{2m} = R^m \sum_{k=0}^m A_{m,k} \left( \delta_{k,0} - U_k R^{-k}(1-\delta_{k,0}) \right).}
We now apply the general substitution relation \FndFnexpr, to get
\eqn\Fndbipcomp{F_{2m}^{(d)} =
  (R^{(d)})^m \sum_{k=0}^m A_{m,k} \left( \delta_{k,0} - U_k^{(d)} (R^{(d)})^{-k}(1-\delta_{k,0}) \right)}
where $R^{(d)}$ and $U_k^{(d)}$ are the series obtained by substituting, for all $j$ between $1$ and $b$,
the formal variable $x_{2j}$ by the series $X_{2j}^{(d)}$ in $R$ and $U_k$ respectively (recall that
$X_j^{(d)}=0$ for $j$ odd in the bipartite case). In particular, since the variables $x_{2j}$ with $j>b$ are
unaffected by the substitution, we have
\eqn\Ujbig{U_k^{(d)} = \sum_{j \geq k+1} {2j-1 \choose j+k} x_{2j} (R^{(d)})^{j+k} \qquad \hbox{for $k \geq b$}.}
We are therefore left with the $b$ unknown quantities $R^{(d)}$ and $U_k^{(d)}$, $k=1,\ldots,b-1$, 
which replace the original unknowns $X_{2j}^{(d)}$, $j=1,\ldots, b$.

We then observe that the condition \Fsmalln\ implies that
$F_{2m}^{(d)} = {\rm Cat}(m)$ for $1 \leq m \leq b-1$, which may be
viewed as a system of $b-1$ {\it linear} equations for 
$U_1^{(d)},\ldots ,U_{b-1}^{(d)}$ (note that the condition
$F_0^{(d)}=1$ is readily satisfied since $F_0=1$). It may easily be
solved by introducing the inverse $B=(B_{n,m})_{n,m \geq 0}$ of the semi-infinite unitriangular
matrix $A=(A_{m,k})_{m,k \geq 0}$, whose coefficients read explicitly
\eqn\Bnm{B_{n,m} = (-1)^{n+m} {n+m \choose 2m}.}
Multiplying \Fndbipcomp\ by $B_{n,m} (R^{(d)})^{-m}$ and summing over $m$ from $0$ to $n$, we obtain that
\eqn\Usol{U_n^{(d)} = - \sum_{m=0}^n B_{n,m} {\rm Cat}(m)  (R^{(d)})^{n-m} \qquad \hbox{for $1 \leq n \leq b-1$}.}
Plugging \Ujbig\ and \Usol\ into \Fndbipcomp\
yields
\eqn\Fndbipplug{F_{2m}^{(d)} =
  \sum_{k=0}^{b-1} \sum_{\ell=0}^k A_{m,k} B_{k,\ell}  {\rm Cat}(\ell) (R^{(d)})^{m-\ell}
  - \sum_{k=b}^m \sum_{j \geq k+1} A_{m,k} {2j-1 \choose j+k} x_{2j} (R^{(d)})^{m+j}.}
This expression may be further simplified using the two hypergeometric identities
\eqn\hypergeomone{\sum_{k=\ell}^{b-1} A_{m,k} B_{k,\ell} =
  (-1)^{b-\ell-1} {b-\ell \over m-\ell} {2m \choose m-b} {b+\ell \choose 2\ell},}
\eqn\hypergeomtwo{\sum_{k=b}^m A_{m,k} {2j-1 \choose j+k} =
  {b+j \over m+j} {2m \choose m-b} {2j-1 \choose j+b},}
which are routinely obtained via Gosper's algorithm and may
easily be checked by induction on $b$, resulting in the expression
\eqn\Fndbipfin{\eqalign{F_{2m}^{(d)} = {2m \choose m-b}
  & \Bigg(\sum_{\ell=0}^{b-1}(-1)^{b-\ell-1} {b-\ell \over m-\ell} {b+\ell \choose 2\ell} {\rm Cat}(\ell) (R^{(d)})^{m-\ell} \cr
  & \qquad  - \sum_{j \geq b+1} {b+j \over m+j} {2j-1 \choose j+b} x_{2j} (R^{(d)})^{m+j} \Bigg).\cr} }
Interestingly, this expression involves only the (so far unknown) series $R^{(d)}$. This quantity is
in turn determined by the last case of condition \Fsmalln, namely that $F_{d}^{(d)} = {\rm Cat}(b) + z$, which may be rewritten in the form
\eqn\eqforRdbip{z + \sum_{\ell=0}^{b} (-1)^{b-\ell} {b+\ell \choose 2\ell} {\rm Cat}(\ell) (R^{(d)})^{b-\ell}
+ \sum_{j \geq b+1} {2j-1 \choose j+b} x_{2j} (R^{(d)})^{b+j}=0.}
In particular, if we impose a bound on the face degrees
(i.e. we take $x_{2j}=0$ for $j$ larger than some fixed $M$), then
$R^{(d)}$ hence $F_m^{(d)}$ are algebraic. This is notably the case
for irreducible $d$-angular dissections, where we specialize
$x_{2j}=0$ for all $j$, keeping $z$ as the only formal variable, to get
\eqn\fndbipfin{f_{2m}^{(d)}(z) = {2m \choose m-b} \sum_{\ell=0}^{b-1}(-1)^{b-\ell-1} {b-\ell \over m-\ell} {b+\ell \choose 2\ell} {\rm Cat}(\ell) (r^{(d)})^{m-\ell}}
where $r^{(d)}(z)$, the corresponding specialization of $R^{(d)}$, satisfies the algebraic equation
\eqn\eqforrdbip{z + \sum_{\ell=0}^{b} (-1)^{b-\ell} {b+\ell \choose 2\ell} {\rm Cat}(\ell) (r^{(d)})^{b-\ell} = 0.}

A remarkable formula follows from differentiating \Fndbipfin\ wrt $R^{(d)}$ (keeping the $x_{2j}$, $j>b$, fixed):
observe that
\eqn\Fndbipdiff{{\partial F_{2m}^{(d)} \over \partial R^{(d)}} = {2m \choose m-b} (R^{(d)})^{m-b} {\partial F_{d}^{(d)} \over \partial R^{(d)}}.}
Since $F_d^{(d)}={\rm Cat}(b) + z$, it follows, by multiplying both sides by ${\partial R^{(d)} \over \partial z}$, that we have the {\it pointing formula}
\eqn\Fndbippoint{{\partial F_{2m}^{(d)} \over \partial z} = {2m \choose m-b} (R^{(d)})^{m-b}}
and, in particular, for $m=b+1$,
\eqn\Fndbippointpar{{\partial F_{d+2}^{(d)} \over \partial z} = (d+2) R^{(d)}.}
These formulas generalize, in some sense, the expression \PMCF
\eqn\Fnbippointprev{F_{2m}^{\bullet} = {2m \choose m} (R^{(d)})^{m}}
for the generating function of {\it pointed} rooted bipartite
(non necessarily irreducible) planar maps, which we recover in the
case $d=0$ (upon understanding $z$ as a weight per vertex). 

\fig{An example of irreducible hexangular dissection of the octagon.}{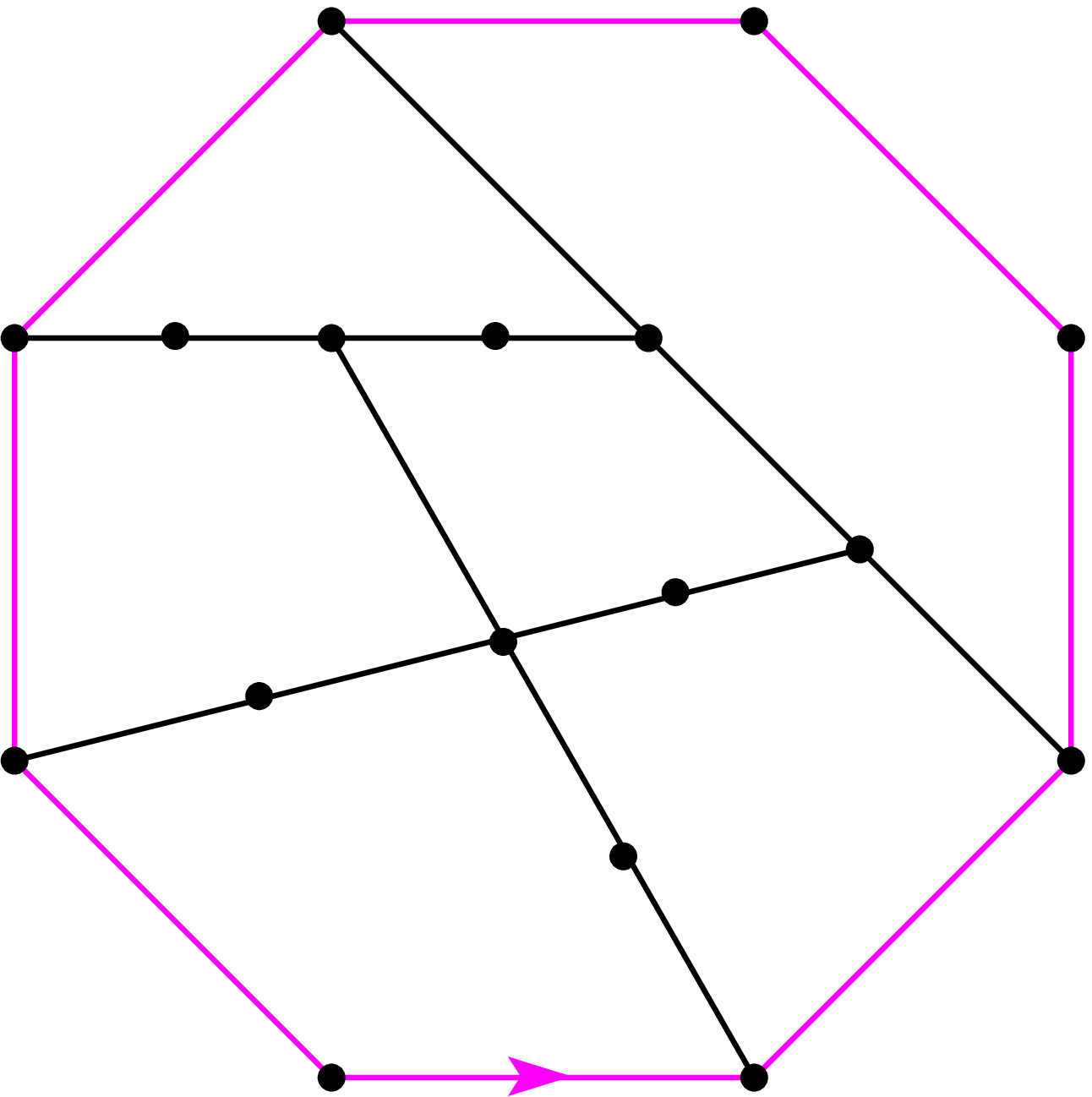}{4.5cm}
\figlabel\hexdiss

Let us now give some examples. Taking $b=2$ and $x_i=0$ for all $i$,
we recover the case of irreducible quadrangular dissections discussed
in Sect.~\irredquad: \fndbipfin\ yields the expression \fnquad\ for
the generating function of irreducible quadrangular dissections of
the $2m$-gon and \eqforrdbip\ yields the algebraic equation
\algebRquad\ for $r^{(4)}$. Taking now $b=3$ and still $x_i=0$ for all
$i$, we obtain generating functions of irreducible hexangular
dissections (see Fig.~\hexdiss), namely
\eqn\fnhex{f_{2m}^{(6)}(z) = {2m \choose m-3} \left( {10 \over m-2} (r^{(6)})^{m-2} - {12 \over m-1} (r^{(6)})^{m-1} + {3 \over m} (r^{(6)})^{m} \right)}
where $r^{(6)}$ satisfies
\eqn\eqforrhex{z-(r^{(6)})^3 + 6\, (r^{(6)})^2 - 10\, r^{(6)} + 5 = 0.}
The first few terms read
\eqn\hexfft{\eqalign{
  r^{(6)} &= 1 + z + 3\,z^2 + 17\,z^3 + 120\,z^4 + 948\,z^5 + 8022\,z^6 + \cdots\cr
  f_8^{(6)} &= 14 + 8\,z + 4\,z^2 + 8\,z^3 + 34\,z^4 + 192\,z^5 + 1264\,z^6 + 9168\,z^7 + \cdots\cr
  f_{10}^{(6)} &= 42 + 45\,z + 45\,z^2 + 105\,z^3 + 450\,z^4 + 2547\,z^5 + 16785\,z^6 + 121815\,z^7 + \cdots\cr}}
consistently with the pointing formula.

\subsec{Elimination in the general case}
\subseclab\elimgen

We now repeat the same strategy in the general case, i.e.\ when the
maps are non necessarily bipartite. Our starting point is the
expression found in \PMCF\ for the generating function $F_n$ of
general maps with outer degree $n$. 
It involves generating functions of {\it three-step paths}, i.e.\ lattice paths
in $Z^2$ made of three types of steps: up-steps $(1,1)$, down-steps $(1,-1)$ and
level-steps $(1,0)$. We denote by
\eqn\pathgf{\eqalign{
  P_k(n;r,s) &= \sum_{j=0}^{\lfloor {n-k \over 2} \rfloor} {n! \over j! (j+k)! (n-2j-k)!} r^j s^{n-2j-k}  \cr
  P^+_k(n;r,s) &= \sum_{j=0}^{\lfloor {n-k \over 2} \rfloor} {(k+1)\, n! \over j! (j+k+1)! (n-2j-k)!} r^j s^{n-2j-k}}}
the  generating polynomials of respectively arbitrary and
nonnegative three-step paths from $(0,0)$ to $(n,k)$, counted with a weight $r$ per down-step and $s$ per
level step (a path is said nonnegative if it only visits vertices with nonnegative ordinates). 
Note that, for $k \leq 0$,  $P_k(n;r,s)$ is well defined (the sum over $j$ then starts in practice at $j=-k$) and
equals $r^{-k} P_{-k}(n;r,s)$. Then, we have \PMCF
\eqn\Fngen{F_n = P_0^+(n;R,S) - \sum_{k \geq 1} P_k^+(n;R,S) V_k}
where $V_k$ is defined by
\eqn\Vjdef{V_k = \sum_{j \geq k+2} x_j P_{-k-1}(j-1;R,S) }
and $R,S$ are formal power series determined by the equations
\eqn\eqsforRS{R=1+V_0 \qquad S=V_{-1}.}
Clearly, \Fngen\ yields \Fnbipcomp\ in the bipartite case where $S=0$,
$P_{2k}^+(2m;R,0)=A_{m,k} R^{m-k}$, $V_{2k}=U_k$ , $V_{2k+1}=0$.

Applying the general substitution relation \FndFnexpr, we get
\eqn\Fndgen{F_n^{(d)} = P_0^+(n;R^{(d)},S^{(d)}) - \sum_{k \geq 1} P_k^+(n;R^{(d)},S^{(d)}) V_k^{(d)}}
where $R^{(d)}$, $S^{(d)}$ and $V_k^{(d)}$ are the series obtained by
substituting, for all $j$ between $1$ and $d$, the formal variable
$x_j$ by the series $X_j^{(d)}$ in $R$, $S$ and $V_k$ respectively, namely
\eqn\RSVjdexpr{\eqalign{
R^{(d)} & =R(X_1^{(d)},X_2^{(d)},\ldots,X_d^{(d)},x_{d+1},\ldots)\cr
S^{(d)} & =S(X_1^{(d)},X_2^{(d)},\ldots,X_d^{(d)},x_{d+1},\ldots)\cr
V_k^{(d)}& =V_k(X_1^{(d)},X_2^{(d)},\ldots,X_d^{(d)},x_{d+1},\ldots)\ .\cr}}
Note that, in particular,
\eqn\Vjdbig{V_k^{(d)} = \sum_{j \geq k+2} x_j P_{-k-1}(j-1;R^{(d)},S^{(d)}) \qquad \hbox{for $k \geq d-1$}.}
We are left with $d$ unknowns $R^{(d)}$, $S^{(d)}$ and $V_1^{(d)},\ldots, V_{d-2}^{(d)}$ which 
replace the original unknowns $X_1^{(d)},\ldots,X_d^{(d)}$.
Similarly to the previous section, we may determine these quantities using the conditions \Fsmalln. Here,
we need the inverse of the semi-infinite unitriangular matrix
$(P_k^+(n;r,s))_{n,k\geq 0}$, which is denoted by $(Q_{n,k}(r,s))_{n,k \geq 0}$
with the explicit form
\eqn\Qdef{Q_{n,k}(r,s) = \sum_{j=0}^{\lfloor {n-k \over 2} \rfloor}
  {(n-j)! \over k! j! (n-2j-k)!} (-r)^j (-s)^{n-2j-k}.}
Multiplying \Fndgen\ by $Q_{k,n}(R^{(d)},S^{(d)})$ and summing over $n$, the conditions \Fsmalln\ amount to
\eqn\Vjdsol{V_k^{(d)} = - \sum_{j=0}^{\lfloor k/2 \rfloor} Q_{k,2j}(R^{(d)},S^{(d)})\, {\rm Cat}(j) - \delta_{k,d} z\qquad
\hbox{for $1\leq k \leq d$}}
(recall that ${\rm Cat}(j)=0$ for noninteger $j$). In particular,
comparing this expression with \Vjdbig\ for $k=d-1$ or $d$ yields 
the equations
\eqn\eqsforRSd{\eqalign{
  \sum_{j=0}^{\lfloor {d-1 \over 2} \rfloor} Q_{d-1,2j}(R^{(d)},S^{(d)})\, {\rm Cat}(j) + \sum_{j \geq d+1} x_j P_{-d}(j-1;R^{(d)},S^{(d)}) &=0 \cr
  z + \sum_{j=0}^{\lfloor d/2 \rfloor} Q_{d,2j}(R^{(d)},S^{(d)})\, {\rm Cat}(j) + \sum_{j \geq d+2} x_j P_{-d-1}(j-1;R^{(d)},S^{(d)}) &=0 \cr }}
which determine the power series $R^{(d)}$ and $S^{(d)}$. As
previously, if we impose a bound on the face degrees ($x_j=0$ for
$j>M$) then $R^{(d)}$ and $S^{(d)}$ are algebraic. It follows that
$F_n^{(d)}$ is also algebraic, since by \Vjdsol\ (for $k \leq d-2$)
and \Vjdbig\ (for $k \geq d-1$), \Fndgen\ is a polynomial in $R^{(d)}$
and $S^{(d)}$. This polynomial admits an expression similar to
\Fndbipplug\ ({\it mutatis mutandis}) which we do not find particularly
illuminating: no simplification as nice as \Fndbipfin\ has been found.
Remarkably however, a pointing formula stills holds in the form
\eqn\Fndpointgen{{\partial F_n^{(d)} \over \partial z} = P_d(n;R^{(d)},S^{(d)})}
and may be viewed as a generalization of the formula \PMCF
\eqn\Fnpointgenprev{F_n^{\bullet} = P_0(n;R,S)}
for the generating function of pointed rooted maps.
We will prove \Fndpointgen\ in Sect.~7 below using a combinatorial argument.

We now conclude this section by some examples. Taking $d=3$ and
$x_i=0$ for all $i$, we recover the case of irreducible triangular
dissections discussed in Sect.~\irredquad: \eqsforRSd\ yields the
algebraic equations $1 - r^{(3)} + (s^{(3)})^2 = 0$ and $z + 2 r^{(3)}
s^{(3)} - (s^{(3)})^3 - 3 s^{(3)} = 0$, which amount to
\algebRStri. The pointing formula \Fndpointgen\ yields \dffourtridz\
for $n=4$, which in turn implies the expression \ffourtrizn\ for the
number of irreducible triangular dissections of the square. Taking now
$d=5$ and still $x_i=0$ for all $i$, we obtain the case of irreducible
pentagular dissections. From \eqsforRSd\ we find that $r^{(5)}$ and
$s^{(5)}$ are determined by the algebraic equations
\eqn\SRfive{\eqalign{&(r^{(5)})^2 - 3 r^{(5)} (s^{(5)})^2 +
(s^{(5)})^4 - 3 r^{(5)} + 6 (s^{(5)})^2 + 2 = 0\cr
z - s^{(5)} & \left( 3 (r^{(5)})^2 - 4 r^{(5)} (s^{(5)})^2 + (s^{(5)})^4 - 12 r^{(5)} + 10 (s^{(5)})^2 + 10 \right) = 0\cr}}
and their first few terms read
\eqn\SRfiveexp{\eqalign{r^{(5)}&= 1 + 3 z^2 + 73 z^4 + 3015 z^6 + 151842 z^8 + 
    8493934 z^{10} + 507165545 z^{12}+\ldots \cr
    s^{(5)} &= z + 12 z^3 + 422 z^5 + 19780 z^7 + 1062275 z^9 + 
    61781482 z^{11} + 3786534059 z^{13} +\ldots. \cr}}
By the pointing formula $df_6^{(5)}/dz = 6 s^{(5)}$, the
numbers of irreducible pentagular dissections of the hexagon with
up to 14 inner faces are read off
\eqn\fsixfive{f_6^{(5)} = 5 + 3 z^2 + 18 z^4 + 422 z^6 + 14835 z^8 + 637365 z^{10} + 30890741 z^{12} + 1622800311 z^{14} + \ldots.}

\newsec{Combinatorial interpretation via slices}

In this section, we provide a combinatorial interpretation for some of
the quantities that appear in the approach by substitution. We are led
to define some particular classes of maps called {\it slices},
generalizing the notion introduced in [\xref\PMCF, Appendix A].

\subsec{General slices}
\seclab\genslices
\fig{Schematic picture of slices of type $p/p$ (left) and type $p/p+1$ (right).}{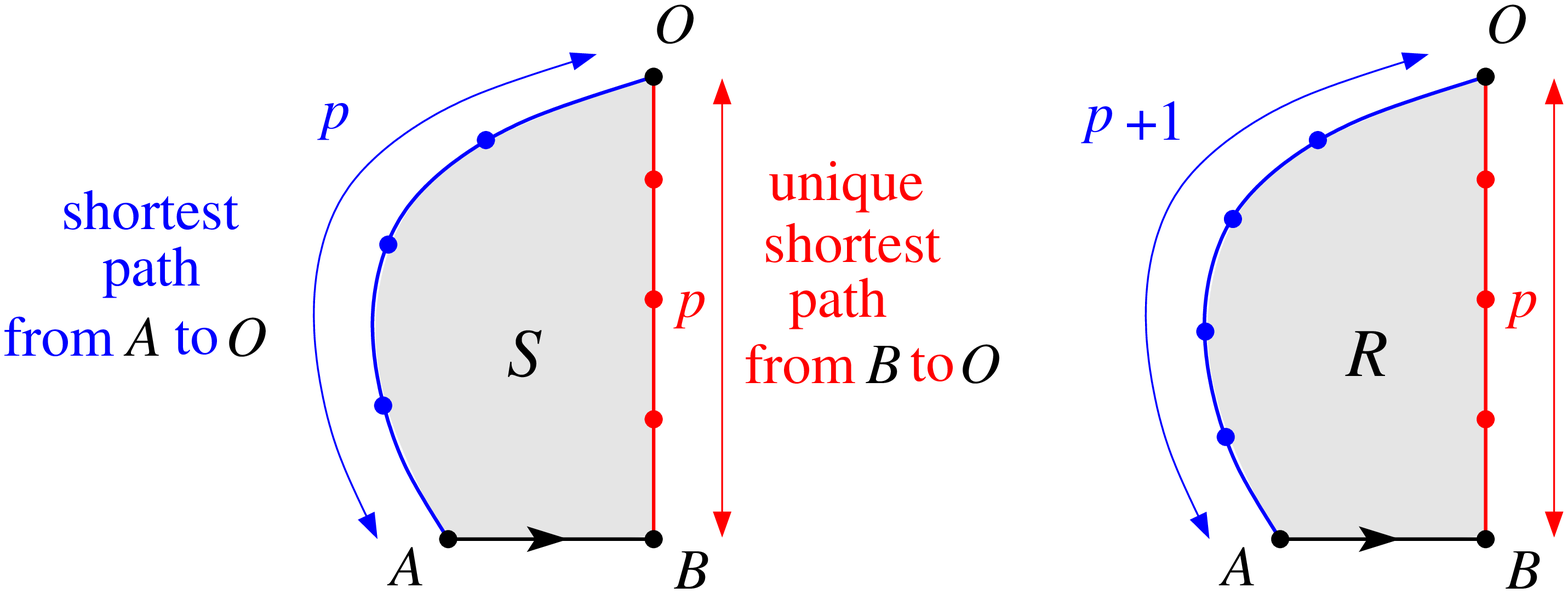}{11.cm}
\figlabel\slices
As explained in \PMCF, eqs.~\Fngen\ and \Fnpointgenprev\ have a direct combinatorial interpretation resulting from a decomposition of the maps enumerated 
respectively by $F_n$ and $F_n^{\bullet}$ into more primitive components called {\it slices}, of which $R$ and $S$ are generating functions.
More precisely, a slice is defined as a rooted map with a marked vertex $O$ (later called the apex) incident to its outer face, and which is of the type displayed
in Fig.~\slices, namely satisfies:
\item{-} the right boundary of the map, defined as the path joining the endpoint of the root edge to $O$ by following the outer face
counterclockwise around the rest of the map, is the unique shortest path in the map between these two points;
\item{-} the left boundary of the map, defined as the path joining the origin of the root edge to $O$ by following the outer face
clockwise around the rest of the map, is a shortest path in the map between these two points;
\item{-} the vertex $O$ is the only vertex common to both the right and left boundaries.
\par\noindent
Clearly, if we denote by $p\geq 0$ the length of the right boundary, that of the left boundary is either $p$ or $p+1$ (it cannot be $p-1$ as otherwise,
the right boundary would not be the unique shortest path). We shall refer to these slices as being of type $p/p$ or $p/p+1$ accordingly.
It was shown in \PMCF\ that $R(x_1,x_2,\ldots)$ is precisely the generating function of slices of type $p/p+1$ for some (unfixed) $p$, counted with weights 
$x_i$ per inner face of degree $i$, while $S(x_1,x_2,\ldots) $ is the generating function of slices of type $p/p$ for some $p$. 
Note that $R$ incorporates a term $1$ accounting for the slice of type $0/1$ reduced to the root edge (with no inner face and with outer face of degree 2),
$O$ being the endpoint of this edge.
\fig{Schematic picture of a $k$-slice of type $p/p+k+1$.}{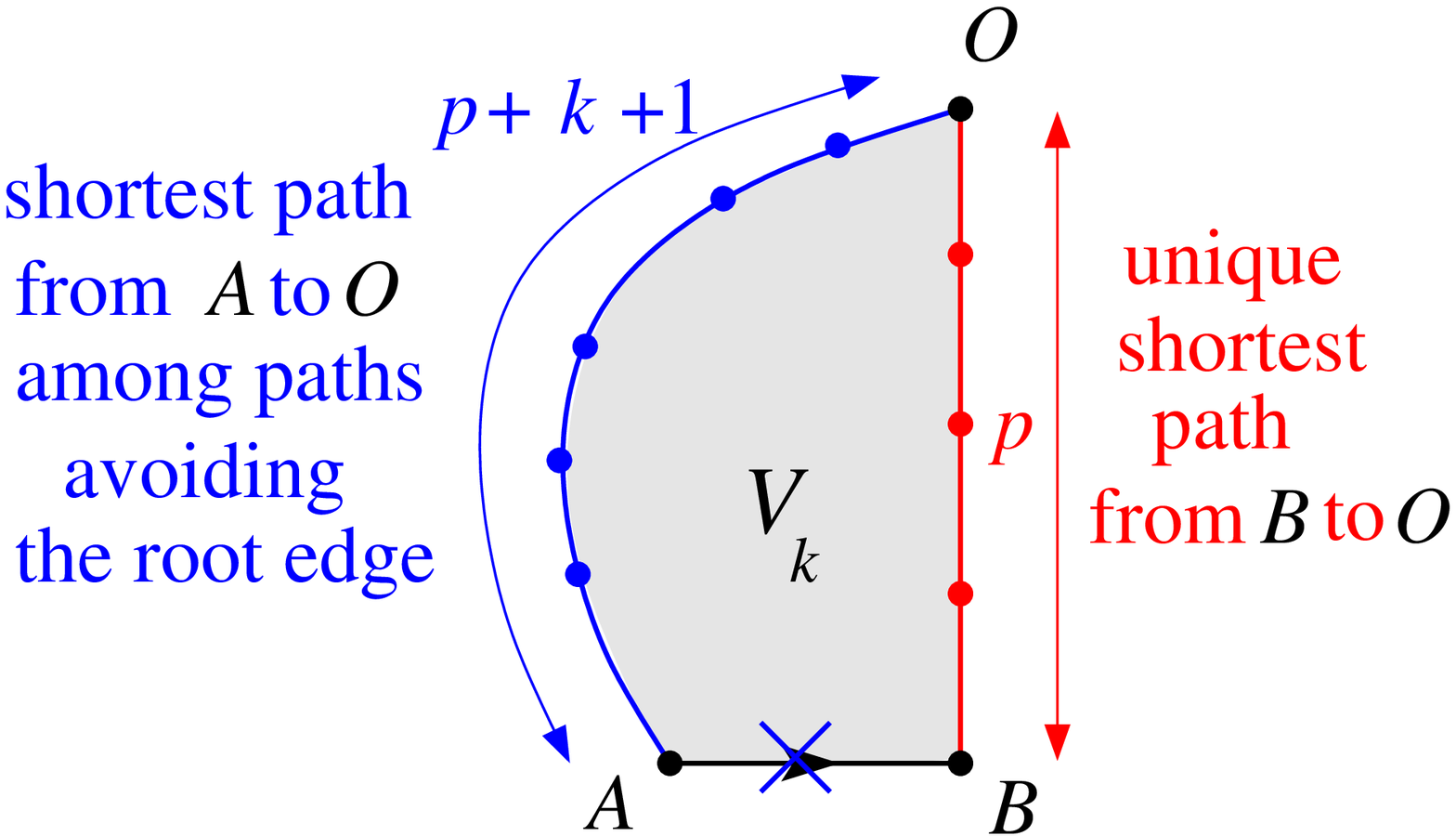}{8.cm}
\figlabel\kslices
In this paper,  we extend this notion of slices to what we shall call {\it $k$-slices}, defined as follows: as displayed in Fig.~\kslices, a $k$-slice is, for a given integer $k$,
a rooted map with a marked vertex $O$ (the apex) incident to its outer face, satisfying the following requirements:
\item{-} the right boundary of the map is the unique shortest path in the map between the endpoint of the root edge and $O$;
\item{-} the left boundary of the map is a shortest path in the map among all paths which join the origin of the root edge to $O$ and 
do not pass via the root edge;
\item{-} the difference of the lengths of the left and right boundaries  is $k+1$. More precisely, if the right boundary has length $p$ for some 
$p\geq 0$, the left boundary then has length $p+k+1$ and we say that the $k$-slice is of type $p/p+k+1$;
\item{-} the vertex $O$ is the only vertex common to both the right and left boundaries;
\item{-} the slice has at least one inner face.
\par\noindent
The interest of this definition lies in the following:

\item{} {\bf Lemma:} For all $k \geq -1$, the generating function of
$k$-slices, counted with weights $x_i$ per inner face of degree $i$,
is equal to the quantity $V_k$ defined by $\Vjdef$.

\noindent {\it Remark:} the statement does not hold for $k<-1$, for
which there are no $k$-slices (if the left boundary has length
strictly less than $p$, the right boundary, of length $p$, cannot be
the unique shortest path between the endpoint of the root edge and
$O$), while $V_k$ is nonzero.

\noindent {\it Proof:} It is easily seen that the lemma holds in the
cases $k=-1$ and $k=0$: the second requirement for a $k$-slice is then
equivalent to demanding that the left boundary itself be a shortest
path among all paths in the map between the origin of the root edge
and $O$, since any path between these two points passing via the root
edge has length larger than $p+1$.  We immediately deduce that the
notion of slices of type $p/p$ and $p/p+1$ introduced above matches
precisely that of $-1$-slices of type $p/p$ and $0$-slices of type
$p/p+1$ respectively (with the slight discrepancy that, due to the
last requirement that a $k$-slice has at least one inner face, the
slice of type $0/1$ reduced to the root edge is not considered as a
$0$-slice). We conclude that the generating function of $-1$-slices is
equal to $S$ and that of $0$-slices is equal to $R-1$. By \eqsforRS,
these are equal to $V_{-1}$ and $V_0$ respectively.
Actually, that $R$ and $S$ satisfy \eqsforRS\ may itself be proved via 
via slices [\xref\PMCF, Appendix A], and we now adapt the argument
to the case of an arbitrary $k\geq -1$.

\fig{Schematic picture of the slice decomposition of a $k$-slice by cutting along the leftmost 
shortest paths from the vertices incident to the face on the left of the root
edge (of degree $j=6$ here) to the apex $O$. Recording the lengths of these paths
creates a three-step path of length $j-1$ with height decrease $k+1$ (bottom right)
with $-1$-slices (resp.\ $0$-slices) associated with level-steps (resp.\ down-steps).
Since we represent the slice with $O$ on top, the three-step path appears vertically reflected in 
the slice. }{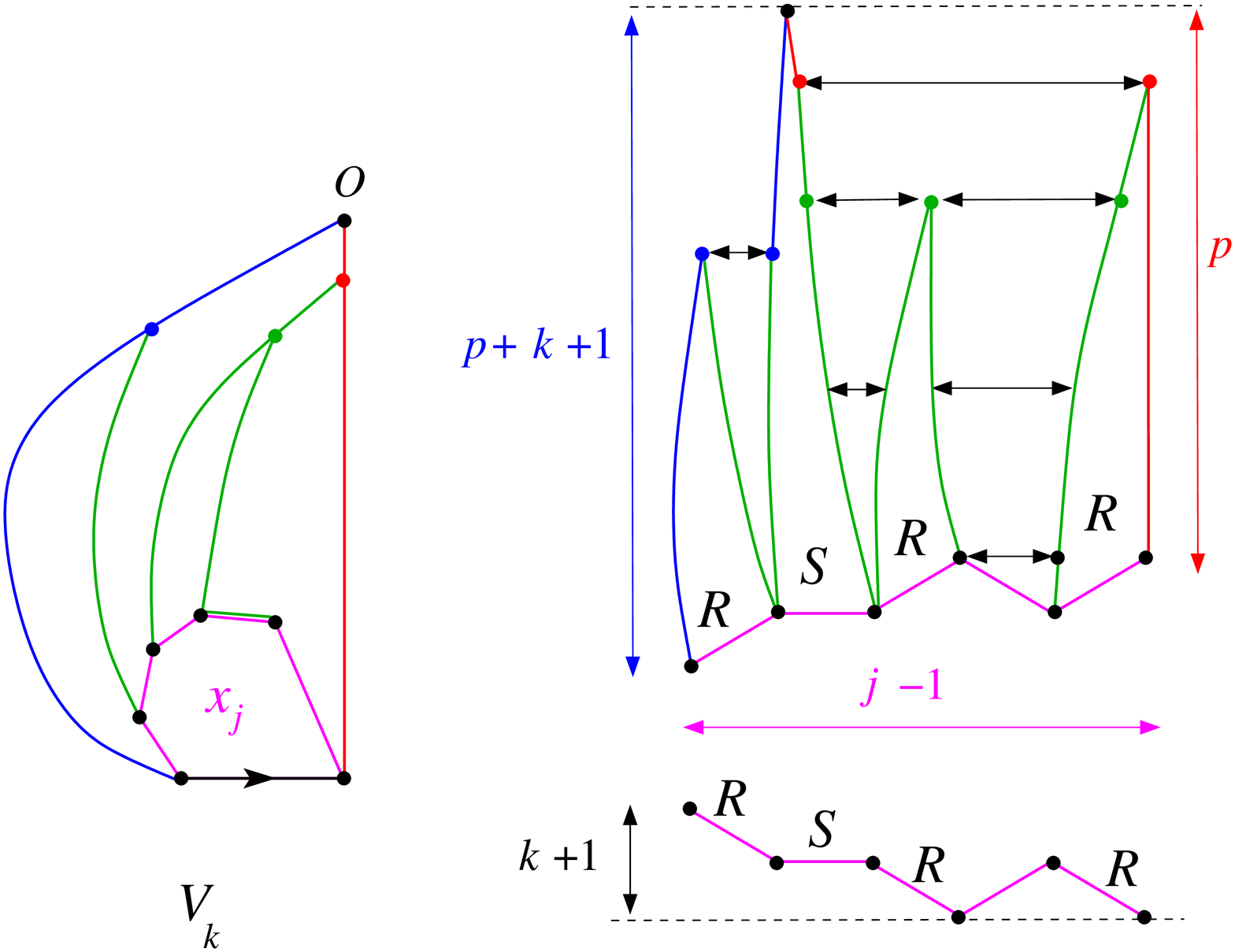}{11.cm}
\figlabel\slicedecomp

Starting from a $k$-slice ($k \geq -1$),
we consider the face to the left of the
root edge: this face has degree $j\geq k+2$ as otherwise, a path of
length $j-1+p$ (obtained by going around the face at hand), hence
strictly smaller than $p+k+1$ and avoiding the root edge, would
join the origin of the root edge to the apex. 
Considering the successive vertices clockwise around this face
and recording their distance to the apex in the $k$-slice with the
root edge removed
creates a three-step path of length $j-1$ starting at height $p+k+1$ ($=$
distance from the origin of the root edge to the apex) and ending at
height $p$ ($=$ distance from the endpoint of the root edge to the
apex), see Fig.~\slicedecomp. We shift all heights
down by $-p-k-1$ so as to obtain a three-step path starting at height $0$ and
ending at height $-k-1$, as counted by $P_{-k-1}(j-1;r,s)$.  Now, for
each vertex around the face, we draw the leftmost shortest path from
it to the apex. Cutting along all these shortest paths results into a
decomposition of the map into connected components, each component
being either a $0$-slice attached to a down-step, since the lengths of
the boundaries differ by $1$ in this case, or a $-1$-slice attached to
a level-step, since the length of the boundaries are then
identical. Note that some steps do not give rise to a slice: up-steps
never do as the leftmost shortest path begins by following the
boundary of the face counterclockwise (see Fig.~\slicedecomp), while
for some down-steps it might occur that a leftmost shortest path
follows the boundary of the face clockwise (in contrast, every
level-step gives rise to a nonempty $-1$-slice).  It is not difficult
to check that, conversely, gluing a sequence of $0$- and $-1$-slices
attached respectively to some down- and all level-steps of a path of
length $j-1$ with total height decrease of $k+1$, and closing the path
by adding an extra root edge, thus creating a new face of degree $j$,
rebuilds a $k$-slice.
Translating this construction in the language of
generating function, it follows that the generating function of
$k$-slices ($k \geq -1$) is equal to $\sum_{j \geq k+2} x_j
P_{-k-1}(j-1;R,S)=V_k$, Q.E.D.


\subsec{$d$-irreducible slices}

\item{} {\bf Lemma:} for all $k \geq 1$, the generating function of
$d$-irreducible $k$-slices, counted with a weight $z$ per inner face
of degree $d$ and, for all $i \geq d+1$, a weight $x_i$ per inner face
of degree $i$, is equal to the quantity $V_k^{(d)}$ defined by
\RSVjdexpr.

\noindent {\it Proof:} let ${\tilde
V}^{(d)}_k(z;x_{d+1},x_{d+2},\ldots)$ be the generating function of
$d$-irreducible $k$-slices. Let us show that the basic substitution
relation \substF\ may be adapted to $k$-slices, namely that we have
\eqn\substVk{{\tilde V}_k^{(d-1)}(0;x_d,x_{d+1},\ldots)={\tilde V}_k^{(d)}(G_d(x_d,x_{d+1},\ldots);x_{d+1},\ldots)}
with $G_d$ as in \Gd.  Indeed, a $k$-slice of type $p/p+k+1$ is a
nothing but a particular instance of rooted map with outer face of
degree $n=2p+k+2$ to which we may apply our substitution approach,
keeping the apex unchanged under substitution.  By following the same
line of arguments as in Sect.~\basicsub, we relate $k$-slices of
type $p/p+k+1$ and girth at least $d$ (l.h.s) to $d$-irreducible
$k$-slices of the same type (r.h.s).  That the type $p/p+k+1$ of slice
remains the same is clear since the substitution does not alter the
lengths of the right and left boundaries of the map.  The only
nontrivial property which must be verified is that these boundaries
remain shortest paths within the desired path sets. Here we use again
the no-shortcut lemma: starting from a $d$-irreducible $k$-slice of
type $p/p+k+1$, we easily deduce from the lemma that substituting
faces of degree $d$ by rooted maps with girth $d$ and outer degree $d$
does not alter the distances between the originally existing vertices
(we apply the same idea as in the proof of the encircling lemma: given
a simple oriented path between two such vertices, we may construct a
path in the original map lying to its left, having the same endpoints,
and having a smaller or equal length). It follows that the right
boundary remains the unique shortest path (of the same length) between its
endpoints after substitution (indeed,
having another path of smaller or equal length after substitution would imply 
having to the left of this path another path of smaller or equal length
which existed before substitution, a contradiction). As for
the left boundary, it remains a shortest path (of the same length) 
between its endpoints among all paths which avoid the root edge
(indeed the existence of a strictly shorter path avoiding the root edge
after substitution would imply the existence to the left of this path -- hence
also avoiding the root edge -- 
of a path with even shorter length already present in the original map,
a contradiction).
Conversely, erasing the
interior of the outermost cycles of length $d$ in a $k$-slice of type
$p/p+k+1$ with girth at least $d$ does not modify the distances
between the remaining vertices so the conditions on the two boundaries
remain satisfied and these boundaries keep the same lengths.

Taking $x_d=X_d(z;x_{d+1},\ldots)$ as in \substVk\ and iterating, we obtain that
\eqn\substVkrecip{\eqalign{{\tilde V}_k^{(d)}(z;x_{d+1},\ldots) &= {\tilde
V}_k^{(d-1)}(0;X_d,x_{d+1},\ldots) = \cdots =
{\tilde V}_k^{(0)}(0;X_1^{(d)},\ldots,X_d^{(d)},x_{d+1},\ldots)}}
where $X_k^{(d)}(z;x_{d+1},\ldots)$, $1 \leq k \leq d$, is defined as
in \Xkddef. But ${\tilde V}_k^{(0)}(0;x_1,x_2,\ldots)$ is the generating
function of $0$-irreducible, i.e.\ arbitrary $k$-slices which, from Sect.~4.1,
is nothing but $V_k(x_1,x_2,\ldots)$. Thus, by \RSVjdexpr, ${\tilde V}_k^{(d)}$ and
$V_k^{(d)}$ coincide, Q.E.D.

To conclude this section, let us discuss the simplifications arising in the case of bipartite maps, i.e.\ when $d$ is even and $x_i=0$ for 
odd $i$.
Then the degree $2p+k+2$ of the outer face of a $k$-slice of type $p/p+k+1$ is necessarily even, which implies that $V_{k}^{(d)}=0$ for odd $k$,
and in particular $S^{(d)}=V_{-1}^{(d)}=0$.
We then have
\eqn\Uk{V_{2k}^{(d)}=U_{k}^{(d)}}
where $U_k^{(d)}$ is the quantity introduced in Sect.~\elimbip.

\newsec{Recursive decomposition of $d$-irreducible slices}

In this section, we shall show that, as generating functions of $d$-irreducible $k$-slices, the quantities $V_k^{(d)}$ for $k<d-1$ satisfy a system of nonlinear equations 
which express a recursive decomposition of the $k$-slices into smaller components, themselves $m$-slices of some kind.
This provides an alternative route for computing them without recourse to the substitution procedure. 
Here, we shall assume that 
\eqn\asum{d\geq 1\qquad \hbox{and}\qquad  -1\leq k \leq d-2.}
Then, since the outer face of a $k$-slice of type $p/p+k+1$ has degree $2(p+1)+k$, the requirement 
of girth at least $d$ implies that $p\geq (d-2-k)/2$. In particular the value $p=0$ may be realized in the range \asum\ that we consider only when $k=d-2$. 
In this case the condition of $d$-irreducibility implies that the only $d$-irreducible $(d-2)$-slice of type $0/d-1$ has a single inner face of degree $d$. 
We now describe how to decompose slices of type $p/p+k+1$ with $p \geq 1$.

\subsec{The binary decomposition procedure}
\subseclab\bindecompproc

\fig{Binary decomposition of a $d$-irreducible $k$-slice of type $p/p+k+1$ with root edge $AB$, for $p \geq 1$: we cut along the leftmost
shortest path from $B$ to the apex $O$ which avoids the first edge $BC$ of the right boundary. This path may
either merge first with the left boundary (left) or with the right boundary (right).}{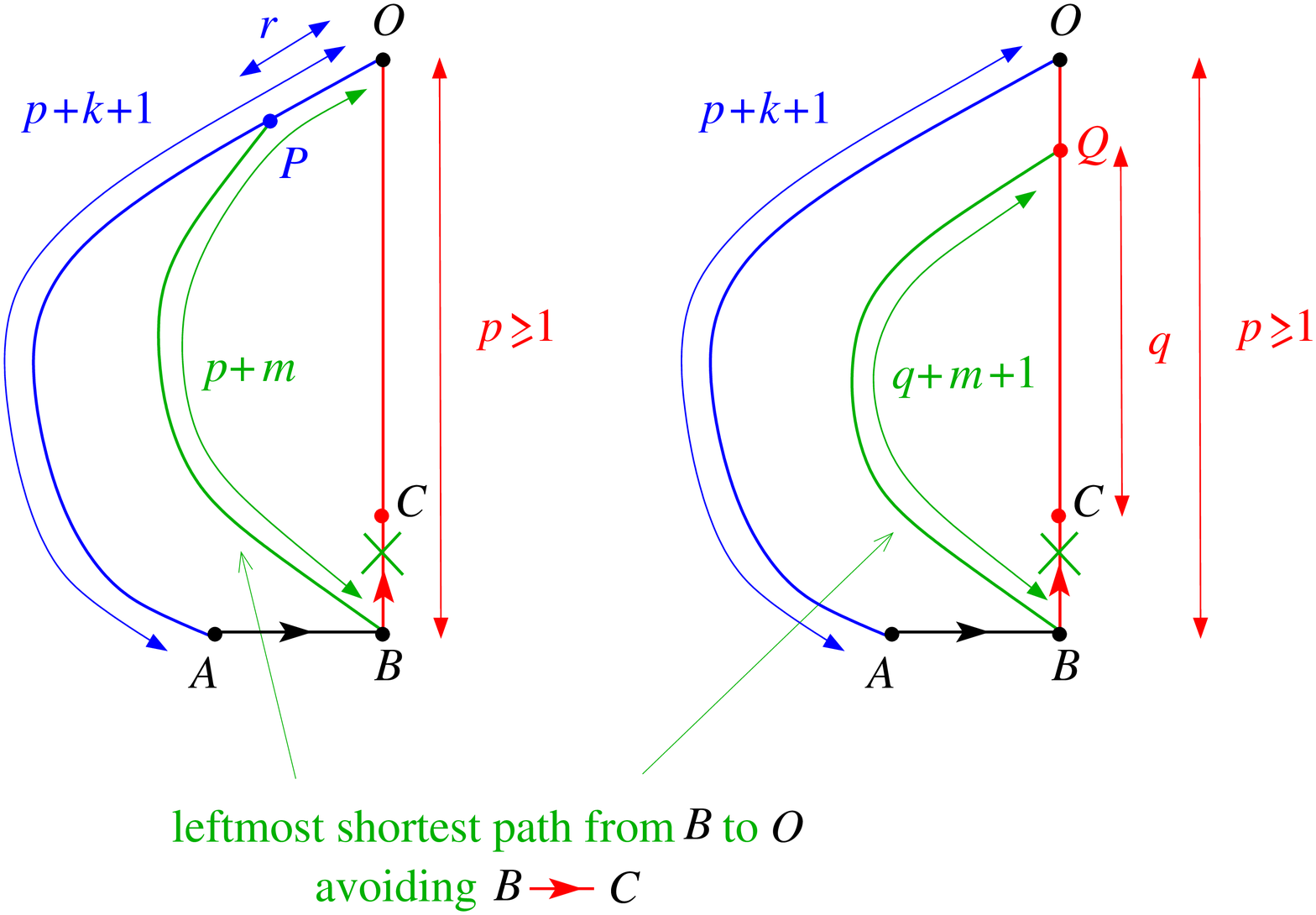}{12.cm}
\figlabel\slicerecur

Let us consider a $d$-irreducible slice of type $p/p+k+1$ with $p \geq
1$, and let $AB$ be its root edge, see Fig.~\slicerecur. We then
singularize the first edge $BC$ of the right boundary, which serves as
root for a new map obtained by cutting the $k$-slice along the
leftmost shortest path from $B$ to the apex $O$ which avoids the edge
$BC$. This {\it marked path} has a length $p+m$, with $m\geq 1$, since
the right boundary is the unique shortest path from $B$ to $O$, and
$m\leq k+2$, since the left boundary gives rise to a path of length
$p+k+2$ between $B$ and $O$ and avoiding $BC$. Clearly the part of the
map lying in-between the marked path and the right boundary is a
$d$-irreducible $m$-slice (with root edge $BC$). This $m$-slice may be
of type $q/q+m+1$ for any $q\leq p-1$ as the marked path may hit the
right boundary at a point $Q$ between $C$ and $O$ (note that the
marked path ``sticks'' to the boundary after hitting it).  If $m=k+2$,
the marked path starts with $BA$ and necessarily coincides with the left boundary of the
$k$-slice, thus the decomposition amounts to rerooting the slice on
$BC$, changing its type to $p-1/p+k+2$.
If $m\leq k+1$, the part lying in-between the marked path and the left
boundary is a $d$-irreducible $(k-m)$-slice. It is of type
$q'/q'+k-m+1$, where $q'$ is the length of the portion of marked path
from $B$ to the point $P$ where it hits the left boundary (thus
$q'=p+m-r$ where $r\geq 0$ is the distance $PO$, see again
Fig.~\slicerecur).
Since the marked path either hits the right or the left boundary, we 
have either $q=p-1$ or $q'=p+m$, i.e.\ $p=\max(q+1,q'-m)$.

To summarize, a $d$-irreducible slice of type $p/p+k+1$ with $p \geq
1$ is decomposed into either a $d$-irreducible slice of type
$p-1/p+k+2$ or into a pair of $d$-irreducible slices, one of type
$q/q+m+1$ and one of type $q'/q'+k-m+1$, where $m$ is an integer
between $1$ and $k+1$ and $q,q'$ are two nonnegative integers such
that $p=\max(q+1,q'-m)$.

\fig{The gluing of a $d$-irreducible $m$-slice of type $q/q+m+1$ with a $d$-irreducible $(k-m)$-slice of
type $q'/q'+k-m+1$ yields a $d$-irreducible $k$ slice of type $p/p+k+1$ with $p=\max(q'-m,q+1)$.}{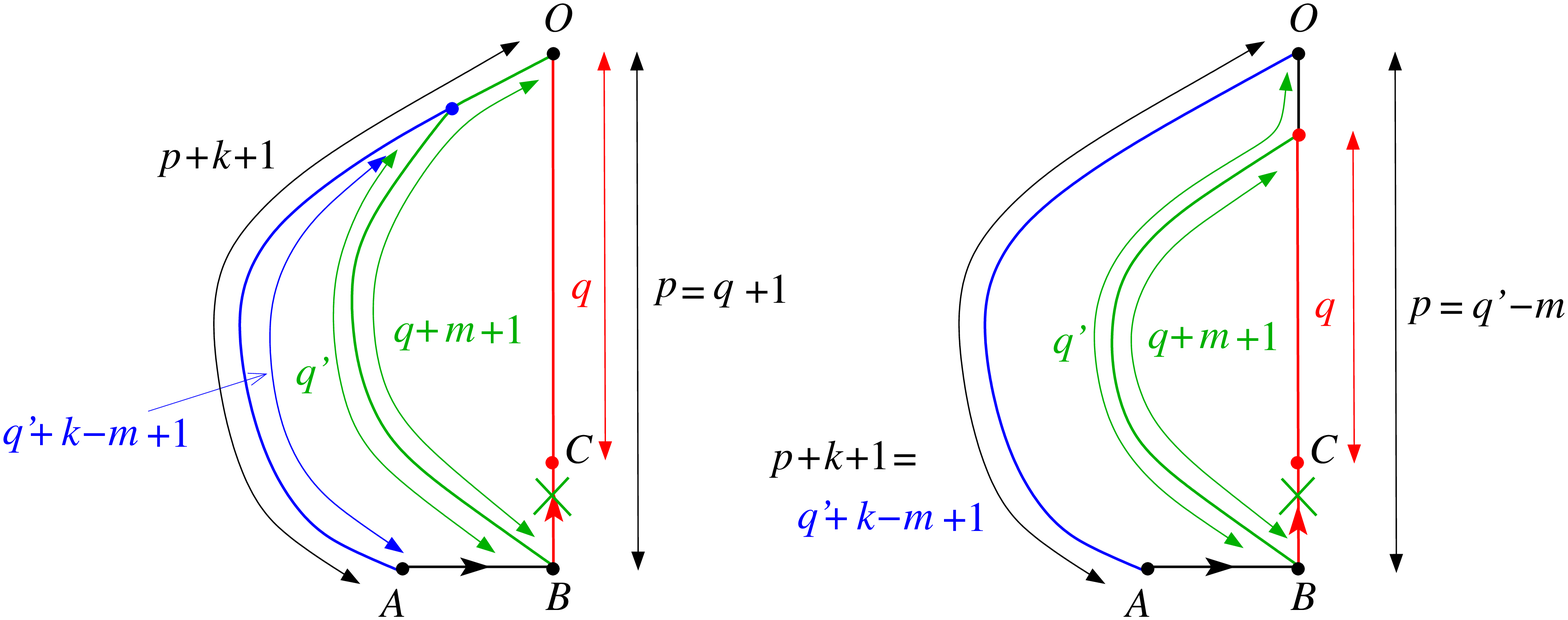}{14.cm}
\figlabel\slicerecurback

Conversely, gluing as in Fig.~\slicerecurback\ an arbitrary  $d$-irreducible $m$-slice of type $q/q+m+1$, $1\leq m\leq k+1$, with root edge $BC$, and an arbitrary
 $d$-irreducible $(k-m)$-slice of type $q'/q'+k-m+1$, with root-edge $AB$, creates a 
rooted map of type $p/p+k+1$ where $p={\rm max}(q'-m,q+1)>0$. Here the apex $O$ of the concatenated map is chosen as being the apex of the $m$-slice if $q'\leq q+m+1$
and that of the $(k-m)$-slice otherwise. This map is clearly $d$-irreducible: indeed, in a $d$-irreducible slice, any simple path joining two vertices of the left (resp.\ right) boundary and not entirely included in the boundary 
has length at least (resp.\ strictly larger than) $d/2$ (as otherwise closing this path with the, necessarily shorter, portion of the boundary between the two vertices would create a cycle of length
 strictly less than $d$), which ensures that the concatenation of $d$-irreducible slices along their boundaries is still $d$-irreducible.
 Moreover, since $m\geq 1$, the right boundary of the resulting map is clearly the unique shortest path in the map from the point $B$ to the apex.  
 As for the left boundary of the 
 concatenated map, it is clearly a shortest path from $A$ to $O$ among all paths which avoid both the $AB$ and $BC$ edges. To make sure we have a $k$-slice, we must
 guarantee that this is also true among paths which avoid $AB$ but pass via $BC$, which requires that the length of any path from $A$ to $B$ avoiding
 $AB$ has length at least $k+1$. This is again guaranteed by the condition of girth at least $d$: any path from $A$ to $B$ 
 avoiding $AB$ has length at least $d-1$ which is
 larger than or equal to $k+1$ in the range \asum. This shows that our decomposition is a bijection (in the case 
 of a slice of type $p-1/p+k+2$, its unique pre-image of type $p/p+k+1$ is recovered by a simple rerooting).

\subsec{The iterated decomposition procedure}
\subseclab\itdecompproc

\fig{Iterated decomposition of a $d$-irreducible $k$-slice with root $AB$ and apex $O$. 
We start with the binary decomposition of Fig.~\slicerecur, i.e.\ pick the leftmost shortest path $P_1$ 
from $B$ to $O$ avoiding the first edge $BC_1$ of the right boundary. We then iterate the procedure by picking the leftmost shortest path $P_2$  from $B$ to $O$
staying to the left of $P_1$ 
and avoiding the first edge $BC_2$ of $P_1$. The process ends  after a finite number of iterations  whenever the selected leftmost shortest path (here $P_3$)
matches the left boundary of the slice.}{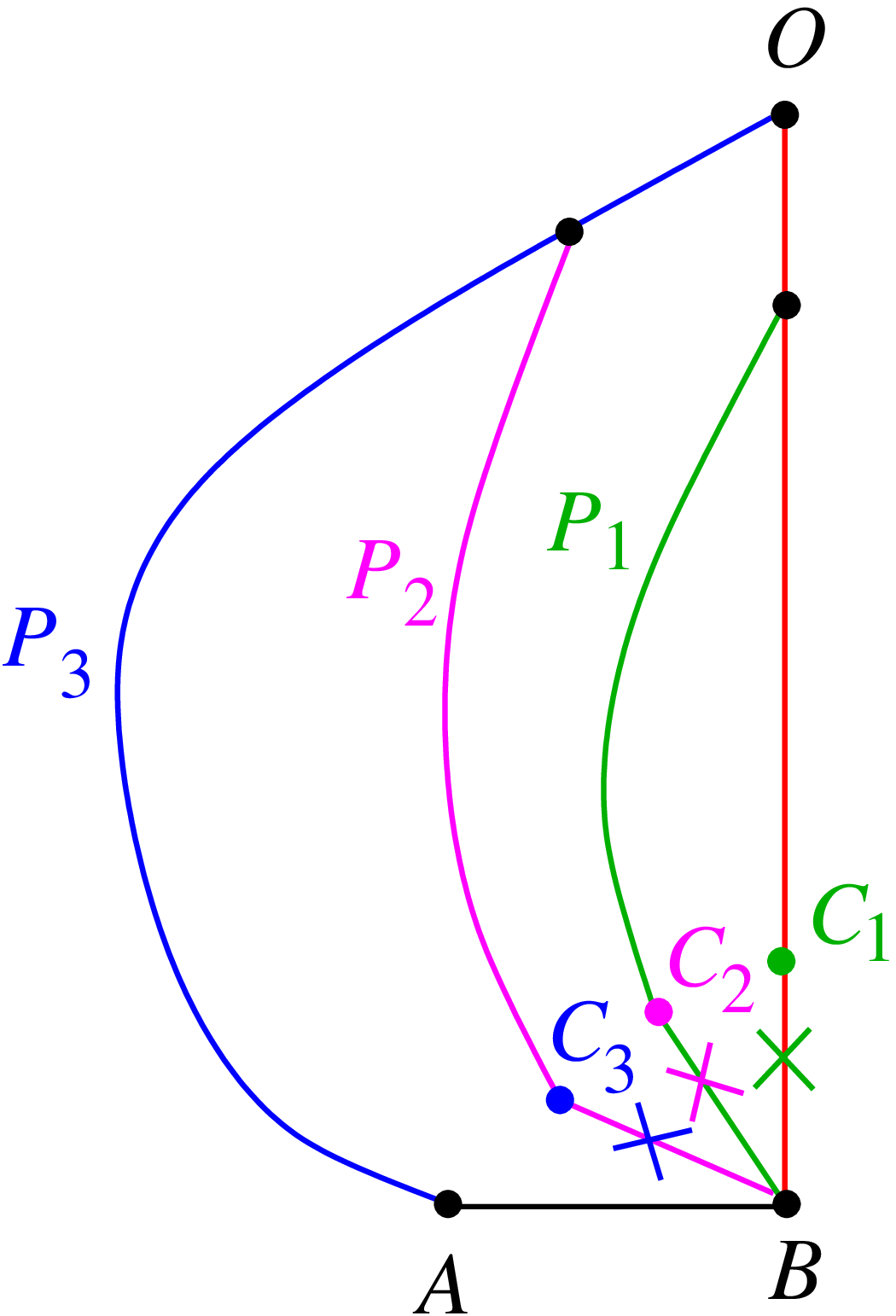}{4.5cm}
\figlabel\iterateddecomp

It is useful to consider a variant of the above binary decomposition
procedure (we dub it ``binary'' since it splits a slice into at most
two parts). Simply put, it consists in iterating the decomposition on
all subslices with the same root edge as the original slice. More
precisely, starting again from a slice of type $p/p+k+1$ ($p \geq
1$) with root edge $AB$, we decompose it as follows:
\item{-} as in the binary decomposition, we pick the leftmost shortest
path $P_1$ among all paths from $B$ to the apex $O$ which do not pass
via the first edge $BC_1$ of the right boundary.  Since the right
boundary is the only shortest path between $B$ and $O$, the length of
$P_1$ is $p+m_1$ for some $m_1\geq 1$.  We call $C_2$ the extremity of
the first edge of $P_1$ (see Fig.~\iterateddecomp).
\item{-} we then iterate by picking the leftmost shortest path $P_2$
among all paths from $B$ to the apex $O$ which stay to the left of
$P_1$ and do not pass via the edge $BC_2$. Since $P_1$ is a leftmost
shortest path between $B$ and $O$, the length of $P_2$ is $p+m_1+m_2$
for some $m_2\geq 1$. We call $C_3$ the extremity of the first edge of
$P_2$.
\item{-} we continue the process until the sum $\sum_i m_i$, which
increases strictly at each step, reaches $k+2$. The corresponding
leftmost path $P_q$, of length $p+k+2$, then necessarily follows the
left boundary of the slice. In other words, the extremity $C_{q+1}$ of
the first edge of $P_q$ coincides with the origin $A$ of the root
edge.
\par
\noindent Cutting along all the paths $P_j$, $j=1,\cdots,q-1$, results
in a decomposition of the slice into $q$ pieces. Upon rooting the
$i$-th piece on the edge $BC_i$, this piece is clearly, by
construction, a $d$-irreducible $m_i$-slice (the lengths of its right
and left boundaries differ by $m_i+1$).  Note again that, if $q=1$,
i.e. $m_1=k+2$, the slice is kept uncut but rerooted at the edge
$BC_1$ so that we obtain a $(k+2)$-slice whose right boundary length is
reduced by $1$ and whose left boundary length is increased by $1$.
To summarize, the iterated decomposition transforms a slice of type
$p/p+k+1$ for any $p \geq 1$ into a sequence of $q$ slices, the $i$-th
slice being a $m_i$-slice, with $m_i \geq 1$ and $m_1 + \cdots + m_q =
k+2$.

If we now define the {\it size} of a slice as being the number of its
{\it non-left edges}, i.e.\ those edges which do not belong to the
left boundary, then it is clear that, in the iterated decomposition,
the sum of the sizes of all the subslices is exactly $1$ less than the
size of the original slice. This is because all non-left edges of the
original slice are non-left edges of exactly one subslice, except the
root edge.

 \subsec{Recursive equations for the generating functions}

We now translate the above decompositions into equations. At this
stage, we are interested in generating functions of $d$-irreducible
$k$-slices, irrespectively of their precise type (i.e. we disregard
the value of $p$ for now). Then, the binary decomposition
states that a $k$-slice (distinct from the trivial $(d-2)$-slice
reduced to a single $d$-valent face) is in one-to-one correspondence
with either a $(k+2)$-slice or an ordered pair formed by a $m$-slice
and $(k-m)$-slice, with $1 \leq m \leq k+1$. This yields immediately
\eqn\recurV{V_k^{(d)}=z \delta_{k,d-2}+\sum_{m=1}^{k+1}V_m^{(d)}
V_{k-m}^{(d)}+ V_{k+2}^{(d)}\ , \qquad -1\leq k\leq d-2}
where the term $z \delta_{k,d-2}$ is the contribution from the only slice
of type $0/k+1$ which, as noted above, is reduced to a single
$d$-valent face, with $k=d-2$. 
As for the iterated decomposition, it yields
the equivalent form
\eqn\recurVbis{V_k^{(d)}=z \delta_{k,d-2}+\sum_{q\geq 1} \sum_{
m_i \geq 1,\  i=1,\cdots ,q \atop m_1+\cdots+m_q=k+2 } \prod_{i=1}^{q} V_{m_i}^{(d)}, \qquad -1\leq k\leq d-2.}
Both systems of equations determine $V_k^{(d)}$ for all $-1\leq k\leq d-2$
from the data of $V_{d-1}^{(d)}$ and $V_d^{(d)}$. These quantities
are themselves obtained from \Vjdbig, namely
\eqn\Vkdminusoned{\eqalign{V_{d-1}^{(d)}&=\sum_{j\geq d+1} x_j P_{-d}(j-1;R^{(d)},S^{(d)})\cr
V_{d}^{(d)}&=\sum_{j\geq d+2} x_j P_{-d-1}(j-1;R^{(d)},S^{(d)}).\cr}
} 
This allows to determine $V_k^{(d)}$ for all $-1\leq k\leq d$  in terms of $z$, the $x_i$ ($i\geq d+1$), $R^{(d)}$ and $S^{(d)}$. Equating
the obtained expressions for $V_{-1}^{(d)}$ and $V_0^{(d)}$ to the values
\eqn\SRVd{V_{-1}^{(d)}=S^{(d)}, \qquad V_0^{(d)}=R^{(d)}-1,}
obtained by specializing \eqsforRS, we get algebraic equations which determine $R^{(d)}$ and $S^{(d)}$ themselves.

It is interesting to note that Eq.~\recurV\ can be extended to the value $k=d-1$ in the case $z=0$, i.e.\ when we consider maps of girth
at least $(d+1)$ but not necessarily $(d+1)$-irreducible. Examining the decomposition of $(d-1)$-slices, we arrive at
\eqn\recurVdminusone{\eqalign{V_{d-1}^{(d)}(0;x_{d+1},\ldots)& =G_{d+1}(x_{d+1},\ldots)+\sum_{m=1}^{d}V_m^{(d)}(0;x_{d+1},\ldots) V_{d-1-m}^{(d)}(0;x_{d+1},\ldots)\cr & \qquad \qquad \qquad \qquad \qquad \qquad \qquad \qquad \qquad \qquad+
V_{d+1}^{(d)}(0;x_{d+1},\ldots)\cr}}
with $V_{d+1}^{(d)}(0;x_{d+1},\ldots)=\sum_{j\geq d+3} x_j P_{-k-1}(j-1;R^{(d)}(0;x_{d+1},\ldots),S^{(d)}(0;x_{d+1},\ldots))$. Here we simply used the fact that the $p=0$ contribution to $V_{d-1}^{(d)}(0;x_{d+1},\ldots)$ is precisely the generating function $G_{d+1}(x_{d+1},\ldots)$  of maps of outer degree $d+1$ and of girth  $d+1$
(hence the first term on the r.h.s.). Eq.~\recurVdminusone\ allows eventually to obtain $G_{d+1}(x_{d+1},\ldots)$ and recover the result of \BFb.
\bigskip

As an illustration, let us consider the case of  $5$-angular irreducible dissections ($d=5$ and $x_i=0$ for all $i\geq 6$).
Using, as before, lower case letters for generating functions specialized at $x_i=0$, Eq.~\recurV\ gives
\eqn\recurfive{\eqalign{
v_{-1}^{(5)} & = v_1^{(5)} \cr
v_0^{(5)} & = v_1^{(5)}v_{-1}^{(5)}+v_2^{(5)} \cr
v_1^{(5)} & = v_1^{(5)}v_0^{(5)}+v_2^{(5)} v_{-1}^{(5)}+v_3^{(5)}\cr
v_2^{(5)} & = v_1^{(5)}v_1^{(5)}+v_2^{(5)}v_0^{(5)}+v_3 ^{(5)}v_{-1}^{(5)}+v_4^{(5)}\cr
v_3^{(5)}& =z+ v_1^{(5)}v_2^{(5)}+v_2^{(5)}v_1^{(5)}+v_3^{(5)}v_0^{(5)}+v_4^{(5)} v_{-1}^{(5)}+v_5^{(5)}\cr}}
with $v_4^{(5)}=v_5^{(5)}=0$ from \Vkdminusoned\ and, from \SRVd, $v_{-1}^{(5)}=s^{(5)}$ and $v_0^{(5)}=r^{(5)}-1$.
This leads, by elimination, to the system of equations 
\eqn\SRfivebis{\eqalign{0&= 2+(r^{(5)})^2+6 (s^{(5)})^2 +(s^{(5)})^4 -3 r^{(5)}(1+(s^{(5)})^2)\cr
z &= s^{(5)}\left(2(r^{(5)})^2+4(2+(s^{(5)})^2)-r^{(5)}(9+(s^{(5)})^2)\right)\cr}}
which is equivalent to \SRfive.
    
\subsec{Solving the recursive equations in the bipartite case}

In the bipartite case, setting $d=2b$, Eq.~\recurV\ translates into
\eqn\recurU{U_k^{(d)}=z \delta_{k,b-1}+\sum_{m=1}^{k}U_m^{(d)} U_{k-m}^{(d)}+
U_{k+1}^{(d)}\ , \qquad 0\leq k\leq b-1}
which determines $U_k^{(d)}$ for all $0\leq k\leq b-1$  from the data of $U_b^{(d)}$. The latter is itself determined from \Ujbig\ by
\eqn\Ukb{U_{b}^{(d)}=\sum_{j\geq b+1} x_{2j} {2j-1 \choose j+b}(R^{(d)})^{j+b}} 
with moreover the relation $U^{(d)}_0=R^{(d)}-1$. By a simple extension of the argument leading to \recurVdminusone, we find also in the special case  $z=0$
\eqn\recurUb{\eqalign{U_{b}^{(d)}(0;x_{d+2},\ldots)& =G_{d+2}(x_{d+2},\ldots)+\sum_{m=1}^{b}U_m^{(d)}(0;x_{d+2},\ldots) U_{b-m}^{(d)}(0;x_{d+2},\ldots)\cr & \qquad \qquad \qquad \qquad \qquad \qquad \qquad \qquad \qquad \qquad+
U_{b+1}^{(d)}(0;x_{d+2},\ldots)\cr}}
with $U_{b+1}^{(d)}=\sum_{j\geq b+2} x_{2j} {2j-1 \choose j+b+1}(R^{(d)})^{j+b+1}$. This allows to deduce the generating function $G_{d+2}(x_{d+2},\ldots)$
of maps of outer degree $d+2$ and girth $d+2$ \BFb.
 
It turns out that it is possible to obtain a single 
algebraic equation for $R^{(d)}$ as follows: we observe that the system \recurU\ is triangular in the sense that the $(k+1)$-th 
(with $k\geq 0$) equation allows to express $U_{k+1}^{(d)}$ in terms of the $U_\ell^{(d)}$ with $\ell \leq k$, therefore in terms of $U_{0}^{(d)}$.
Therefore, we may consider the semi-infinite triangular system of equations 
\eqn\recurUtilde{{\tilde U}_k= \sum_{m=1}^{k}{\tilde U}_m {\tilde U}_{k-m}+
{\tilde U}_{k+1}, \quad k\geq 0}
determining the ${\tilde U}_k$, $k>0$, in terms of ${\tilde U}_0$ and then \recurU\ is recovered by identifying
${\tilde U}_0=U^{(d)}_0=R^{(d)}-1$ and taking formally ${\tilde U}_b=z+U^{(d)}_b$. Introducing the generating function
\eqn\gfUtilde{{\tilde U}(t)=\sum_{k\geq 1} {\tilde U}_k t^k\ ,}
Eq.~\recurUtilde\ translates immediately into
\eqn\eqUtilde{{\tilde U}_0+{\tilde U}= {\tilde U}\left({\tilde U}_0+{\tilde U}+{1\over t}\right)}
hence 
\eqn\eqt{t= {{\tilde U}\over (1-{\tilde U})({\tilde U}_0+{\tilde U})}\ .}
The Lagrange inversion formula [\xref\FS-\xref\GouldJack] states that, for $k>0$,
\eqn\restildeUk{{\tilde U}_k= {1\over k}[{\tilde U}^{k-1}]\left( (1-{\tilde U})({\tilde U}_0+{\tilde U})\right)^k= -{1\over k}\sum_{p=1}^{k} 
{k\choose p}{k \choose p-1} (-{\tilde U}_0)^p\ .}
In particular, taking $k=b$, the above identification leads to the algebraic equation for $R^{(d)}$:
\eqn\alpeq{ -{1\over b} \sum_{p=1}^{b}  {b\choose p}{b \choose p-1} (1-R^{(d)})^p = z+ \sum_{j\geq b+1} x_{2j} {2j-1 \choose j+b}(R^{(d)})^{j+b}\ .}
We let the reader verify that this coincides with \eqforRdbip\ via yet another hypergeometric identity.
Recall that the generating functions $F_{2n}^{(d)}$ are related to $R^{(d)}$ by the pointing formula \Fndbippoint. 

Let us now give a few examples. Irreducible quadrangular dissections are obtained for $b=2$ and $x_{2j}=0$ for all $j\geq 3$, hence 
Eq.~\alpeq\ reads
\eqn\algeqirrquad{(r^{(4)}-1)-(r^{(4)}-1)^2=z\ ,}
which matches \algebRquad.
The case of irreducible hexangular dissections corresponds to $b=3$ and $x_{2j}=0$ for all $j\geq 4$, hence Eq.~\alpeq\  
now reads
\eqn\algeqirrhex{(r^{(6)}-1)-3(r^{(6)}-1)^2+(r^{(6)}-1)^3=z\ ,}
which matches \eqforrhex.
Finally, the case of $2$-irreducible quadrangulations, i.e.\ quadrangulations without multiple edges corresponds to $b=1$, $z=0$ and 
$x_{2j}=0$ for all $j\geq 3$, in which case Eq.~\alpeq\ gives
\eqn\algeqtwoirrquad{(R^{(2)}-1)=x_4 \left( R^{(2)} \right)^3\ ,} 
in agreement with \DSQ. From \recurUb, we then deduce that the  generating function $G_4(x_4)$ of quadrangulations of outer degree $4$ and girth 4
(i.e.  without multiple edges) is $G_4(x_4)=(R^{(2)}-1)(2-R^{(2)})$. 

\subsec{Solution of the recursive equations in the general case}

Let us now see how to extend to the general (non necessarily bipartite) case the strategy of Sect.~5.4. Our
goal is again to obtain a system of algebraic equations for $R^{(d)}$ and $S^{(d)}$ without recourse to substitution,
i.e.\ using the system  \recurV\ as starting point. Again we observe that this system is triangular, with its $(k+2)$-th 
(with $k\geq -1$) equation allowing to express $V_{k+2}^{(d)}$ in terms of the $V_\ell^{(d)}$ with $\ell \leq k$, hence eventually
in terms of $V_{-1}^{(d)}$ and $V_{0}^{(d)}$ only. Again we introduce the semi-infinite triangular system of equations 
\eqn\recurVtilde{{\tilde V}_k= \sum_{m=1}^{k+1}{\tilde V}_m {\tilde V}_{k-m}+
{\tilde V}_{k+2}, \quad k\geq -1}
which determines the ${\tilde V}_k$, $k>0$, in terms of ${\tilde V}_{-1}$ and ${\tilde V}_0$. Then we may obtain algebraic equations
for $R^{(d)}$ and $S^{(d)}$ by simply identifying
${\tilde V}_{-1}=V^{(d)}_{-1}=S^{(d)}$,  ${\tilde V}_{0}=V^{(d)}_{0}=R^{(d)}-1$ and by taking formally ${\tilde V}_{d-1}=V^{(d)}_{d-1}$ and 
${\tilde V}_{d}=z+V^{(d)}_{d}$ so as to fulfill \recurV\ at $k=d-3$ and $k=d-2$.  If we now introduce the generating function
\eqn\gfVtilde{{\tilde V}(t)=\sum_{k\geq 1} {\tilde V}_k t^k\ ,}
Eq.~\recurV\ translates into
\eqn\eqVtilde{{{\tilde V}_{-1}\over t}+{\tilde V}_0+{\tilde V}= {\tilde V}\left({{\tilde V}_{-1}\over t}+{\tilde V}_0+{\tilde V}+{1\over t^2}\right)}
hence 
\eqn\eqtbis{t\, {\tilde V}_{-1}+ t^2 ({\tilde V}_0+{\tilde V})= {{\tilde V}\over 1-{\tilde V}}\ .}
For fixed ${\tilde V}_{-1}$ and ${\tilde V}_{0}$, this determines $t$ as a function of ${\tilde V}$, hence in principle 
all ${\tilde V}_{k}$ for $k>0$ in terms of  ${\tilde V}_{-1}$ and ${\tilde V}_{0}$ via a Lagrange inversion  ${\tilde V}_k= {1\over k}[{\tilde V}^{k-1}]( {\tilde V}/t)^k$.
We have not been able to perform the computation but, by inspection of the first terms, we conjecture that
\eqn\tildeVexpr{\eqalign{&{\tilde V}_{2j-1}=-\sum_{k\geq 0 ,m \geq 0\atop k+m\leq j-1}  {j+m-1\choose k+2m}{j+m\choose k+2m}{ {k+2m \choose 2m}\over 2m+1}(-{\tilde V}_{-1})^{2m+1} (-{\tilde V}_0)^k\cr
&{\tilde V}_{2j}=-\sum_{k\geq 0 ,m \geq 0\atop 1\leq k+m\leq j} {j+m-1\choose k+2m-1}{j+m\choose k+2m-1}{ {k+2m-1 \choose 2m-1}\over 2m}(-{\tilde V}_{-1})^{2m} (-{\tilde V}_0)^k\cr}
}
for $j\geq 1$ (in the second equation, the coefficient ${k+2m-1 \choose 2m-1}/(2m)$ for $m=0$ and $k\geq 1$, should be understood as its $m\to 0$ limit $1/k$). Assuming that these expressions are valid, we deduce, upon making the identifications mentioned above for
${\tilde V}_{d-1}$ and ${\tilde V}_{d}$, the equations 
\eqn\eqsforRSdeven{\eqalign{
\sum_{k\geq 0 ,m \geq 0\atop k+m\leq {d\over 2}-1}  {{d\over 2}+m-1\choose k+2m}{{d\over 2}+m\choose k+2m}& {{k+2m \choose 2m}\over 2m+1}(-S^{(d)})^{2m+1}(1-R^{(d)})^k 
\cr & + \sum_{j \geq d+1} x_j P_{-d}(j-1;R^{(d)},S^{(d)}) =0 \cr
  z + \sum_{k\geq 0 ,m \geq 0\atop 1\leq k+m\leq {d\over 2}}  {{d\over 2}+m-1\choose k+2m-1}{{d\over 2}+m\choose k+2m-1} & { {k+2m-1 \choose 2m-1}\over 2m}(-S^{(d)})^{2m} (1-R^{(d)})^k \cr
  &  + \sum_{j \geq d+2} x_j P_{-d-1}(j-1;R^{(d)},S^{(d)}) =0 \cr }}
if $d$ is even and
\eqn\eqsforRSdodd{\eqalign{
 \sum_{k\geq 0 ,m \geq 0\atop 1\leq k+m\leq {d\!-\!1\over 2}}  {{d\!-\!1\over 2}+m-1\choose k+2m-1}{{d\!-\!1\over 2}+m\choose k+2m-1}&{ {k+2m-1 \choose 2m-1}\over 2m}(-S^{(d)})^{2m} (1-R^{(d)})^k \cr
 & + \sum_{j \geq d+1} x_j P_{-d}(j-1;R^{(d)},S^{(d)}) =0 \cr
 z + 
 \sum_{k\geq 0 ,m \geq 0\atop k+m\leq {d\!+\!1\over 2}-1}  {{d\!+\!1\over 2}+m-1\choose k+2m}{{d\!+\!1\over 2}+m\choose k+2m}&{ {k+2m \choose 2m}\over 2m+1} (-S^{(d)})^{2m+1} (1-R^{(d)})^k  \cr 
 &  + \sum_{j \geq d+2} x_j P_{-d-1}(j-1;R^{(d)},S^{(d)}) =0 \cr }}
if $d$ is odd. These equations determine $R^{(d)}$ and $S^{(d)}$ and are an alternative to the system \eqsforRSd. 

\newsec{Bijection between slices and trees}

In this section, we exhibit a bijection between slices and some trees,
which enjoy the same decomposition structure. For the sake of
simplicity, we restrict ourselves to irreducible $d$-angular
dissections, with $d \geq 3$. There seems to be no conceptual
difficulty in extending the forthcoming discussion to general
$d$-irreducible maps but the corresponding trees would become more
complicated.
 
\subsec{$d$-oriented $k$-trees}

Recall that, by the iterated decomposition of Sect.~\itdecompproc,
the generating functions of irreducible $d$-angular $k$-slices $v_k^{(d)}$
satisfy 
\eqn\recurv{v_k^{(d)}=z \delta_{k,d-2}+\sum_{q\geq 1} \sum_{
1 \leq m_1,\cdots,m_q \leq d-2 \atop m_1+\cdots+m_q=k+2 } \prod_{i=1}^{q} v_{m_i}^{(d)}, \qquad -1\leq k\leq d-2}
as seen by specializing \recurVbis, and noting that $v_k^{(d)}=0$ for
$k \geq d-1$ by \Vjdbig\ at $x_i=0$. Note that we may restrict the
range to $1 \leq k \leq d-2$ and still obtain a closed system for
$v_1^{(d)},\ldots,v_{d-2}^{(d)}$. This system clearly specifies some trees,
which we now describe.

\fig{Rules for the construction of $d$-oriented trees (see text). Black (resp.\ white) dots represent inner or root vertices 
(resp.\ non-root leaves). }{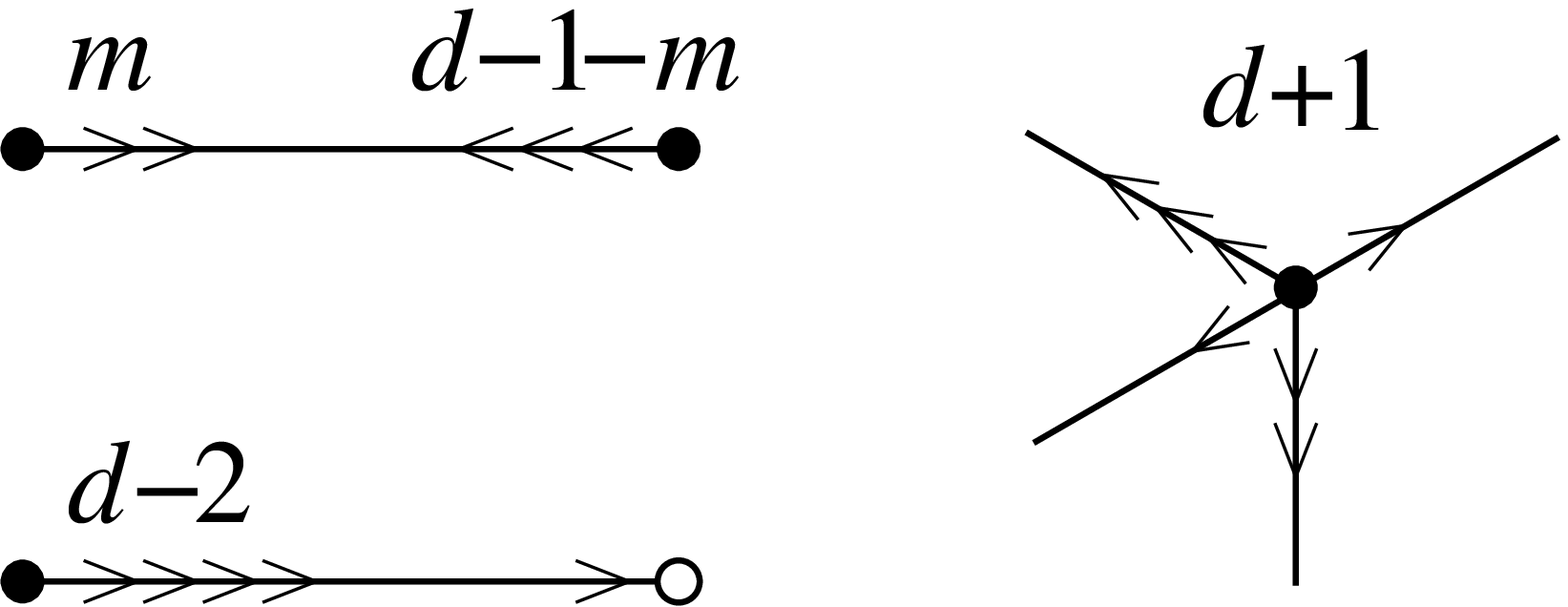}{9.cm}
\figlabel\treerules
For $1\leq k\leq d-2$, we define a {\it $d$-oriented $k$-tree} as a planted 
(with a marked univalent root vertex) plane tree such that:
\item{-} the edges of the tree carry arrows and are of two types (see Fig.~\treerules):
\itemitem{-} inner edges of type $m/(d-1-m)$, with $1\leq m\leq d-2$, whose two half-edges carry respectively $m$ and
$d-1-m$ arrows pointing away from the associated edge extremity. These edges connect only inner vertices or the
root vertex.
\itemitem{-} leaf edges with one half-edge carrying $d-2$ arrows pointing away from the associated edge extremity and
the other half-edge carrying a single arrow pointing toward the associated extremity (this arrow plays no role in the following
but we decided to introduce it so as to recover trees similar to those of \FPS). The first extremity of the
edge is necessarily an inner vertex or the root vertex, and the second extremity is a leaf of the tree.
\item{-} the out-degree of any inner vertex, defined as the total number of arrows pointing away from that vertex among all
the incident half-edges, is $d+1$.
\item{-} the out-degree of the root vertex is $k$.
\par
\noindent 
In a $d$-oriented $k$-tree, the edge emerging from the root vertex can be a leaf edge only if $k=d-2$, resulting in
a tree with a single edge and a single leaf. If it is instead an inner edge, it must be of type $k/(d-1-k)$ linking the
root vertex to some inner vertex. The descending subtrees attached to this vertex form a sequence of a number $q\geq 1$ of trees, 
the $i$-th one being a $d$-oriented $m_i$-tree for some $m_i$ between $1$ and $d-2$. From the out-degree condition on 
inner vertices, we deduce that $\sum_i m_i=d+1-(d-1-k)=k+2$. We immediately deduce that the generating function of $d$-oriented $k$-trees, counted with a weight $z$ per leaf, is equal to $v_k^{(d)}$ since it satisfies the same
Eq.~\recurv.
Otherwise stated,  $d$-irreducible $k$-slices with a given number $n$ of inner 
faces are in one-to-one correspondence with $d$-oriented $k$-trees with $n$ leaves.

\fig{Inductive definition of the bijection between $d$-irreducible
$k$-slices and $d$-oriented $k$-trees. Left: the $(d-2)$-slice reduced
to a single inner face is associated with the $(d-2)$-tree reduced to
a single leaf edge. Right: given a $k$-slice of size $>1$, we perform
its iterated decomposition: the corresponding $k$-tree is obtained by
keeping the root edge $AB$, which we decorate into an edge of type
$k/(d-1-k)$, and then inductively constructing the trees associated
with the subslices delimited by $P_1,P_2,\ldots$ (thus the edges
$BC_1,BC_2,\ldots$ will be kept in the tree at the next
stage).}{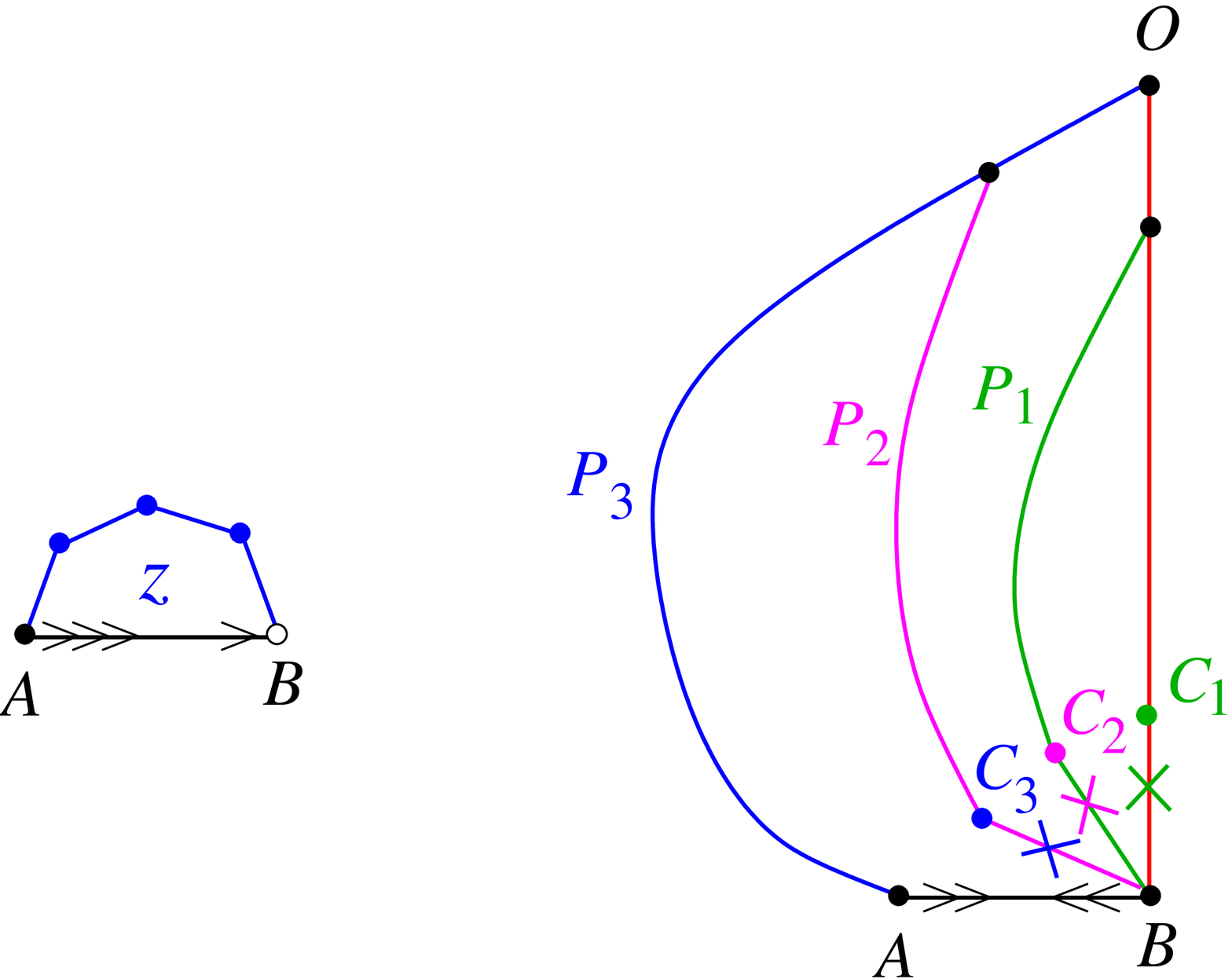}{8.cm}
\figlabel\slicetotree
\fig{Example of the full construction (c) of a $d$-oriented $k$-tree (b) from a $d$-irreducible $k$-slice (a) (here $d=5$ and $k=2$).}{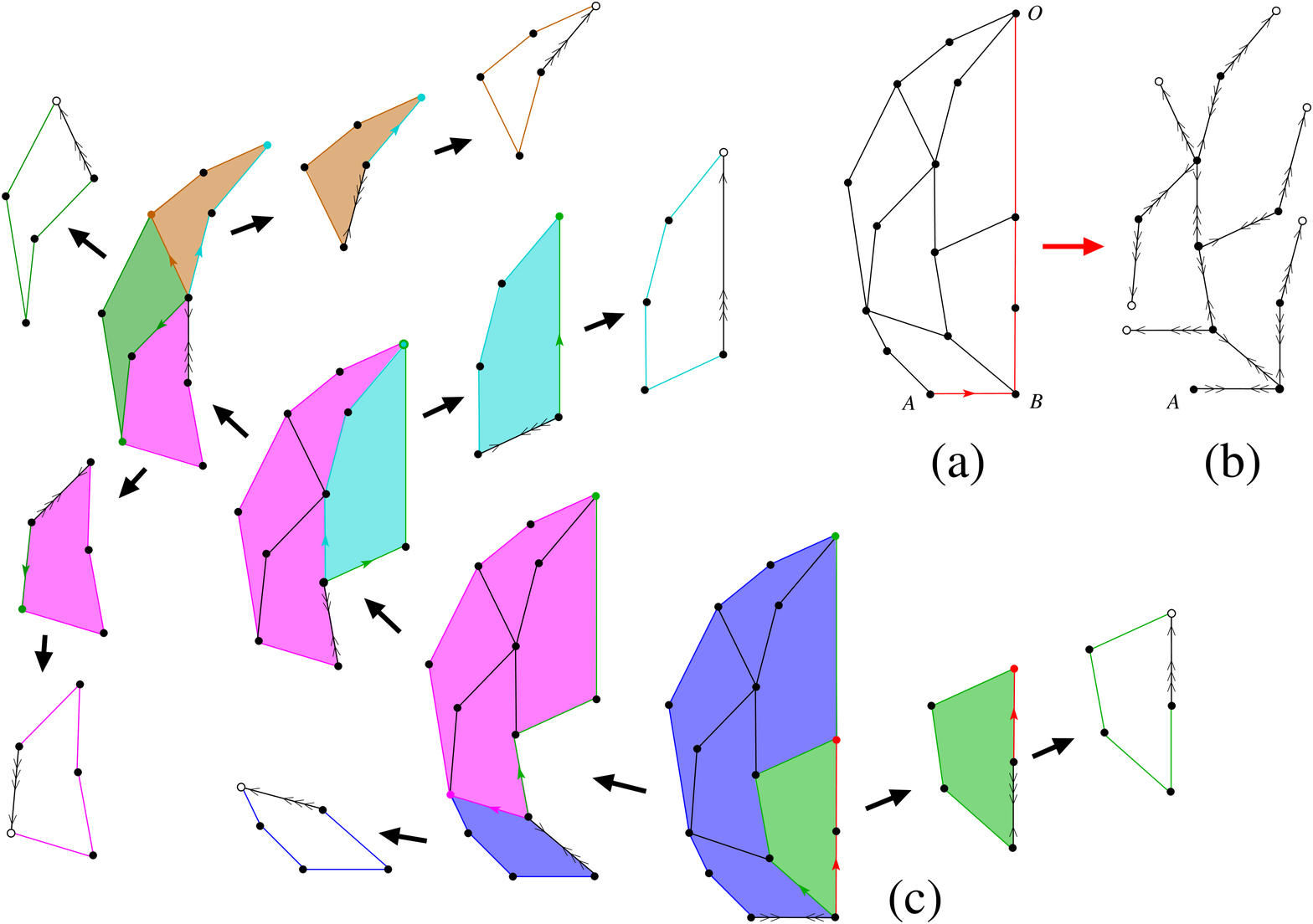}{14.cm}
\figlabel\opening
Obviously, an explicit bijection between these two sets can be defined
through their recursive decompositions. Recall from 
Sect.~\itdecompproc\ that the size of a $d$-irreducible $k$-slice is defined
as the number of its non-left edges (rather than that of its inner
faces). Similarly, we define the size of a $d$-oriented $k$-tree as
its number of edges (rather than that of leaves). The only $k$-slice of size $1$ is the
$(d-2)$-slice reduced to a single $d$-valent face, and we associate it
with the $d$-oriented $(d-2)$-tree reduced to a single leaf
edge. Suppose now that we have defined the bijection up to size
$N$. Given a $k$-slice of size $N+1$, we perform its iterated
decomposition to obtain a sequence of $q \geq 1$ of subslices whose
sizes add up to $N$. By induction hypothesis we may associate each
subslice with a $d$-oriented tree, and we merge the root vertices of
these trees together, and add an inner edge of type $k/(d-1-k)$ to
obtain the $k$-tree of size $N+1$ corresponding to the $k$-slice at
hand. Interestingly, the $k$-tree may be represented as a subgraph of
the $k$-slice by identifying at each step the added tree edge with 
the root edge of the slice at hand, see Fig.~\slicetotree. 
Fig.~\opening\ displays an
example of the construction of the $5$-oriented $2$-tree associated with 
a $5$-irreducible $2$-slice.

\fig{The rules of Fig.~\treerules\ in the particular case $d=3$. The out-degree of all inner vertices is
always $1$, so that inner vertices have necessarily degree $4$. We deduce that $3$-oriented $1$-trees 
reduce to ternary trees. }{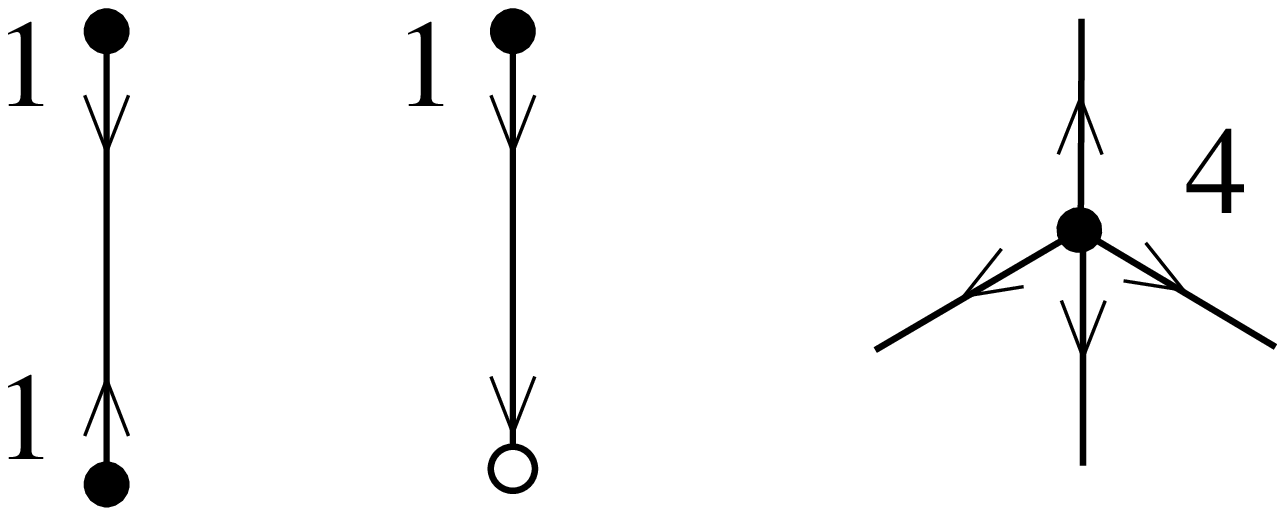}{7.cm}
\figlabel\ternarytrees

To conclude this section, let us discuss the particular case $d=3$. Taking the general rules of Fig.~\treerules\ at $d=3$, we see that
the edges of $3$-oriented trees are of the type displayed in Fig.~\ternarytrees, with an out-degree $1$ for each half-edge leaving an inner vertex. The 
constraint of out-degree  
$d+1=4$ at each inner vertex implies that all inner vertices have degree $4$, hence the tree is a ternary tree. In other words, irreducible triangular slices are in bijection with ternary trees, and we recover the bijection of \FusyT. The generating function $v^{(3)}_1$  may be identified with the generating function  of
planted ternary trees with a weight $z$ per leaf, and satisfies
\eqn\ternaryeq{v^{(3)}_1=z+(v^{(3)}_1)^3}
accordingly. This is nothing but Eq.~\recurv\ at $k=1$ while at $k=-1$ and $k=0$, this equation yields $v^{(3)}_{-1}=v^{(3)}_1$ and
$v^{(3)}_0=v^{(3)}_{-1}v^{(3)}_1$. Setting $s^{(3)}=v^{(3)}_{-1}=v^{(3)}_1$ and $r^{(3)}=1+v^{(3)}_0$, we recover Eq.~\algebRStri. 
Since the number of leaves
in a planted ternary tree is twice the number of inner vertices plus one, we see that, upon
setting $s^{(3)}=z T$ as in Sect.~2.2,  $T$ may be interpreted as the generating function for ternary trees with a weight $z^2$ per 
inner vertex, and satisfies  \eqforTbis\ accordingly.

\subsec{Alternative description of the mapping from trees to slices}
\subseclab\altertreeslice

\fig{The reconstruction of the $d$-irreducible $k$-slice of Fig.~\opening\ from its associated $d$-oriented $k$-tree 
(with $d=5$, $k=2$). Left: we first label the corners of the tree according to the rules explained in the text and call $-p-1$
the minimal label (here $p=3$). Middle: we then add $k+p+1$ new edges (blue edges), thus creating new corners
which we label $k,k-1,\cdots,-p-1$. Right: the slice is obtained by connecting each leaf to the first corner with same label 
encountered counterclockwise around the tree. This slice is of type $p/k+p+1$.}{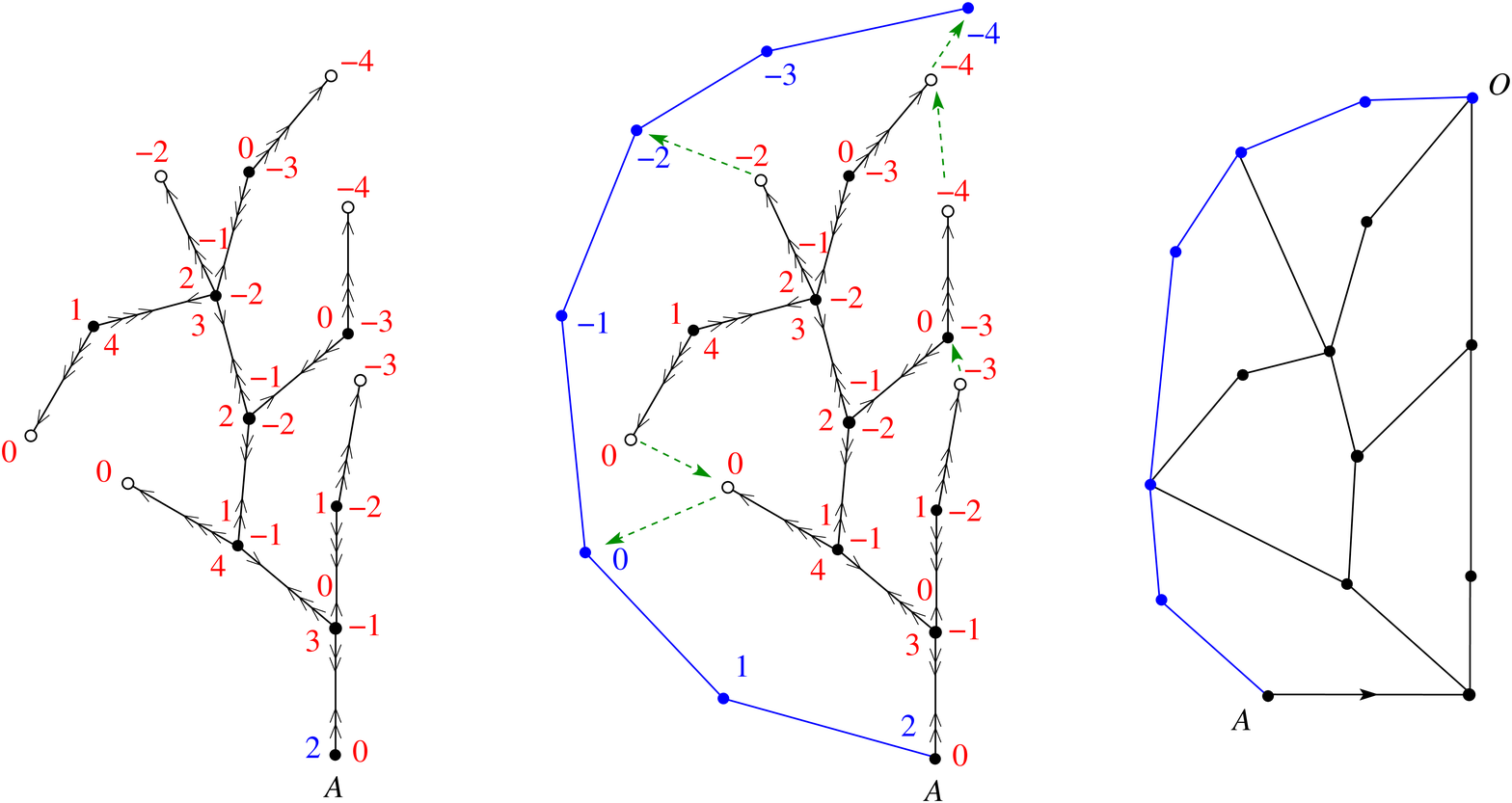}{11.cm}
\figlabel\closure

While the $k$-slice associated with a $k$-tree may be obtained
inductively by following the above construction backwards (we
decompose the $k$-tree into subtrees, construct inductively their
associated slices and glue them together), we find it worthwile to
describe an alternative (but equivalent) construction. It is a closing
procedure similar to that of \FPS, which consists in going
counterclockwise around the tree from its root vertex and, whenever a
leaf is followed by at least $d$ edge-sides before encountering a new
leaf, connecting it to the corner following immediately the $d$-th
encountered edge-side.
The procedure must in general be repeated several times, and a left boundary must be restored, before the slice is recovered. 

Let us now give a more precise definition of this procedure:
given a $d$-oriented $k$-tree, we consider its {\it contour walk} obtained by going counterclockwise around the tree
from its root vertex, thus visiting all corners of the tree successively upon following edge sides. Each corner receives a height equal to the height of the 
preceding corner minus $1$ unless the two corners are separated by an edge side that belongs to a leaf edge and corresponds to the second visit of this edge,
i.e.\ while going away from the leaf: in this case the height of the second corner is incremented by $d-1$ instead (see Fig.~\closure-left). Starting from height $0$,
we reach after a complete exploration of the tree the height $k$. Indeed, it is easily seen by induction that any subtree whose root vertex has out-degree 
$m$ separates two corners whose heights differ by $m$. Calling $-p-1$ (with $p \geq 0$) the minimum height obtained along the contour walk, necessarily
attained on a leaf, we complete
the tree by a sequence of $k+p+1$ new edge sides creating new corners with respective heights $k,k-1,\cdots, -p-1$ and accordingly extend the contour walk 
so as to end at the minimal height $-p-1$ (see Fig.~\closure-middle).  We define the {\it closure} of the tree as the slice obtained by simply connecting each leaf with 
height $m$ to its successor, which is the first corner
with the same height $m$ encountered along the contour walk (see Fig.~\closure-right). The apex of the slice is the (unique after connection) vertex with smallest height $-p-1$.
We leave it as an exercice to the reader to check that the closure of
the tree indeed coincides with the $k$-slice obtained by the inductive
bijection. 

Two remarks are in order: first, let us observe that our construction (in particular the
fact that we add new edges to obtain the left boundary) is somewhat
reminiscent of the construction of the ``discrete map with geodesic
boundaries'' associated with a labeled tree \LG\ (actually slices and
DMGBs are essentially the same objects). It is however unclear
that the two constructions could be unified:  the involved trees
($d$-oriented vs labeled) are quite different and, in the DMGB construction,
the slice is obtained by closing the tree on its right boundary
rather than the left one here . Second, our construction differs from
those of [\xref\FusyT,\xref\FPS] on the fact that constructing a slice requires only
a ``partial closure" (using the terminology of [\xref\FusyT,\xref\FPS]) whereas constructing a true dissection requires instead a ``complete closure" performed on some unrooted tree.
While it seems possible to define a complete closure of $d$-oriented trees,
we do not explore this direction here.

\subsec{Trees in the bipartite case} 
\subseclab\treebipcase
\fig{Simplification of the rules of Fig.~\treerules\ for $d$-oriented trees in the case $d=2b$.  Upon descending along the tree, the original 
edges (black) are converted into simpler edges (blue) as shown. The original rule of total out-degree $d+1$
for inner vertices translates after conversion into a constraint of total out-degree $b+1$ (bottom).}{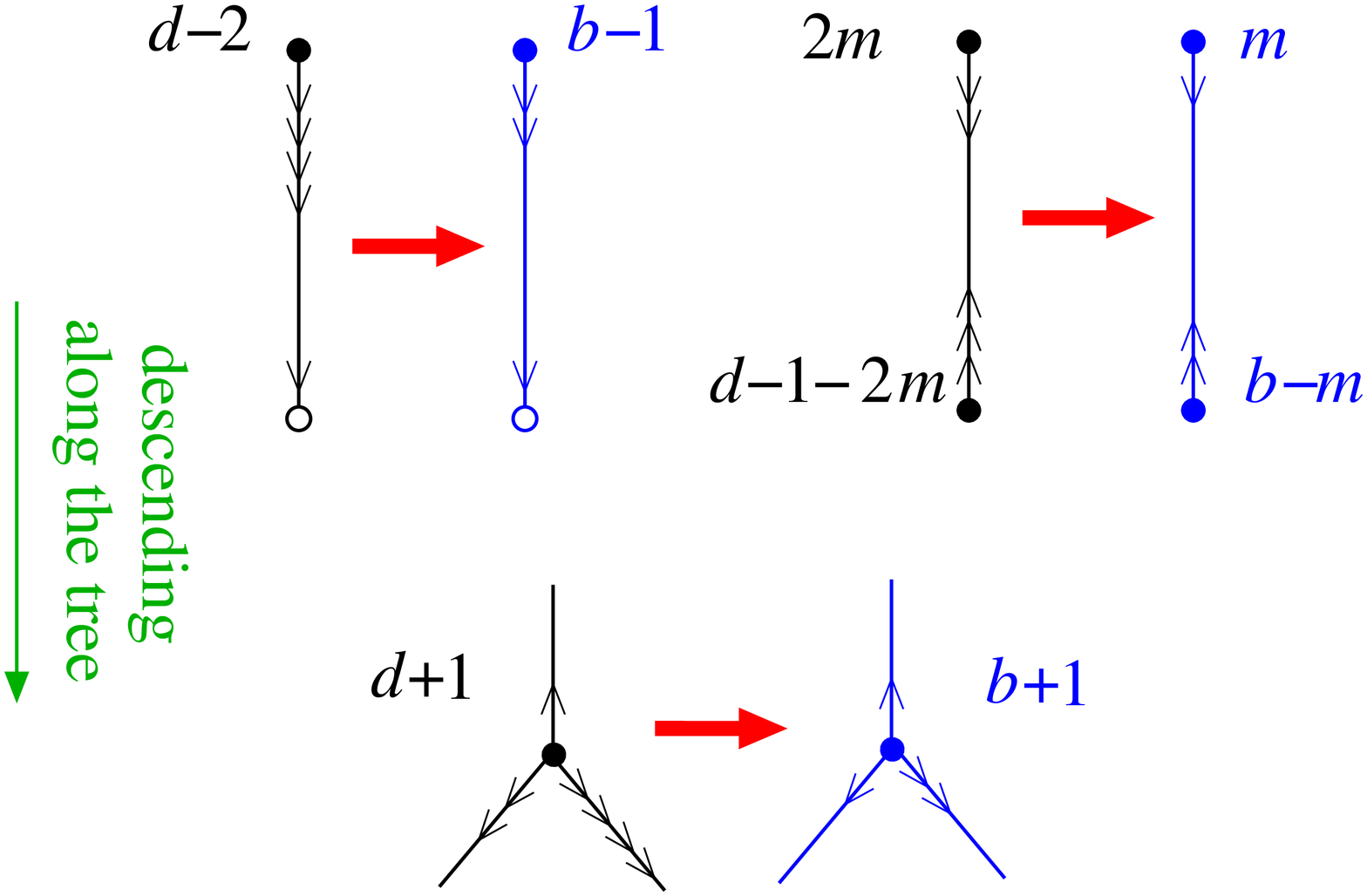}{10.cm}
\figlabel\treeconversion
We now consider the case of irreducible $d$-angulations when $d$ is even. As before, we then set $b=d/2$.
Let us first look at $d$-oriented $(2k+1)$-trees ($k\geq 0$), i.e\ trees with an odd root out-degree. Then the root edge is necessarily an inner edge
connecting the root vertex to an inner vertex $B$. The constraint $\sum_i m_i=2k+3$ for the sum of the out-degrees of 
the descending subtrees at $B$ implies that there is an odd number of such subtrees with odd root out-degree $m_i$.
In particular, there is at least one such subtree with odd root out-degree. By the same reasoning, this subtree has itself at least a descending subtree
with odd root out-degree and, by iteration, the $d$-oriented $(2k+1)$-tree cannot be finite. This is consistent with $v^{(d)}_{2k+1}=0$ as generating functions
count finite trees.
As for 
 $d$-oriented $2k$-trees, they are for the same reason build only out of leaf edges or 
 inner edges of type $2m/(d-1-2m)$, i.e.\  having an even number of arrows $2m$ followed by the complementary odd number of arrows $2b-1-2m$ when descending
 along the tree. We may then simplify the trees by converting these edges of type  $2m/(d-1-2m)$ into simpler edges of type
 $m/(b-m)$ (see Fig.~\treeconversion). Consistently we convert the leaf edges into edges with one half-edge carrying $b-1$ arrows pointing away from the 
 associated edge extremity and
the other half-edge carrying a single arrow pointing toward the associated extremity. Once this conversion is done, the out-degree of an inner vertex dangling from an edge of type $m/b-m$ becomes
$\sum_i m_i +b-m=b+1$ (since before conversion, we had $\sum_i 2m_i +2b-1-2m=2b+1$). The generating functions $u^{(d)}_k=v^{(d)}_{2k}$, viewed as generating
functions of $d$-oriented $2k$-trees, may thus alternatively be interpreted as counting trees, with a weight $z$ per leaf, such that:
\item{-} the edges of the tree carry arrows and are of two types:
\itemitem{-} inner edges of type $m/(b-m)$, with $1\leq m\leq b-1$, connecting only inner vertices or the
root vertex.
\itemitem{-} leaf edges with one half-edge, connected to an inner vertex or to the
root vertex, carrying $b-1$ arrows pointing away from this vertex and
the other half-edge carrying a single arrow pointing toward a leaf of the tree.
\item{-} the out-degree of any inner vertex is $b+1$.
\item{-} the out-degree of the root vertex is $k$.
\par
\noindent As such,  the $u^{(d)}_k$ satisfy 
\eqn\recuru{u_k^{(d)}=z \delta_{k,b-1}+\sum_{q\geq 1} \sum_{
1 \leq m_1,\ldots,m_q \leq b-1 \atop m_1+\cdots+m_q=k+1 } \prod_{i=1}^{q} u_{m_i}^{(d)} \qquad 1\leq k\leq b-1 ,}
with the convention $u^{(d)}_m=0$ for $m\geq b$. This system of equations is easily seen to follow from \recurU\ in the same way that \recurv\ follows from \recurV.

\fig{The rules of Fig.~\treeconversion\ (after conversion) in the particular case $d=4$ ($b=2$). The out-degree of all inner vertices is
always $1$, so that inner vertices have necessarily degree $3$. We deduce that $4$-oriented $2$-trees 
reduce to binary trees.}{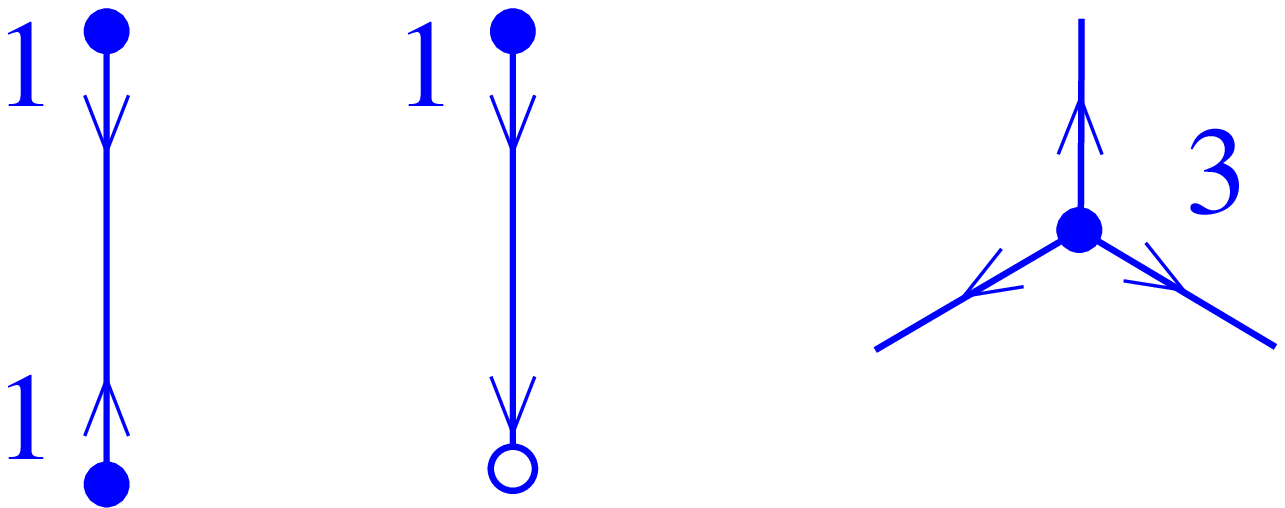}{7.cm}
\figlabel\binarytrees

To conclude this section, let us discuss the particular case $d=4$, i.e.\ $b=2$. After conversion, taking the general rules of 
Fig.~\treeconversion\ at $b=2$, we see that
the edges of the obtained trees are of the type displayed in Fig.~\binarytrees, with an out-degree $1$ for each half-edge leaving an inner vertex. The 
constraint of out-degree  
$b+1=3$ at each inner vertex implies that all inner vertices have degree $3$, hence the tree is a binary tree. In other words, irreducible quadrangular slices are in bijection with binary trees, and we recover the bijection of \FPS.
The generating function $u^{(4)}_1$ 
is nothing but the generating function of planted binary
trees with a weight $z$ per leaf, and satisfies
\eqn\ternaryeq{u^{(4)}_1=z+(u^{(4)}_1)^2}
accordingly. This is nothing but Eq.~\recuru\ at $k=1$ while, from \recurVbis\ at $k=0$, we deduce $u^{(4)}_0=v^{(4)}_0=v^{(4)}_2=u^{(4)}_1$. 
Using $r^{(4)}=1+u^{(4)}_0=1+u^{(4)}_1$, we then recover Eq.~\algebRquad. 
Setting $r^{(4)}=1+z T$ as in Sect.~2.1, i.e. $u^{(4)}_1=z T$, we see immediately that
$T$ may now be interpreted as the generating function of binary trees with a weight $z$ per 
inner vertex, and satisfies  \eqforT\ accordingly.

\newsec{Bijective proof of the pointing formulas}
\seclab\pointing

The purpose of this section is to give a combinatorial proof of the
general pointing formula \Fndpointgen\ which, by specialization,
implies the bipartite pointing formula \Fndbippoint, alternatively
obtained as a consequence of the computations of Sect.~\elimbip. 
Observe that the l.h.s.\ of \Fndpointgen, $\partial
F_n^{(d)}/\partial z$, is the generating function of $d$-irreducible
maps with outer degree $n$ and a marked inner face of degree $d$
(which does not receive a weight $z$). Such a map is said {\it
annular} and its marked face is called the {\it central face}.  Our
purpose is then to show that annular $d$-irreducible maps with outer
degree $n$ are in bijection with the objects naturally counted by
$P_d(n;R^{(d)},S^{(d)})$, namely three-step paths from $(0,0)$ to
$(n,d)$ endowed with the data of:
\item{-} for each down-step, a slice of type $p/p+1$ for some
arbitrary $p$ (i.e. either a $0$-slice or the map reduced to a single
root edge if $p=0$),
\item{-} for each level-step, a slice of type $p/p$ (i.e. a $-1$ slice).
\par\noindent
This bijection, which should preserve the total number of
(non central) inner faces of each degree, will of course be a variant
of the slice decomposition introduced in \PMCF\ and already encountered
in Sect.~\genslices.

\subsec{From paths to annular maps}
\subseclab\pathtoannul

\fig{Illustation of the gluing procedure. Left: the partial gluing of
a slice-decorated three-step path from $(0,0)$ to $(n,d)$. Right: the
folding of the left boundary onto the right boundary yielding an annular
map. (The heights of the path represent the relative distance of the
lower boundary vertices to the apex in the partial gluing, hence with the convention of
representing the slices pointing upwards the path appears vertically
reflected in the map.)}{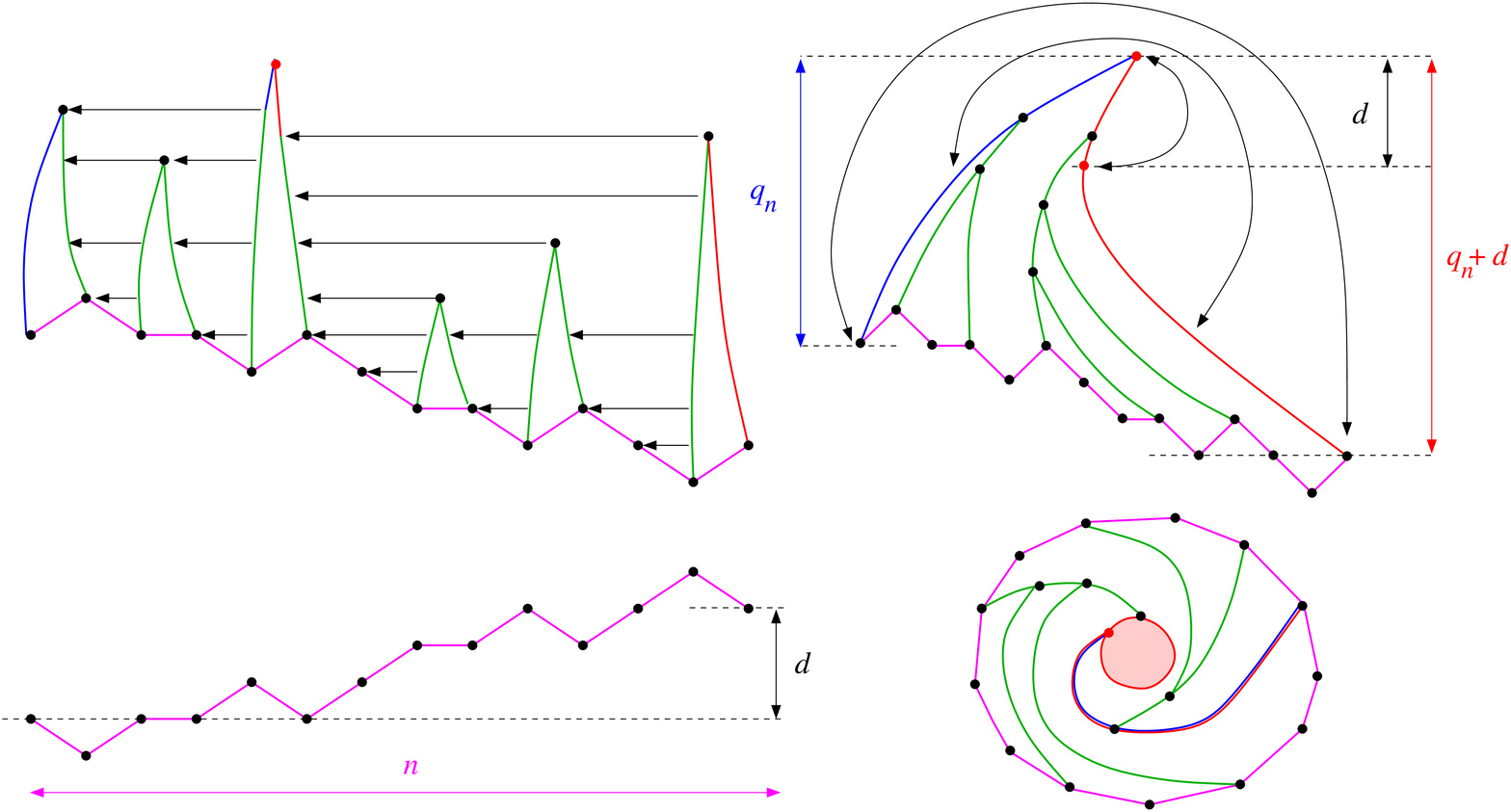}{14.cm}
\figlabel\pathtoannular

Here, it is simpler to first describe the mapping from slice-decorated
three-step paths to annular maps, see Fig.~\pathtoannular. It again
consists in gluing the slices together but we shall be precise about
the procedure. Let us denote by $h_0=0,h_1,\ldots,h_n=d$ the
successive heights of the three-step path and, for $1 \leq i \leq n$,
let ${\cal S}_i$ be the slice attached to the $i$-th step $h_{i-1} \to
h_i$ (for an up-step, ${\cal S}_i$ is by convention the map reduced to
a single edge).

We first define the {\it partial gluing} ${\tilde {\cal S}}_n$ of
${\cal S}_1,\ldots,{\cal S}_n$ inductively. We take ${\tilde {\cal
S}}_0$ to be the vertex-map (the map reduced to a single vertex with
no edge, which we view as a rooted map of outer degree $0$). Let us
now assume by induction that we have defined ${\tilde {\cal S}}_i$,
the partial gluing of ${\cal S}_1,\ldots,{\cal S}_i$, as a
$d$-irreducible map of outer degree $i+h_i+2 q_i$ for some $q_i \geq
\max(0,-h_i)$, satisfying the following constraint. When turning
around ${\tilde {\cal S}}_i$ in counterclockwise direction starting
from the root edge, we divide its outer boundary into three parts:
\item{-} a first part of length $i$, called the lower boundary,
\item{-} a second part of length $h_i+q_i$, called the right boundary,
which is the unique shortest path between its endpoints,
\item{-} a third part of length $q_i$, called the left boundary, which
is a (non necessarily unique) shortest path between its endpoints.
\par\noindent
Recall that the outer boundary of the slice ${\cal S}_{i+1}$ (assuming
that it corresponds to a level- or down-step) is also split into three
parts:
\item{-} a first path of length $1$, the root edge,
\item{-} a right boundary of length $p_i$, for some $p_i \geq 0$,
which is the unique shortest path between its endpoints,
\item{-} a left boundary of length $p_i+h_i-h_{i+1}$ (i.e. $p_i$ or
$p_i+1$ depending on whether the $i$-step is a level- or down-step),
which is a shortest path between its endpoints.
\par\noindent
For an up-step, we view ${\cal S}_{i+1}$ as having a right boundary of
length $1$ and a left boundary of length $0$. We then naturally define
${\tilde {\cal S}}_{i+1}$ by gluing the right boundary of ${\tilde
{\cal S}}_i$ to the left boundary of ${\cal S}_{i+1}$, identifying the
origin of the root edge of ${\cal S}_{i+1}$ with the first vertex of
the right boundary of ${\tilde {\cal S}}_i$ (in particular, ${\tilde
{\cal S}}_1={\cal S}_1$) .  Note that these boundaries do not
necessarily have the same length: for instance if $q_i>p_i-h_{i+1}$
then there remain some unmatched edges on the right boundary of ${\tilde
{\cal S}}_i$, which become part of the right boundary of ${\tilde
{\cal S}}_{i+1}$. It is easily seen that ${\tilde {\cal S}}_{i+1}$ has
outer degree $i+1+h_{i+1}+2q_{i+1}$, with
$q_{i+1}=\max(q_i,p_i-h_{i+1})$,
 and we take its root edge to be that of ${\tilde {\cal
S}}_i$. Furthermore, using the data that ${\tilde {\cal S}}_{i}$ and
${\cal S}_{i+1}$ are $d$-irreducible and that their left/right
boundaries are shortest paths, it is not difficult to check that the
same properties hold for ${\tilde {\cal S}}_{i+1}$. Thus, we fulfill
the induction hypothesis.

Having defined the partial gluing ${\tilde {\cal S}}_n$, we finish by
gluing its right boundary (of length $d+q_n$) to its left boundary (of
length $q_n$), identifying the origin of the root edge to the
first vertex of the right boundary. The $d$ unmatched edges of the
right boundary yield a central inner face of degree $d$, and the lower
boundary yields an outer face of degree $n$, thus we obtain an annular
map, which is easily shown to be $d$-irreducible. In particular, to
see that the only cycle of length $\leq d$ winding around the central
face is its boundary, we use the fact that the right boundary is the
unique shortest path between its endpoints, so that each vertex of the
left boundary of ${\tilde {\cal S}}_n$ except the first one is
identified with a vertex of the right boundary at distance $>d$ in
${\tilde {\cal S}}_n$.

Finally, let us note the cyclic invariance of the construction. More
precisely, we let $\epsilon_i=h_i-h_{i-1} \in \{-1,0,1\}$, $1 \leq i
\leq n$, denote the $i$-th {\it increment} of our three-step
path. Then, for any $m$ between $1$ and $n$, we consider the
circularly shifted sequence
$\epsilon_m,\epsilon_{m+1},\ldots,\epsilon_n,\epsilon_1,\dots,\epsilon_{m-1}$,
which is the sequence of increments of another three-step path from
$(0,0)$ to $(n,d)$.  Attaching to this path the slice sequence ${\cal
S}_m,{\cal S}_{m+1},\ldots,{\cal S}_n,{\cal S}_1,\ldots,{\cal
S}_{m-1}$, we obtain another slice-decorated three-step path to which
we may apply the gluing procedure. It is not difficult to see that the
resulting annular map is the same, except for the position of the root
edge which is moved by $m-1$ steps along the outer face.

\subsec{Slice decomposition of annular maps}
\subseclab\annultopath

We now explicit the inverse mapping from annular maps to
slice-decorated three-step paths. Because of the cyclic invariance
noted above, it is desirable to have a construction in which the root
edge does not play a specific role. Such a construction turns out to
be naturally described on the {\it lift} of the annular map which we
define as follows. Let ${\cal M}$ be a $d$-irreducible annular map of
outer degree $n$, which we think of as being drawn on the complex plane,
the origin being in the interior of the central face. Then, we define
its lift ${\tilde {\cal M}}$ as the preimage of ${\cal M}$ by the
mapping $z \mapsto \exp(2i\pi z)$, i.e.\ the map whose vertices and
edges are the preimages of those of ${\cal M}$ (it is not difficult to
convince oneself that the notion of lift is well behaved with respect
to continuous deformation). The map ${\tilde {\cal M}}$ is infinite
but locally finite, i.e. each vertex has finite degree (equal to the
degree of its image-vertex in ${\cal M}$). Furthermore, the preimage
of the central (resp.\ outer) face of ${\cal M}$ forms a single face
of infinite degree, the {\it upper} (resp.\ {\it lower}) face of
${\tilde {\cal M}}$, while each inner non central face of ${\cal M}$
yields infinitely many faces of finite degree.
The translation $z \mapsto z+1$ induces a natural automorphism $T$ of
${\tilde {\cal M}}$, and we endow ${\tilde {\cal M}}$ with the graph
distance ${\tilde D}(\cdot,\cdot)$.

\fig{Illustration of the slice decomposition of an annular map ${\cal M}$ (top). We pass to its lift $\tilde{\cal M}$ (bottom) and for 
each vertex $\ell_i$ of the lower boundary, we draw the leftmost shortest path from $\ell_i$ to $u_{K(\ell_i)}$,
where $K(\ell_i)$ is defined by the coalescence lemma (intuitively speaking, this path is the ``leftmost shortest path
from $\ell_i$ to $-\infty$''). Clearly the figure is invariant by the translation $T$ and we obtain the wanted slice decomposition 
of ${\cal M}$.}{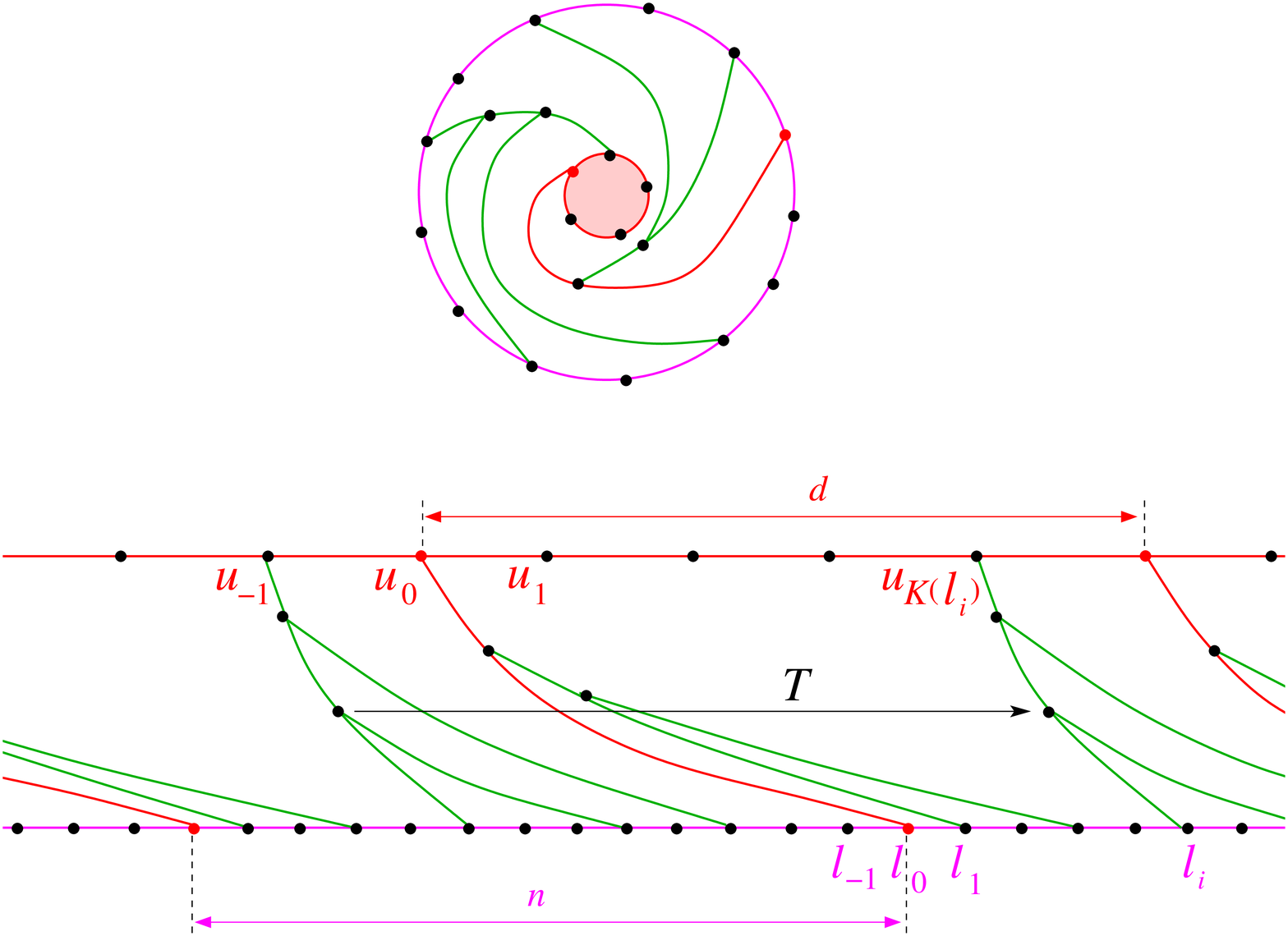}{12.cm}
\figlabel\lift

Since we are in the plane, the notion of leftmost shortest path from a
given vertex of ${\tilde {\cal M}}$ to another is defined in an
obvious manner. Let us denote by
$\ldots,u_{-2},u_{-1},u_0,u_1,u_2,\ldots$ the vertices incident to the
upper face, read in successive order by following the boundary from
$-\infty$ to $+\infty$ (the choice of $u_0$ will turn out to be
irrelevant). Note that $T(u_k) = u_{k+d}$ for all $k$. We then have
the following: \smallskip

\item{} {\bf Coalescence lemma:} for each vertex $v$ of ${\tilde {\cal
M}}$, there exists a unique integer $K(v)$ such that, for all $k \leq
K(v)$, any shortest path from $v$ to $u_k$ passes through
$u_{K(v)}$ (and hence follows the upper boundary from $u_{K(v)}$ to
$u_k$), but some of the shortest paths do not pass through $u_{K(v)+1}$. 
In particular we have ${\tilde
D}(v,u_k)-(K(v)-k)={\tilde D}(v,u_{K(v)}) \leq {\tilde D}(v,u_{K(v)+1})$.
\smallskip

\noindent In other words, for $k \to -\infty$, all shortest paths to
$u_k$ eventually coalesce with the upper boundary. The proof of the
lemma, which relies crucially on the $d$-irreducibility of ${\cal M}$,
is postponed to the end of this section. Note that $K(T(v))=K(v)+d$ by
translation invariance.

We now denote by $\ldots,\ell_{-1},\ell_0,\ell_1,\ell_2,\ldots$ the
successive vertices incident to the lower face, again read from
$-\infty$ to $+\infty$ so that $T(\ell_i)=\ell_{i+n}$ for all $i$ (see Fig.~\lift).
Here we pick $\ell_0$ as an arbitrary preimage of the origin of the
root edge of ${\cal M}$, and, at the price of relabeling the vertices of the
upper boundary $u_k\to u_{k-K(\ell_0)}$, we may assume that $K(\ell_0)=0$. 
We then let
\eqn\epsi{\epsilon_i = {\tilde D}(\ell_{i},u_{K(\ell_{i-1})}) - {\tilde
D}(\ell_{i-1},u_{K(\ell_{i-1})})}
which clearly belongs to $\{-1,0,1\}$ by the triangular inequality. By
translation invariance we have
$\epsilon_i=\epsilon_{i+n}$. Furthermore, it is easily seen that by
planarity $K(\ell_{i-1}) \leq K(\ell_i)$ so that, using the
coalescence lemma,
\eqn\epsisum{\eqalign{\sum_{i=1}^n \epsilon_i &= \sum_{i=1}^n \left(
{\tilde
D}(\ell_i,u_{K(\ell_i)}) + K(\ell_i) - K(\ell_{i-1})-{\tilde D}(\ell_{i-1},u_{K(\ell_{i-1})}) \right) \cr &=
{\tilde D}(\ell_n,u_{K(\ell_n)}) - {\tilde D}(\ell_0,u_{K(\ell_0)}) +
K(\ell_n) - K(\ell_0) = d}}
since $T(\ell_0)=\ell_n$ so that $K(\ell_n) = K(\ell_0)+d$ and
${\tilde D}(\ell_0,u_{K(\ell_0)}) = {\tilde
D}(\ell_n,u_{K(\ell_n)})$. Hence the sequence
$\epsilon_1,\ldots,\epsilon_n$ forms the increments of a three-step
path from $(0,0)$ to $(n,d)$. Let us now consider the bounded region
delimited by the edge $(\ell_{i-1},\ell_i)$ and the two leftmost
shortest paths from $\ell_{i-1}$ and $\ell_i$ to $u_{K(\ell_{i-1})}$. When
$K(\ell_{i-1})=K(\ell_i)$ these two paths may merge before their
endpoint and we then remove their common part, letting $p_i$ be the
length of the proper part of the path with starting point $\ell_i$.
Clearly, the submap enclosed within this region is nothing but a slice
of type $p_i/p_i-\epsilon_i$, which we denote by ${\cal S}_i$ (when
$\epsilon_i=+1$ we obtain a slice reduced to a single edge as
wanted). Note that ${\cal S}_{i+n}$ is simply a translate of ${\cal
S}_i$, so no information is lost by restricting to the interval
$\{1,\ldots,n\}$.

To summarize, starting from a $d$-irreducible annular map of outer
degree $n$, we have constructed a slice-decorated three-step path from
$(0,0)$ to $(n,d)$. It is clear that applying then the gluing
procedure of the previous section restores the original annular
map. To establish that we have a bijection, it remains to check that a
slice-decorated three-step path is indeed recovered as the slice
decomposition of its gluing.

Consider such a path and its partial gluing ${\tilde {\cal S}}_n$ as
defined in the previous section. The lift ${\tilde {\cal M}}$ of the
corresponding annular map ${\cal M}$ is obtained by gluing infinitely
many copies $\ldots,{\tilde {\cal S}}_n^{(-1)},{\tilde {\cal
S}}_n^{(0)},{\tilde {\cal S}}_n^{(1)},\ldots$ of ${\tilde {\cal S}}_n$
to each other along their left/right boundaries, see
Fig.~\lift. Consistently with the previous notations, let
$\ell_{ni},\ell_{ni+1},\ldots,\ell_{n(i+1)}$ be the vertices of the
lower boundary of ${\tilde {\cal S}}_n^{(i)}$, and
$u_{di},u_{di+1},\ldots,u_{d(i+1)}$ be the vertices of the upper
boundary, i.e. the last $d+1$ vertices of the right boundary in
counterclockwise direction (note indeed that $\ell_{ni}$ and $u_{di}$
are identified with vertices of ${\tilde {\cal S}}_n^{(i-1)}$ in the
gluing procedure). Using the data that the left and right boundary of
${\tilde {\cal S}}_n$ are leftmost shortest paths, it is not difficult
to check that the leftmost shortest path in ${\tilde {\cal M}}$ from
$\ell_{ni}$ to $u_k$ follows the left boundary of ${\tilde {\cal
S}}_n^{(i)}$ for all $k \leq di$, in particular it passes through
$u_{di}$ but not $u_{di+1}$, so that $K(\ell_{ni})=di$ where
$K(\cdot)$ is defined as in the coalescence lemma. This implies that
the slice decomposition of ${\tilde {\cal M}}$ cuts precisely along
the boundaries of the ${\tilde {\cal S}}_n^{(i)}$ for all
$i$. Furthermore, within each ${\tilde {\cal S}}_n^{(i)}$, the
leftmost shortest paths from $\ell_{ni+k}$, $0 \leq k \leq n$, to
$u_{di}$ precisely delimit the same slices as those attached to the
original path, so the slice decomposition is indeed the inverse
mapping of the gluing procedure.

It remains to prove the coalescence lemma. We start with the
intermediate:
\smallskip

\item{} {\bf Wrapping lemma:} for $v$ a vertex of ${\tilde {\cal M}}$
and $m$ an integer, we have ${\tilde D}(v,T^m(v)) \geq |m|\, d$. For
$m \neq 0$, equality holds if and only if $v$ is incident to the upper
face, and in this case the unique shortest path from $v$ to $T^m(v)$
follows the upper boundary. \smallskip

\noindent {\it Proof:} we may assume $m>0$ without loss of generality,
upon exchanging the role of $v$ and $T^m(v)$. Let us consider a path
from $v$ to $T^m(v)$: in ${\cal M}$ its image is a closed path of the
same length whose winding number around the central face is $m$. The
case $m=1$ is simpler: by possibly removing some edges we obtain a
cycle of ${\cal M}$ which, by $d$-irreducibility, has length $\geq d$
with equality iff it coincides with the boundary of the central
face. For $m>1$, note that the closed path cannot be simple (since the
possible winding numbers for a cycle are $0$ or $\pm 1$), so it has a
multiple vertex. Splitting at this vertex, we obtain two closed paths
of smaller length, whose winding numbers (around the central face) add
up to $m$. If one of these subpaths has winding number $>1$, we
further split it, and so on until we are left with a collection of
closed paths of winding numbers $\leq 1$. In particular we have at
least $m$ subpaths of winding number $1$, each of them of length $\geq
d$. Thus the original path has length $\geq m d$, and equality holds
iff all subpaths coincide with the boundary of the central face,
Q.E.D. \smallskip

A corollary of the wrapping lemma is that, for all integers $i,j$, we
have ${\tilde D}(u_i,u_j)=|j-i|$ and the upper boundary is the unique
shortest path from $u_i$ to $u_j$: consider a path of length $L$
between them and, assuming without loss of generality that $i<j$, pick
$m>0$ such that $i+m d \geq j$.  By appending the upper boundary from
$u_j$ to $u_{i+m d}=T^m(u_i)$ we obtain a path of length $L+i+md-j$
from $u_i$ to $T^m(u_i)$, so that $L+i+md-j \geq m d$, i.e.\ $L \geq j-i$ and equality holds iff
the path follows the upper boundary.  \smallskip

\noindent {\it Proof of the coalescence lemma:} let $v$ be again an
arbitrary vertex of ${\tilde {\cal M}}$. By the triangular inequality
${\tilde D}(v,u_k) \geq |k| - {\tilde D}(v,u_0) \to +\infty$ for $k
\to \pm \infty$, so in particular there exists an integer $K'(v)$ such
that, for all $k \leq K'(v)$, ${\tilde D}(v,u_k)$ is larger than the
number $M$ of vertices of ${\cal M}$. By the pigeonhole principle, any
path in ${\cal M}$ of length $\geq M$ passes necessarily twice through
the same vertex.  We apply this principle to the image of a shortest
path $\gamma$ from $v$ to $u_k$, and conclude that $\gamma$
necessarily passes successively through two distinct vertices $v'$ and
$v''$ of ${\tilde {\cal M}}$ having the same image in ${\cal M}$,
i.e. $v''=T^m(v')$ for some $m \neq 0$. Actually, we may assume that
${\tilde D}(v,v') < M$ by applying the principle to the prefix of
length $M$ of $\gamma$. Cutting $\gamma$ at $v'$ and $v''$, we obtain
three subpaths $\gamma'$, $\gamma''$ and $\gamma'''$. Then, by
concatenating $\gamma'$ with $T^{-m}(\gamma''')$ and then with the
upper boundary from $u_{k-md}$ to $u_{k}$, we obtain another path from
$v$ to $u_k$ which has length ${\tilde D}(v,u_k)- {\tilde
D}(v',v'')+|m|d$. Its length should be not lesser than ${\tilde
D}(v,u_k)$ and we deduce from the wrapping lemma that $v'$ is incident
to the upper face, i.e.\ $v'=u_{k'}$ for some $k'$. From ${\tilde
D}(v,v') < M$ we deduce that $k'>K'(v)$, and since the concatenation
of $\gamma''$ and $\gamma'''$ is a shortest path between $u_{k'}$ and
$u_k$ we deduce from the corollary of the wrapping lemma that it
passes through $u_{K'(v)}$, hence $u_{K'(v)}$ belongs to any shortest
path from $v$ to $u_k$ for all $k \leq K'(v)$. This shows that the set
of the integers $K$ such that, for all $k \leq K$, $u_K$ belongs to
the leftmost shortest path from $v$ to $u_k$ is nonempty. It is
bounded from above since ${\tilde D}(v,u_k) \to +\infty$ for $k \to
+\infty$ thus contains a maximal element $K(v)$ which clearly
satisfies the conditions of the coalescence lemma, Q.E.D.

\newsec{Discrete integrable equations}

It has been noted in several occasions that maps are related to
integrable systems, for instance, map generating functions are
tau-functions of the KP hierarchy \GJKP. In the planar setting, a
slightly different connection has been uncovered: generating functions
of maps with marked points at a prescribed distance have been shown to
satisfy a hierarchy of ``discrete integrable equations'' [\xref\GEOD-\xref\PDFRaman], whose combinatorial meaning is now quite understood \PMCF.
Remarkably, the integrability phenomenon subsists in the context of
irreducible maps.

Integrable equations are obtained by adding a new parameter in the
game: at the combinatorial level it consists in controlling the
maximal length of the slices. Both the substitution and the slice
decomposition approaches still work in this case. Here, we choose to
emphasize the cases of irreducible quadrangular and
triangular dissections, as they are related to naturally embedded trees. For
simplicity, we follow the substitution approach, with an analysis
parallel to that of Sect.~\firstsimp. The general case, and the
slice decomposition, are then discussed in the last subsection.

\subsec{Integrable equations from irreducible quadrangular dissections}

In the case of quadrangular dissections, the starting point is some refined version of Eq.~\Fnquad\ of Sect.~2.1, which 
states that [\xref\DFG,\xref\PMCF]
\eqn\Fnquadi{\eqalign{F_2&=R_i - \left(x_4 R_{i-1}R_iR_{i+1}\right),\cr
F_4&= R_i(R_i+R_{i+1})- (R_i+R_{i+1}+R_{i+2}) \left(x_4 R_{i-1}R_iR_{i+1}\right),\qquad \qquad i\geq 1,\cr
}}  
where $R_i$ is the solution of the equation
\eqn\eqRi{ R_i=1+x_2 R_i+x_4 R_i (R_{i-1}+R_i+R_{i+1}),\qquad i\geq 1,}
with $R_0=0$ (here, besides squares weighted by $x_4$, we also possibly allow bivalent faces with weight $x_2$).
More precisely, whenever \eqRi\ is satisfied, the expressions in the r.h.s.\ of \Fnquadi\ are {\it conserved quantities},
i.e.\ their value does not depend on $i$ and may moreover be identified with $F_2$ and $F_4$ respectively. 
This is a sign of integrability for Eq.~\eqRi, which indeed admits the explicit solution \GEOD
\eqn\solforRi{R_i=R {(1-y^i)(1-y^{i+3})\over (1-y)^{i+1}(1-y)^{i+2}}, \qquad i\geq 0, \qquad \hbox{where}\ y+{1\over y}+1={1\over x_4 R^2},}
with $R$ solution of the homogeneous (i.e.\ without indices) version of \eqRi, $R=1+x_2 R+3 x_4 R^2$.
Performing the same specialization $(x_2,x_4)\to (X_2(z),X_4(z))$ as in Sect.~2.1, we deduce that
the generating functions $f_2^{(4)}(z)$ and $f_4^{(4)}(z)$ of irreducible quadrangulations with 
a boundary of length $2$ satisfy a generalization of \Fnquadd, namely
\eqn\Fnquaddi{\eqalign{f_2^{(4)}&=1=r^{(4)}_i - \left(X_4(z) r^{(4)}_{i-1}r^{(4)}_{i}r^{(4)}_{i+1}\right),\cr
f_4^{(4)}&=2+z=r^{(4)}_i(r^{(4)}_i+r^{(4)}_{i+1})-\, (r^{(4)}_i+r^{(4)}_{i+1}+r^{(4)}_{i+2} )\left(X_4(z) r^{(4)}_{i-1}r^{(4)}_{i}r^{(4)}_{i+1}\right),\qquad  i\geq 1,\cr
}}
where $r^{(4)}_i=R_i(X_2(z),X_4(z))$. Eliminating $X_4(z)$, we see that $r_i^{(4)}$ is now fully 
determined by the equation
\eqn\inteqrifour{z+r^{(4)}_i r^{(4)}_{i+2}- (r^{(4)}_i+r^{(4)}_{i+1}+r^{(4)}_{i+2})+2=0, \quad i\geq 1,}
with initial conditions $r^{(4)}_1=1$ and $r^{(4)}_2=1+z$ (from \Fnquaddi\ at $i=1$ with $r^{(4)}_0=0$). In practice,
we may equivalently extend Eq.~\inteqrifour\ to include the case $i=0$ and use as initial conditions $r^{(4)}_0=0$ and 
$r^{(4)}_1=1$.
By substituting $(x_2,x_4)\to (X_2(z),X_4(z))$ in \solforRi, we readily
obtain the explicit expression
\eqn\solforrifour{r^{(4)}_i=r^{(4)} {(1-y^i)(1-y^{i+3})\over (1-y^{i+1})(1-y^{i+2})}, \qquad i\geq 0, \qquad \hbox{where}\ y+{1\over y}={1\over r^{(4)}-1},}
with $r^{(4)}$ solution of \algebRquad. Here, we simplified the equation for $y$, 
upon using the relation $X_4(z) (r^{(4)})^2=(r^{(4)}-1)/r^{(4)}$ read off \Fnquadd. 
We observe that the same expression appears in [\xref\AlFuPo,
Proposition 4.5], where it is interpreted as a generating function of
symmetric irreducible quadrangular dissections (such objects are
indeed obtained by gluing several copies of a same slice counted by
$r^{(4)}_i$).

As in Sect.~2.1, Eq.~\inteqrifour\  is made more transparent upon setting 
\eqn\rtoTi{r^{(4)}_i(z)=1+z \,T_{i-1}(z),\qquad i\geq 1,}
as it then reads
\eqn\eqforTi{T_i=1+z\, T_{i-1}T_{i+1}, \qquad i\geq 1,}
with initial condition $T_0=0$. Here we recognize the equation determining the generating function of ``naturally embedded binary trees" (NEBT) introduced
in \BM, whose integrability remained so far quite mysterious. 
Now, from \solforrifour\ and the first line of \Fnquaddi, we simply deduce
\eqn\solforTi{\eqalign{T_i={r^{(4)}_{i+1}-1\over z}={X_4(z)\over z} r^{(4)}_{i}r^{(4)}_{i+1}r^{(4)}_{i+2}& =T{(1-y^i)(1-y^{i+5})\over (1-y^{i+2})(1-y^{i+3})}, \qquad
i\geq 0,\cr & \qquad \hbox{where} \ y+{1\over y}={1\over z T}\cr}}
and $T=(r^{(4)}-1)/ z$ is solution of \eqforT. We recover the
particular form of the solution found in [\xref\BM, Prop.~25]. Let us summarize the
combinatorial steps of our derivation. First, we note that
Eq.~\eqforTi, which is easily interpreted in the language of
NEBT, also follows from slice decomposition, see Sect.~8.3 below.
Second, generating functions of irreducible and arbitrary slices are
simply related by a change of variables. Thus, the explicit form
\solforRi, combinatorially explained in [\xref\PMCF, Sect.~6.2], directly
translates into \solforrifour. Finally, the fact that
$T_i=(r^{(4)}_{i+1}-1)/z$ also has a nice factorized form is a
consequence of the first line of \Fnquadi, combinatorially explained
in [\xref\PMCF, Sect.~3.3].

\subsec{Integrable equations from irreducible triangular dissections}

If we now play the same game for triangular dissections, we may use the following expressions \PMCF,
valid for $i\geq 1$:
\eqn\Fntrii{\eqalign{F_1&=S_{i-1} - \left(x_3 R_iR_{i-1}\right),\cr
F_2&= (S_{i-1})^2+R_i-  (S_{i-1}+S_i) \left(x_3 R_i R_{i-1}\right),\cr
F_3&= (S_{i-1})^3+ R_i (2S_{i-1}+S_i)-((S_{i-1})^2+S_{i-1}S_i+(S_i)^2+ R_i+R_{i+1}) \left(x_3 R_iR_{i-1}\right),\cr
}}
where $R_i$ and $S_{i-1}$ ($i\geq 1)$ are determined via
\eqn\eqforRiSi{R_i=1+x_2 R_i + x_3 R_i (S_{i-1}+S_i),\qquad  S_{i-1}=x_1+x_2 S_{i-1}+ x_3((S_{i-1})^2+ R_{i-1}+R_i)}
with $R_0=0$.
These equations are again integrable, with explicit solution
\eqn\solforRiSi{\eqalign{& S_{i-1}=S-x_3 R^2 y^{i-1} {(1-y)(1-y^2)\over (1-y^i)(1-y^{i+1})}, \qquad R_i=R {(1-y^i)(1-y^{i+2})\over (1-y^{i+1})^2}\cr
&  \qquad\qquad i\geq 1 \qquad \hbox{where}\ \ \ y+{1\over y}+2={1\over x_3^2 R^3}\cr}}
with $R$ and $S$ solutions of the homogeneous version of \eqforRiSi:  $R=1+x_2 R + 2 x_3 R S$, $S=x_1+x_2 S+ x_3(S^2+2 R)$.

Specializing these expressions at the particular renormalized values $x_1=X_1(z)$, $x_2=X_2(z)$, $x_3=X_3(z)$ of
Sect.~2.2, we deduce expressions for the generating functions of irreducible triangulations with a boundary of length $1$, $2$ and $3$:
\eqn\Fntridi{\eqalign{f^{(3)}_1&=0=s^{(3)}_{i-1} - \left(X_3(z) r^{(3)}_{i-1} r^{(3)}_i \right),\cr
f^{(3)}_2&= 1=(s^{(3)}_{i-1})^2+r^{(3)}_i-(s^{(3)}_{i-1}+s^{(3)}_i) \left(X_3(z) r^{(3)}_{i-1} r^{(3)}_i \right),\cr
f^{(3)}_3&=z=(s^{(3)}_{i-1})^3+ r^{(3)}_i (2 s^{(3)}_{i-1}+s^{(3)}_i)\cr & \qquad -((s^{(3)}_{i-1})^2+s^{(3)}_{i-1}s^{(3)}_i+(s^{(3)}_i)^2+ r^{(3)}_i+r^{(3)}_{i+1})
 \left(X_3(z) r^{(3)}_{i-1} r^{(3)}_i \right), \cr
}}  
in terms of $r^{(3)}_i(z)=R_i(X_1(z),X_2(z),X_3(z))$ and $s^{(3)}_i(z)=S_i(X_1(z),X_2(z),X_3(z))$.
Eliminating $X_3$, $r^{(3)}_i$ and $s^{(3)}_i$ are fully determined in terms of $z$ via:
\eqn\algebRStrii{r^{(3)}_i=1+s^{(3)}_{i-1}s^{(3)}_i, \qquad z+s^{(3)}_{i-1}s^{(3)}_i s^{(3)}_{i+1}-s^{(3)}_i=0\ ,\qquad i\geq 1}
with initial condition $s^{(3)}_0=0$.
Using $s^{(3)}=X_3(z) (r^{(3)})^2$ from the first line of \Fntrid, formulas \solforRiSi\ specialize to
\eqn\solforrisitri{\eqalign{s^{(3)}_{i-1}& =s^{(3)}\left(1- y^{i-1} {(1-y)(1-y^2)\over (1-y^i)(1-y^{i+1})}\right)=s^{(3)} {(1-y^{i-1})(1-y^{i+2})\over (1-y^i)(1-y^{i+1})}\cr 
 r^{(3)}_i & =r^{(3)} {(1-y^i)(1-y^{i+2})\over (1-y^{i+1})^2}, \qquad i\geq 1\qquad 
 \hbox{where}\ \ \ y+{1\over y}+1={1\over (s^{(3)})^2}\cr}}
with $r^{(3)}$ and $s^{(3)}$ solutions of the system \algebRStri. Here we used $(X_3(z))^2 (r^{(3)})^3=(s^{(3)})^2/r^{(3)}=(s^{(3)})^2/((s^{(3)})^2-1)$.
See also [\xref\AlFuPo, Proposition 5.5] for equivalent expressions
interpreted as generating functions of symmetric irreducible
triangular dissections.

Again the equation for $s^{(3)}_i$ in \algebRStrii\ is more transparent upon setting 
\eqn\stoTi{s^{(3)}_i(z)=z \,T_i(z) ,\qquad i\geq 0}
as it reads
\eqn\eqforTibis{T_i=1+z^2\, T_{i-1}T_iT_{i+1}, \qquad i\geq 1}
with $T_0=0$, allowing to identify $T_i$ with the generating function of some ``naturally embedded" ternary trees \Kuba\ (with a weight $z^2$ per inner vertex). 
Eq.~\eqforTibis\ was already encountered  in the context of quadrangulations without multiple edges \DSQ, and this explained its integrability.
It is remarkable to find it here in the different context of irreducible triangulations, hence providing a new explanation for integrability.
Let us mention that a direct bijection between irreducible
triangulations and non-separable maps (themselves in bijection with
quadrangulations without multiple edges) was given by Fusy \FusyNB.

\subsec{Slices of controlled boundary lengths}
Eqs.~\Fnquadi\ and \Fntrii\ are particular examples of the general formula \PMCF
\eqn\Fngenp{F_n = Z_{i-1,i-1}^+(n;\{R_{m+1},S_m\}_{m\geq 0}) - \sum_{k \geq 1} Z_{i-1,i-1+k}^+(n;\{R_{m+1},S_m\}_{m\geq 0})) V_{k;i-2},}
for $i\geq 1$, where $V_{k;p}$ is defined by 
\eqn\Vkpdef{V_{k;p} = \sum_{j \geq k+2} x_j Z_{p+k+1;p}(j-1;\{R_{m+1},S_m\}_{m\geq 0})), \qquad k,p\geq -1.}
Here $Z_{p,p'}(n;\{R_{m+1},S_m\}_{m\geq 0})$ (with $p,p'\geq -1$) denotes the generating function of three-step paths of length $n$, 
starting at height $p$, ending at height $p'$, where each level-step at height $m$ receives a weight $S_m$ while 
each down-step from height $m+1$ to height $m$ receives a weight $R_{m+1}$. The quantity $Z^+_{p,p'}$ (with $p'\geq p\geq -1$) 
denotes the same generating function limited to paths whose heights remain larger than or equal to $p$.
It will always be assumed that $R_0=0$ so in practice the paths which dip below $0$ do not contribute and $Z_{p,p'}$ and  $Z^+_{p,p'}$
depend only on $R_{m+1}$ and $S_m$ for $m\geq 0$. This assumption  also implies that $V_{k,-1}=0$. 

In the above expressions, $R_m$ and $S_m$ must be taken as the solutions of
\eqn\eqsforRpSp{R_{m+1}=1+V_{0,m}, \qquad S_{m}=V_{-1,m},\qquad m\geq 0 .}
The interpretation of $R_{m+1}$ (respectively $S_m$) is that it is the generating functions of slices of type $m'/m'+1$ (respectively of type $m'/m'$)
with $0\leq m'\leq m$.
By an argument similar to that of Sect.~4.1, the reader will be convinced that $V_{k;p}$ for $k\geq -1$ and $p\geq 0$ may then be understood as the generating
function of $k$-slices of type $p'/p'+k+1$ with $0\leq p'\leq p$ and that \Vkpdef\ simply translates the decomposition of 
Fig.~\slicedecomp.

Applying the general substitution relation \FndFnexpr, we get
\eqn\Fndgenp{F_n^{(d)} = Z_{i-1,i-1}^+(n;\{R_{m+1}^{(d)},S_m^{(d)}\}_{m\geq 0}) - \sum_{k \geq 1} Z_{i-1,i-1+k}^+(n;\{R_{m+1}^{(d)},S_m^{(d)}\}_{m\geq 0})
V_{k,i-2}^{(d)}}
for $i\geq 1$, where $R_m^{(d)}$, $S_m^{(d)}$ and $V_{k,p}^{(d)}$ are the series obtained by
substituting, for all $j$ between $1$ and $d$, the formal variable
$x_j$ by the series $X_j^{(d)}$ in $R_m$, $S_m$ and $V_{k,p}$.
As before, the $V_{k,p}^{(d)}$ for $k\geq d-1$ are easily determined in terms of $R_m^{(d)}$, $S_m^{(d)}$ via \Vkpdef, the $x_j$ for $j\geq d+1$ being kept un-substituted. In particular, in the specialized case
of irreducible $d$-angulations, we have $v_{k,p}^{(d)}=0$ for $k\geq d-1$.
As for the $V_{k,p}^{(d)}$ for $-1\leq k\leq d-2$, a first way to
determining them is to use again the conditions \Fsmalln: \Fndgenp\
yields a linear system for the $V_{k,p}^{(d)}$ with $1\leq k\leq d$,
which we may solve. To that end, we note that the inverse of the
semi-infinite unitriangular matrix
$(Z_{i-1,i-1+k}^+(n;\{r_{m+1},s_m\}_{m\geq 0}))_{n,k \geq 0}$ admits
an explicit expression in terms of monomer-dimers \Viennot, which reduces
to \Qdef\ when $r_m=r$ and $s_m=s$ for all $m$. This means that we may
repeat the same strategy as that of Sect.~\elimgen, by simply adding
a new parameter in the game, and derive a closed system of nonlinear
recurrence equations for $R_m^{(d)}$ and $S_m^{(d)}$, generalizing
\eqsforRSd. However, we find the combinatorial meaning of this
approach to be quite unclear.

A more transparent approach is to understand $V_{k,p}^{(d)}$ as the
generating functions of $d$-irreducible $k$-slices of type $p'/p'+k+1$
with $0\leq p'\leq p$ (recall that the substitution does not modify
the length of the boundaries of a slice).
As such, they satisfy new recursive equations which provide an alternative route to determine them. Indeed, the binary decomposition of Sect.~\bindecompproc\ yields immediately
\eqn\recurVkp{V_{k,p}^{(d)}=z \delta_{k,d-2}+(1-\delta_{p,0})\sum_{m=1}^{k+1}V_{m,p-1}^{(d)} V_{k-m,p+m}^{(d)}+
V_{k+2,p-1}^{(d)}\ , \quad -1\leq k\leq d-2 ,\ p\geq 0,}
which is a refined version of \recurV.
For instance, in the case of irreducible triangulations ($d=3$), these equations read
\eqn\vkptri{\eqalign{& v^{(3)}_{1,p} = z + v^{(3)}_{1,p-1}v^{(3)}_{0,p+1}\cr
& v^{(3)}_{0,p} =  v^{(3)}_{1,p-1}v^{(3)}_{-1,p+1}\cr
& v^{(3)}_{-1,p} =  v^{(3)}_{1,p-1}, \qquad \qquad p\geq 0,\cr
}}
with the convention $v^{(3)}_{1,-1}=0$. These equations reproduce precisely \algebRStrii\ upon identifying, $v^{(3)}_{0,p}=r^{(3)}_{p+1}-1$ and  $v^{(3)}_{-1,p}=s^{(3)}_{p}$.

Using the slice/tree bijection in the specialized case
of irreducible $d$-angulations, the $v_{k,p}^{(d)}$ for $1\leq k\leq d-2$ may as well be understood as generating functions
of $d$-oriented $k$-trees with {\it depth} at most $p$. We say that the tree has depth $p$ ($p\geq 0$) if the minimal height assigned to a corner in the
closure procedure of Sect.~6.3 is $-p-1$.
As such, they satisfy the recursive equations:
\eqn\recurVkpbis{v_{k,p}^{(d)}=z \delta_{k,d-2}+(1-\delta_{p,0}) \sum_{q\geq 1} \sum_{
m_i\geq 1, i=1,\cdots ,q \atop m_1+\cdots+m_q=k+2 } \prod_{i=1}^{q} v_{m_i,p-1+\sum_{j=1}^{i-1}m_{j}}^{(d)} \qquad 1\leq k\leq d-2,\ p\geq 0}
(recall that the height increases by $m_i$ when going around a subtree whose root vertex has out-degree 
$m_i$), while $v_{k,p}^{(d)}=0$ for $k\geq d-1$.
\fig{An example of naturally embedded ternary (left) and binary (right) tree with corners labeled according to the
rules of Sect.~\altertreeslice\ when $d=3$ and $d=4$ respectively (starting with label $0$ at the first encountered inner corner).
We see the first corner encountered at each vertex receives a label (in red) equal to (minus) its horizontal position.}{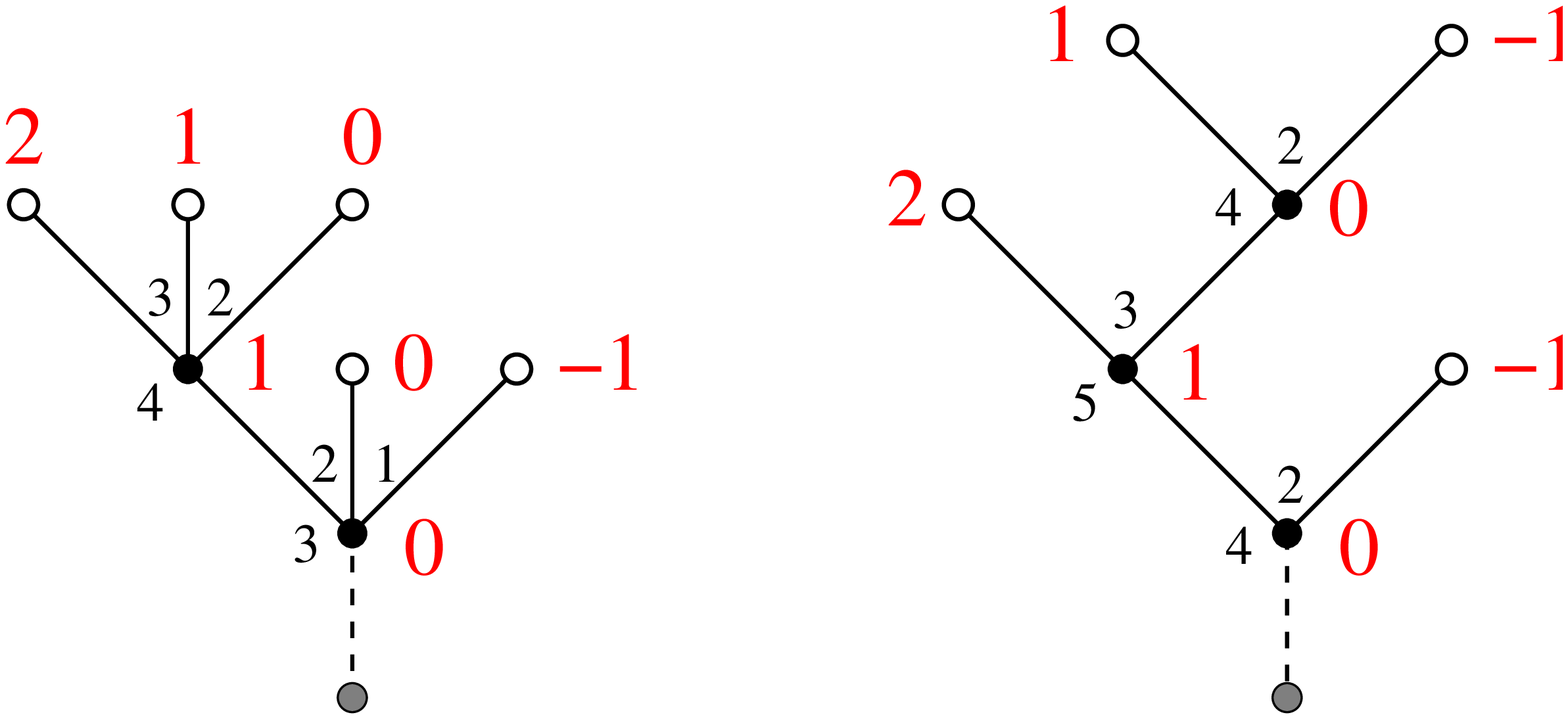}{10.cm}
\figlabel\naturally
For instance, for $d=3$, we get
\eqn\Vkptribis{v^{(3)}_{1,p}=z+(1-\delta_{p,0})\ v^{(3)}_{1,p-1} v^{(3)}_{1,p} v^{(3)}_{1,p+1},\qquad p\geq 0,}
which matches \algebRStrii\ upon identifying
$v^{(3)}_{1,p}=s^{(3)}_{p+1}$. The combinatorial meaning of this
identification is transparent, upon observing that the corner
labelling of a $3$-oriented, i.e. ternary, tree essentially coincides
with the natural embedding of its vertices, see Fig.~\naturally. In the bipartite case ($d$
even), it is easily seen along the same lines as Sect.~\treebipcase\ 
that $v^{(d)}_{k,p}$ vanishes for $k$ odd. In the case
$d=4$, we find that $v^{(4)}_{2,p}$ satisfies
\eqn\Vkpquadbis{v^{(4)}_{2,p} = z+(1-\delta_{p,0})\ v^{(4)}_{2,p-1} v^{(4)}_{2,p+1},\qquad p\geq 0,}
which matches \eqforTi\ upon identifying $v_{2,p}=z T_{p+1}$ (consistenly
with the relations $r^{(4)}_p-1=v^{(4)}_{0,p-1}=v^{(4)}_{2,p-2}$). Again
the combinatorial meaning of the identification is transparent, see again Fig.~\naturally.

Let us conclude by giving a few more integrable equations as obtained from this framework.
If we consider irreducible hexangulations ($d=6$) we find the system
\eqn\vkphex{\eqalign{& v^{(6)}_{4,p} = z + v^{(6)}_{2,p-1}v^{(6)}_{2,p+2}+ v^{(6)}_{4,p-1}v^{(6)}_{0,p+4}\cr
& v^{(6)}_{2,p} =  v^{(6)}_{2,p-1}v^{(6)}_{0,p+2}+ v^{(6)}_{4,p-1}\cr
& v^{(6)}_{0,p} =  v^{(6)}_{2,p-1}, \qquad \qquad p\geq 0,\cr
}}
with the convention $v^{(6)}_{2,-1}=v^{(6)}_{4,-1}=0$. This triangular system yields an equation 
for $v^{(6)}_{0,p}$ which, upon setting
\eqn\vtorhexp{v^{(6)}_{0,p}=r^{(6)}_{p+1}-1}
reads
\eqn\eqforrhexp{\eqalign{z-r^{(6)}_{p}r^{(6)}_{p+2}r^{(6)}_{p+4}+r^{(6)}_{p}r^{(6)}_{p+2} +r^{(6)}_{p}r^{(6)}_{p+3}+r^{(6)}_{p}r^{(6)}_{p+4}
+r^{(6)}_{p+1}r^{(6)}_{p+3}+r^{(6)}_{p+1}r^{(6)}_{p+4}&\cr +r^{(6)}_{p+2}r^{(6)}_{p+4}
-2(r^{(6)}_{p}+r^{(6)}_{p+1}+r^{(6)}_{p+2}+r^{(6)}_{p+3}+r^{(6)}_{p+4})+5 & =0,\quad p\geq 1,\cr}}
with $r^{(6)}_1=r^{(6)}_2=1$, a refined version of \eqforrhex\ which, upon setting
\eqn\Thex{r^{(6)}_p=1+z T_{p-2}}
may itself be rewritten as
\eqn\eqforThexbis{T_p=1+z (T_{p-2}T_{p+1}+T_{p-1}T_{p+1}+T_{p-1}T_{p+2})-z^2 T_{p-2}T_pT_{p+2},\quad p\geq 1}
with initial condition $T_{-1}=T_0=0$. This equation is integrable and an explicit determinantal formula for $T_p$  may easily be obtained 
from the results of \PMCF.

For $d=8$, a similar calculation leads to the integrable equation
\eqn\eqforTocto{\eqalign{T_p& =1+z (T_{p-3}T_{p+1}+T_{p-2}T_{p+1}+T_{p-1}T_{p+1}+T_{p-2}T_{p+2}+T_{p-1}T_{p+2}+T_{p-1}T_{p+3})
\cr &\  -z^2 (T_{p-2}T_pT_{p+2}
+T_{p-3}T_{p-1}T_{p+2}+T_{p-3}T_pT_{p+2}+T_{p-3}T_pT_{p+3}+T_{p-2}T_pT_{p+3}\cr & \qquad\qquad\qquad \qquad+T_{p-2}T_{p+1}T_{p+3})
+z^3 T_{p-3}T_{p-1}T_{p+1}T_{p+3},\qquad p\geq 1\cr }}
with initial condition $T_{-2}=T_{-1}=T_0=0$ by setting $v^{(8)}_{0,p}=r^{(8)}_{p+1}-1$ and $r^{(8)}_p=1+z T_{p-3}$.

\newsec{Discussion}

Let us end this paper by discussing some extensions of our work. The first one concerns
$d$-irreducible maps with outer degree $d$, which are made nontrivial by slightly relaxing the
notion of irreducibility. The second extension concerns $d$-irreducible maps with two marked faces
of degree strictly larger than $d$. The third extension, suggested by one of the referees, concerns
maps with two marked faces, subject to a control on two distinct girth parameters.
 
\subsec{$d$-irreducible maps with outer degree $d$}
\subseclab\weakirred

In this section, we deal with $d$-irreducible maps with outer degree $d$, and
with at least one inner face (i.e.\ maps not reduced to a tree). In the way
we defined $d$-irreducibility so far, we have a unique such
map, made of a single inner face of degree $d$ glued to the external face:
we shall call it the {\it trivial} map of outer degree $d$. Indeed
the boundary of the external face (called {\it external
boundary} in the following) forms a cycle of length $d$ and, as such, has
to be the boundary of an inner face of degree $d$. A weaker and somewhat
more natural definition of $d$-irreducibility among maps with outer degree
$d$ consists in simply picking those maps of girth $d$ such that all cycles of
length $d$ in the map are the boundary of an arbitrary (i.e.\ inner or external)
face of degree $d$. With this definition, the external boundary needs not
surround a single inner face any longer as it is already the boundary of
the external face of degree $d$.

We shall call $H_d(z;x_{d+1},\ldots)$ the generating function of such
(weakly) $d$-irreducible maps with outer degree $d$ counted, as before,
with a weight $z$ per inner face of degree $d$ and weights $x_j$ per face of
degree $j\geq d+1$. We claim that we have the relation:
\eqn\Hdoddeven{\eqalign{ H_d(z;x_{d+1},\ldots)&=2z-X_d(z;x_{d+1},\ldots), \qquad d\ {\rm odd}\cr
 H_d(z;x_{d+1},\ldots)&=2z+{d\over 2}\times {z^3\over 1+z}-X_d(z;x_{d+1}\ldots), \qquad d\ {\rm even}\cr}}
 where $X_d(z;x_{d+1},\ldots)$ is defined via the inversion \Gcompinv.
In practice, $X_d$ may be obtained for instance from the
knowledge of $V_{d-2}^{(d)}$ via the relation
\eqn\Xdpract{X_d={1\over(R^{(d)})^{d-1}}\left(
V_{d-2}^{(d)}-\sum_{j\geq d+1} x_j P_{-d+1}(j-1;R^{(d)},S^{(d)})\right)}
obtained by taking \Vjdef\ at $k=d-2$ and substituting the renormalized weight
$x_d\to X_d=X_d^{(d)}(z;x_{d+1},\ldots)$.

To prove \Hdoddeven, we use the same substitution approach as in Sect.~3.1.
i.e.\ start from an arbitrary map with outer degree $d$ and girth $d$,
and obtain a $d$-irreducible one by replacing each outermost cycle
of length $d$ by a single $d$-valent face. As we shall now see, this
transformation leads to a functional relation between $H_d$ and the generating
function $G_d(x_d,x_{d+1},\ldots)$ (as defined in Sect.~3) of maps with outer
degree $d$ and girth $d$, counted with weights $x_j$ per inner face of degree
$j\geq d$. It is important to note that, with our weaker notion of
$d$-ireducibility, we now have to take the convention that, in
the determination of the outmost cycles, the external boundary itself should
{\it not} be considered as a cycle of length $d$. This leads to two
caveats: (i) the trivial map of outer degree $d$ (which is $d$-irreducible
and contributes $z$ to $H_d$) has a single antecedent which is the trivial map
of outer degree $d$ itself (contributing $x_d$ to $G_d$) and (ii) some
outermost cycles may now overlap so that some extra prescriptions are required
for a well defined replacement.
Let us analyze this overlapping issue more precisely: when repeating the arguments
of Sect.~3, the only problem that we face is, with the notations of
Fig.~\twocycles, that the cycle $C_{12}$ may coincide with the external
boundary so that it creates no contradiction with
the fact that $C_1$ and $C_2$ were considered as outermost in the first place.
Still, since $C_{12}$ has length $d$, this implies that, necessarily,
$n_1=n_2=d/2$ (recall that $n_1$ and $n_2$ are at least $d/2$) so
that $v$ and $v'$ are antipodal around the external boundary while the two
internal paths connecting them form mutually avoiding {\it diagonal paths} of
length $d/2$. Since no paths of length $d/2$ exist for $d$ odd, we conclude
that outermost paths cannot overlap for $d$ odd and the caveat (ii) never
occurs. For odd $d$, we deduce the relation
\eqn\gdeq{G_d(x_d,x_{d+1},\ldots)-x_d=H_d(G_d(x_d,x_{d+1},\ldots);x_{d+1},\ldots)
-G_d(x_d,x_{d+1},\ldots)\qquad d\ {\rm odd}}
for the enumeration of nontrivial maps with outer degree $d$ and girth $d$
(the r.h.s. is obtained by substituting $z\to G_d$ in $H_d-z$, the
generating function of nontrivial $d$-irreducible maps of outer degree $d$).
Now, from \Gcompinv, we may perform the substitution $G_d\to z$ by setting
$x_d=X_d(z;x_{d+1},\ldots)$ in \gdeq.This yields the announced
result \Hdoddeven\ for odd $d$.
\fig{Schematic picture of a configuration of map with outer degree $d$ and girth $d$
having overlapping outermost cycles of length $d$ (by convention, the external boundary
-- in magenta -- is not considered here as a cycle of length $d$). Two antipodal vertices
$v$ and $v'$ along the external boundary are linked by $p\geq 2$ diagonal paths
of length $d/2$. }{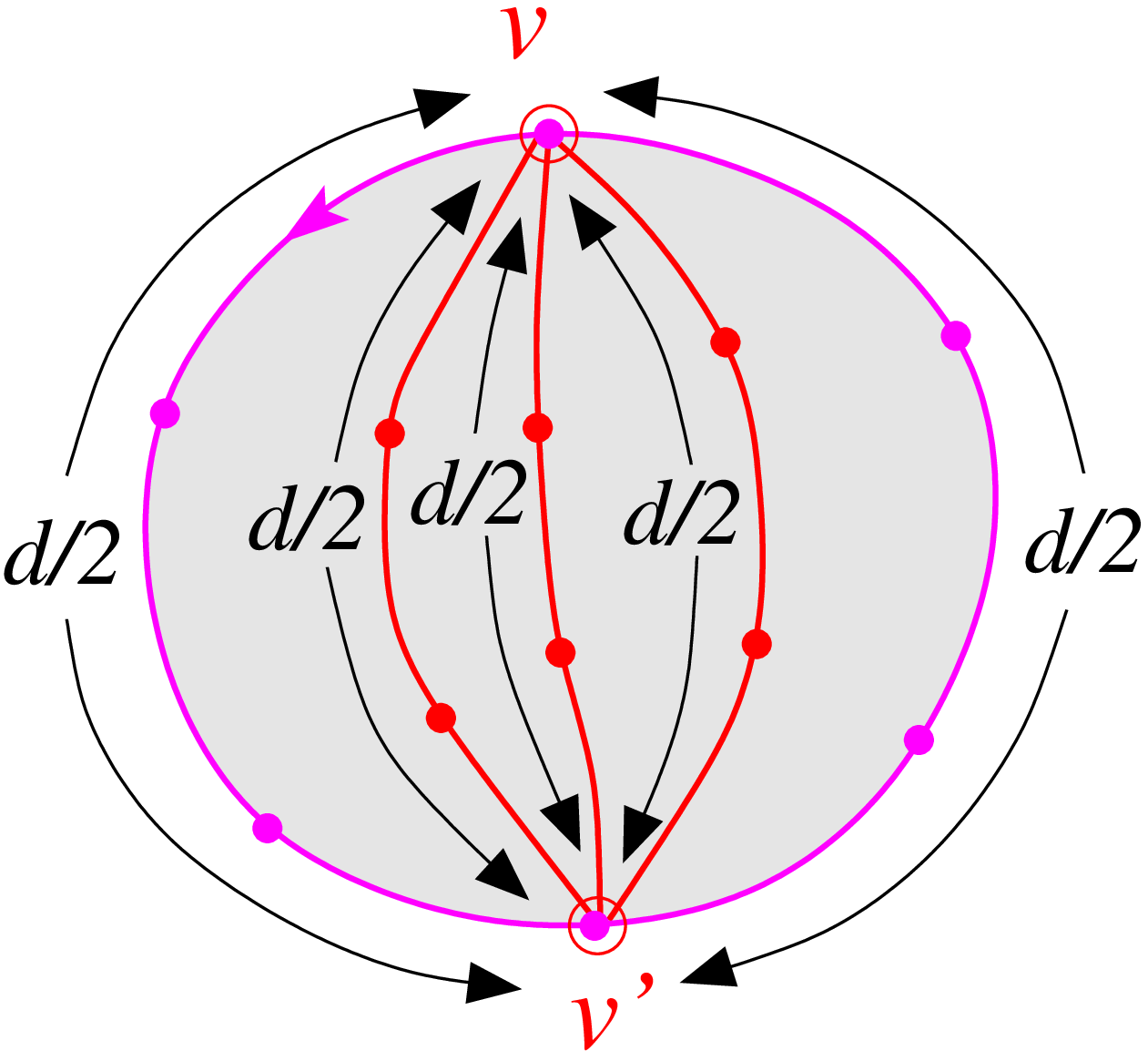}{5.cm}
\figlabel\diagonals

For even $d$, we note that, in a map of girth $d$, there cannot be more than
one pair of antipodal vertices connected by a diagonal path of length $d/2$ as
otherwise, diagonals would have to cross and a cycle of length $<d$ would be
created. Overlapping therefore appears only in situations where exactly one
pair of antipodal vertices $v$ and $v'$ are connected
by an arbitrary number $p\geq 2$ of mutually avoiding diagonal paths of
length $d/2$ (see Fig.~\diagonals). We therefore decide to treat separately and remove from our original set of
maps with girth $d$ the configurations having two antipodal vertices
connected by at least one diagonal path of length $d/2$ (we found simpler
to also include in the removed set situations having a single diagonal although
they do not give rise to overlapping). The set of these removed configurations
clearly displays a $d/2$-fold symmetry by rotation around the external face.
Once these configurations are removed, the $d/2$ $d$-irreducible maps with outer
degree $d$ made of two inner faces of degree $d$ sharing a diagonal of
length $d/2$ cannot be recovered any longer.
Each of these configurations contributes $z^2$ to $H_d$, so we deduce
\eqn\gdeqbis{\eqalign{G_d(x_d &,x_{d+1},\ldots) -x_d
-{d\over 2} G_{d,D}(x_d,x_{d+1},\ldots)\cr &
=H_d(G_d(x_d,x_{d+1},\ldots);x_{d+1},\ldots)-G_d(x_d,x_{d+1},\ldots)-{d\over 2}
(G_d(x_d,x_{d+1},\ldots))^2\cr & \hskip 10.cm d\ {\rm even}\cr}}
for the enumeration of nontrivial maps with outer degree $d$ and girth $d$
without a diagonal of length $d/2$ (the r.h.s. is obtained by substituting 
$z\to G_d$ in $H_d-z-{d\over 2}z^2$). In the l.h.s., $G_{d,D}$ enumerates maps
with outer degree $d$ and girth $d$ with $p\geq 1$ diagonal paths
of length $d/2$ connecting, say the origin of the root edge to its antipodal
vertex. Cutting along all the diagonal paths of length $d/2$, any such map forms
a sequence of $p+1\geq 2$ rooted maps (see Fig.~\diagonals\ for an example with $p=3$). 
These latter maps cannot have any more
diagonal paths of length $d/2$ connecting the origin of the root edge to its antipodal
vertex but are otherwise arbitrary maps of outer degree $d$ and girth $d$
(in particular they may have diagonal paths of length $d/2$ connecting other
pairs of antipodal vertices). Their generating function is therefore
$G_d-G_{d,D}$ and we have
\eqn\gdD{G_{d,D}={(G_d-G_{d,D})^2\over 1-(G_d-G_{d,D})},}
or equivalently
\eqn\gdDbis{G_{d,D}={(G_d)^2\over 1+G_d},}
so that \gdeqbis\ may eventually be written as
\eqn\gdeqter{G_d=x_d-{d\over 2}\times {G_d^3\over 1+G_d}+H_d(G_d;x_{d+1},\dots)
-G_d\qquad d\ {\rm even}\ .}
Setting $x_d=X_d(z;x_{d+1},\ldots)$ in \gdeqter\ yields the announced
result \Hdoddeven\ for even $d$.

Let us end this section by listing expressions for the generating function
$h_d(z)=H_d(z;0,\ldots)$ of (weakly) irreducible $d$-angular dissections
of the $d$-gon at $d=3,4,5,6$, as obtained by a specialization at $x_j=0$,
$j\geq d+1$. Eq.~\Xdpract\ reduces to the simpler relation
\eqn\xdpract{X_d(z;0,\ldots)={v_{d-2}^{(d)}\over(r^{(d)})^{d-1}}\ .}
For $d=3$, using $v_1^{(3)}=s^{(3)}$, we obtain
\eqn\htrois{\eqalign{h_3(z)&=2z-{s^{(3)}\over (r^{(3)})^2}\cr
&=
z+z^3+z^7+3 z^9+12 z^{11}+52 z^{13}+241 z^{15}+1173 z^{17}+5929 z^{19}+30880
   z^{21}\cr & \qquad +164796 z^{23}+897380 z^{25}+4970296 z^{27}+27930828
   z^{29}+O\left(z^{31}\right)
\cr}}
in agreement with \TutteCPT.

\noindent For $d=4$, using $v_2^{(4)}=u_1^{(4)}=r^{(4)}-1$, we obtain
\eqn\hquatre{\eqalign{h_4(z)& =2z+2 {z^3\over 1+z}-{r^{(4)}-1\over (r^{(4)})^3}\cr
&=
z+2 z^2+z^5+4 z^7+6 z^8+24 z^9+66 z^{10}+214 z^{11}+676
   z^{12}+2209 z^{13}+7296 z^{14}
   \cr &\qquad +24460 z^{15}+82926
   z^{16}+284068 z^{17}+981882 z^{18}+3421318
   z^{19}+O\left(z^{20}\right)
   \cr}}
in agreement with \MuSch.

\noindent For $d=5$, using $v_3^{(5)}=3s^{(5)}+(s^{(5)})^3-2s^{(5)}r^{(5)}$, we obtain
\eqn\hcing{\eqalign{h_5(z)&=2z-{3s^{(5)}+(s^{(5)})^3-2s^{(5)}r^{(5)}\over (r^{(5)})^4}\cr
&=
z+5 z^3+46 z^5+1350 z^7+52360 z^9+2382508 z^{11}+119914425 z^{13}\cr & \qquad +6470326059
   z^{15}+367369835490 z^{17}+21686295649075 z^{19}+O\left(z^{21}\right).
\cr}}
For $d=6$, using $v_4^{(6)}=u_2^{(6)}=-2+3 r^{(6)}-(r^{(6)})^2$, we obtain
\eqn\hsix{\eqalign{h_6(z)& =2z+3 {z^3\over 1+z}-{-2+3 r^{(6)}-(r^{(6)})^2\over (r^{(6)})^5}\cr
&=
z+3 z^2+2 z^3+5 z^4+42 z^5+266 z^6+1986 z^7+15552 z^8+127738 z^9+1086998
   z^{10}\cr
   &\qquad +9517362 z^{11}+85291440 z^{12}+779292490 z^{13}+7237661226
   z^{14}+O\left(z^{15}\right).
   \cr}}

\subsec{$d$-irreducible maps with two marked faces of degree strictly larger than $d$}

This section deals with bipartite maps for simplicity. So far we considered maps 
with a single marked face (the external face) of degree, say $2m$ and a marked oriented edge (the root edge)
incident to the external face and oriented so that the external face lies on its right
(in practice we marked the root edge first). We may instead consider maps with 
two marked distinct (and distinguished) faces of respective
degrees $2m$ and $2m'$, and a marked oriented edge incident to each of these
marked face (and having the marked face on its right).  Again we may demand that these maps
be $d$-irreducible, i.e.\ have girth at least $d$ and be such that all cycles of length $d$ 
are the boundary of an inner face (i.e.\ a face different from the marked faces) of degree $d$. 
Assuming $d=2b$ and $m,m'>b$, we may
simply obtain the generating function $F_{2m,2m'}^{(d)}(z;x_{d+1},\ldots)$ of $d$-irreducible
maps with two marked faces of respective degree $2m$ and $2m'$ via
\eqn\Fmmprime{F^{(d)}_{2m,2m'}= 2m' {\partial F_{2m}^{(d)} \over \partial x_{2m'}}}
since the desired maps are obtained from $d$-irreducible maps with a single 
marked face of degree $2m$ (as enumerated by $F^{(d)}_{2m}$) by marking a face of 
degree $2m'$ (via the action of $\partial \cdot/\partial x_{2m'}$) and then marking 
an incident oriented edge (among $2m'$ choices).
Now we may use \Fndbipfin\ to compute 
\eqn\dFmmprime{{\partial F_{2m}^{(d)}  \over \partial x_{2m'}}= {2m \choose m-b} (R^{(d)})^{m-b}
 \left(K(m')- {b+m' \over m+m'} {2m'-1 \choose m'+b}  (R^{(d)})^{b+m'} \right)}
where the quantity 
\eqn\valK{\eqalign{K(m')&= \Bigg(\sum_{\ell=0}^{b-1}(-1)^{b-\ell-1} (b-\ell) {b+\ell \choose 2\ell} {\rm Cat}(\ell) (R^{(d)})^{b-\ell-1}
\cr & \qquad - \sum_{k \geq b+1} (b+k) {2k-1 \choose k+b} x_{2k} (R^{(d)})^{b+k-1} \Bigg)
  {\partial R^{(d)}  \over \partial x_{2m'}} \cr}}
is independent of $m$. Eq.~\dFmmprime\ is valid also for $m=b$, which allows to determine 
$K(m')$ from the identity $\partial F_{d}^{(d)} /\partial x_{2m'}=0$ (since $F_{d}^{(d)}=z+{\rm Cat}(d/2)$),
namely
\eqn\valKbis{K(m')=  {2m'-1 \choose m'+b}  (R^{(d)})^{b+m'}\ .}
Plugging this value in \dFmmprime, we arrive at 
\eqn\Fmmprimebis{F^{(d)}_{2m,2m'}= (m-b)(m'-b) {2m \choose m-b}{2m'\choose m'-b}{(R^{(d)})^{m+m'}
\over m+m'}\qquad m,m'>b .}
This formula may be viewed as a generalization of a similar formula \CF\ for general (non necessarily irreducible)
maps with two marked faces, corresponding to $d=b=0$ here.
The more general formulas of \CF\ for maps with more than two marked faces can also be 
extended to the case of $d$-irreducible maps \WIP.

\subsec{$(d,d')$-irreducible annular maps}

We are indebted to the anonymous referee for pointing out this
extension of our work. Maps with two marked faces are naturally
endowed with two distinct girth parameters: the {\it separating girth}
and the {\it non separating girth}, defined respectively as the
minimum length of cycles separating and not separating the marked
faces from one another. Maps subject to an independant control on both
girth parameters were enumerated in \BFb, and a natural question is
whether these results can be extended to the context of irreducible
maps.

We are led to slightly generalize the notion of annular map,
introduced in Section~\pointing, as follows. For $d,d'$ nonnegative
integers, we define a {\it $(d,d')$-quasi-irreducible annular map} as a
rooted map with a marked inner face of degree $d'$ called the {\it
central face}, subject to the following constraints:
\item{-} any non separating cycle (i.e.~a cycle which does not
separate the outer face and the central face) has length at least $d$,
and if its length is equal to $d$ then it is necessarily the boundary
of an inner face of degree $d$,
\item{-} any separating cycle (i.e.~a cycle which separates the outer
face and the central face) has length at least $d'$.
\par \noindent
If furthermore the only separating cycle of length $d'$ is the
boundary of the central face, then the map is called a {\it
$(d,d')$-irreducible annular map}. In other words, a
$(d,d')$-irreducible annular map has non separating girth $d$,
separating girth $d'$, and all its minimal non separating or separating
cycles are ``trivial'' (in the case of $(d,d')$-quasi-irreducible
annular map, the triviality condition for minimal separating cycles is
relaxed).
Let us denote by $I_n^{(d,d')} \equiv
I_n^{(d,d')}(z;x_{d+1},x_{d+2},\ldots)$ [resp.~${\tilde I}_n^{(d,d')}
\equiv {\tilde I}_n^{(d,d')}(z;x_{d+1},x_{d+2},\ldots)$] the
generating function of $(d,d')$-irreducible
[resp.~$(d,d')$-quasi-irreducible] annular maps whose root face has
degree $n$, counted with a weight $z$ per inner non central face of
degree $d$ and, for all $i \geq d+1$, a weight $x_i$ per inner non
central $i$-valent face (by convention, the central face receives no
weight). Note that, for $d=d'$, $(d,d')$-irreducible annular maps are
precisely the annular maps considered in Section~\pointing, and thus
$I_n^{(d,d)}= {\partial F_n^{(d)} \over \partial z} =
P_d(n;R^{(d)},S^{(d)})$ by virtue of \Fndpointgen. Note also that, for
$z=0$, $d' I_n^{(d,d')}(0;x_{d+1},x_{d+2},\ldots)$ [resp.~$d' {\tilde
I}_n^{(d,d')}(0;x_{d+1},x_{d+2},\ldots)$] coincides with the quantity
denoted ${\widehat G}_{d+1}^{(d',n)}$ [resp.~$G_{d+1,d'}^{(d',n)}$] in
\BFb\ (where the annular maps are doubly rooted, thus the extra factor
$d'$).

We then have the beautiful identities
\eqn\annul{I_n^{(d,d')}=P_{d'}(n;R^{(d)},S^{(d)})}
\eqn\annulbis{{\tilde I}_n^{(d,d')}=P_{-d'}(n;R^{(d)},S^{(d)})}
where we recall that $P_k(n;\cdot,\cdot)$ is a three-step path
generating function, given by \pathgf. These identities generalize
results from [\xref\BFb, Section 6.2], in particular Equation (7) of
that paper coincides up to notations with our second identity at
$z=0$.  We may establish them bijectively along the same lines as in
the proof of the pointing formula \Fndpointgen\ in
Section~\pointing. Let us simply mention which modifications have to
be made, and leave the reader check the details.

The first identity \annul\ is obtained rather straightforwardly.
Starting from a slice-decorated three-step path counted by
$P_{d'}(n;R^{(d)},S^{(d)})$, we define the partial gluing as in
Section~\pathtoannul, but the gluing of the right and left boundaries
now leaves $d'$ unmatched edges, which form the central face. The
slice decomposition of annular maps is still performed by passing to
the lift and cutting along leftmost shortest paths to $-\infty$.

For the second identity \annulbis, given a slice-decorated three-step
path counted by $P_{-d'}(n;R^{(d)},S^{(d)})$, its partial gluing now
has a left boundary longer by $d'$ edges than the right boundary, and
thus the central face is formed by ``wrapping'' on the other side. The
fact that the left boundary is a non necessarily unique shortest path
between its endpoints explains why the resulting annular map is
$(d,d')$-quasi-irreducible but not necessarily $(d,d')$-irreducible.
Conversely, the slice decomposition must be adapted: informally
speaking, we have to cut the lift along leftmost shortest paths to
$+\infty$ instead of $-\infty$. Because there might be several minimal
separating cycles in the annular map, it is no longer true that, in
the lift, all shortest paths to $+\infty$ coalesce with the upper
boundary, still the leftmost ones do.

As a corollary of \annul, \annulbis\ and the path definition of
$P_k(n;\cdot,\cdot)$, we have the relation
\eqn\annulrel{{\tilde I}_n^{(d,d')} = I_n^{(d,d')} \left(R^{(d)}\right)^{d'}}
which may alternatively obtained by decomposing a
$(d,d')$-quasi-annular map of outer degree $n$ along its outermost
minimal separating cycle. Using a similar decomposition, we may easily
obtain an extension of Theorem 32 in \BFb, i.e.\ an expression for the
generating function of $(d,d')$-quasi-irreducible annular maps where
we no longer impose that the central face has degree $d'$ (in other
words, the theorem is still valid with the ``extra variable $z$'').

Finally, let us observe that ${\tilde I}_{d}^{(d,d)} =
\left(R^{(d)}\right)^{d}$ is closely related to the generating
function $H_d$ discussed in Section~\weakirred. Indeed, a
$(d,d)$-quasi-irreducible annular map with outer degree $d$ is ``almost''
a weakly $d$-irreducible map with a marked inner face of degree $d$,
except for the fact that there might exist nontrivial
minimal separating cycles.
More precisely, we have the relation
\eqn\IHrel{\eqalign{ \left(R^{(d)}\right)^{d} &= {1 \over 1 - \left(
  {\partial H_d \over \partial z} - 1\right)}, \qquad d\ {\rm odd}\cr
 \left(R^{(d)}\right)^{d}&={1 \over 1 - \left(
  {\partial H_d \over \partial z} - 1 - d z\right)
  - {d\over 2} {2z+z^2\over 1+2z+z^2}}, \qquad d\ {\rm even}\cr}}
which may be obtained by differentiating \Hdoddeven\ with respect to
$z$, and noting that, by \Gcompinv, ${\partial X_d \over \partial
z}=\left(R^{(d)}\right)^{-d}$ since ${ \partial G_d \over \partial
x_d} = \tilde{I}_d^{(d-1,d)}(0;x_d,x_{d+1},\ldots) =
\left(R^{(d-1)}(0;x_d,x_{d+1},\ldots)\right)^d$ and
$R^{(d-1)}(0;X_d(z,x_{d+1},\ldots),x_{d+1},\ldots)=R^{(d)}(z;x_{d+1},\ldots)$.
Alternatively, a combinatorial proof is obtained by decomposing a
$(d,d)$-quasi-irreducible annular map along its minimal separating
cycles.  This decomposition works straightforwardly for $d$ odd, but
for $d$ even there are some pathologies related to the possible
existence of diagonal paths, already observed in Section~\weakirred.
Working out the details is left as a pleasant exercise to the reader.

\bigskip \noindent {\bf Acknowledgements:} We are grateful to the
referees for their numerous suggestions of improvements. The work of
JB was partly supported by the ANR projects ``Cartaplus''
12-JS02-001-01 and ``IComb'' ANR-08-JCJC-0011.  
\listrefs

\end